\documentclass[11pt, oneside]{article}

\usepackage[dvips]{graphicx}
\usepackage{epsfig}
\usepackage{tikz}
\usepackage{amsmath}
\usepackage{amssymb}
\usepackage{xcolor}
\usepackage{mathtools}

\usepackage{booktabs}
\usepackage{multicol}
\usepackage{multirow}
\usepackage[labelfont=bf]{caption}
\usepackage{subcaption}

\newtheorem{exmp}{Example}

\numberwithin{equation}{section}
\numberwithin{figure}{section}

\begin{document}
\baselineskip=2pc

\vspace{.5in}

\begin{center}

{\large\bf
High-order exponential time differencing multi-resolution alternative finite difference WENO methods for nonlinear degenerate parabolic equations
} 

\end{center}

\vspace{.1in}

\centerline{Ziyao Xu\footnote{Department of Applied and Computational Mathematics and Statistics,
University of Notre Dame, Notre Dame, IN 46556, USA. E-mail: zxu25@nd.edu
} and
Yong-Tao Zhang\footnote{Department of Applied and Computational Mathematics and Statistics,
University of Notre Dame, Notre Dame, IN 46556, USA. E-mail: yzhang10@nd.edu
}
}

\vspace{.3in}

\abstract{
In this paper, we focus on the finite difference approximation of nonlinear degenerate parabolic equations, a special class of parabolic equations where the viscous term vanishes in certain regions. 
This vanishing gives rise to additional challenges in capturing sharp fronts, beyond the restrictive CFL conditions commonly encountered with explicit time discretization in parabolic equations. 
To resolve the sharp front, we adopt the high-order multi-resolution alternative finite difference WENO (A-WENO) methods for the spatial discretization, which is designed to effectively suppress oscillations in the presence of large gradients and achieve nonlinear stability. 
To alleviate the time step restriction from the nonlinear stiff diffusion terms, we employ the exponential time differencing Runge-Kutta (ETD-RK) methods, a class of efficient and accurate exponential integrators, for the time discretization. 
However, for highly nonlinear spatial discretizations such as high-order WENO schemes, it is a challenging problem how to efficiently form the linear stiff part in applying the exponential integrators, since direct computation of a Jacobian matrix for high-order WENO discretizations of the nonlinear diffusion terms is very complicated and expensive. 
Here we propose a novel and effective approach   of replacing the exact Jacobian of high-order multi-resolution A-WENO scheme with that of the corresponding high-order linear scheme in the ETD-RK time marching, 
based on the fact that in smooth regions the nonlinear weights closely approximate the optimal linear weights, while in non-smooth regions the stiff diffusion degenerates.
The algorithm is described in detail, and numerous numerical experiments are conducted to demonstrate the effectiveness of such a treatment and the good performance of our method. The stiffness of the nonlinear parabolic  partial differential equations (PDEs) is resolved well, and large time-step size computations of $\Delta t \sim O (\Delta x)$ are achieved.  }

\vspace{.3in}


\pagenumbering{arabic}

\section{Introduction}
In this paper, we are concerned with the finite difference methods for nonlinear degenerate parabolic equations, which arise in many science and engineering applications. The equations are formulated as follows:
\begin{equation}\label{eq:PDE1D}
u_t + f(u)_x = g(u)_{xx} := (a(u)u_x)_x,
\end{equation}
in one space dimension, where $f(u)$ is the flux of the hyperbolic term, and $g(u)$ is the function of the degenerate parabolic term.
The diffusion coefficient $a(u)=g'(u)$ vanishes for certain regions of $u$.
Similarly, in two space dimensions, the nonlinear degenerate parabolic equations are formulated as:
\begin{equation}\label{eq:PDE2D}
u_t+f_1(u)_x+f_2(u)_y=g_1(u)_{xx}+g_2(u)_{yy},
\end{equation}
where the functions of the parabolic terms $g_1(u), g_2(u)$ may degenerate, i.e., $g_1'(u)$ and $g_2'(u)$ vanish for certain regions of $u$.
Due to the vanishing of parabolic terms, the equations \eqref{eq:PDE1D} and \eqref{eq:PDE2D} exhibit some features of hyperbolic equations, e.g., the existence of non-smooth weak solutions and the propagation of sharp wave fronts with finite speed.

\begin{figure}[!htbp]
 \centering
 \begin{subfigure}[b]{0.4\textwidth}
  \includegraphics[width=\textwidth]{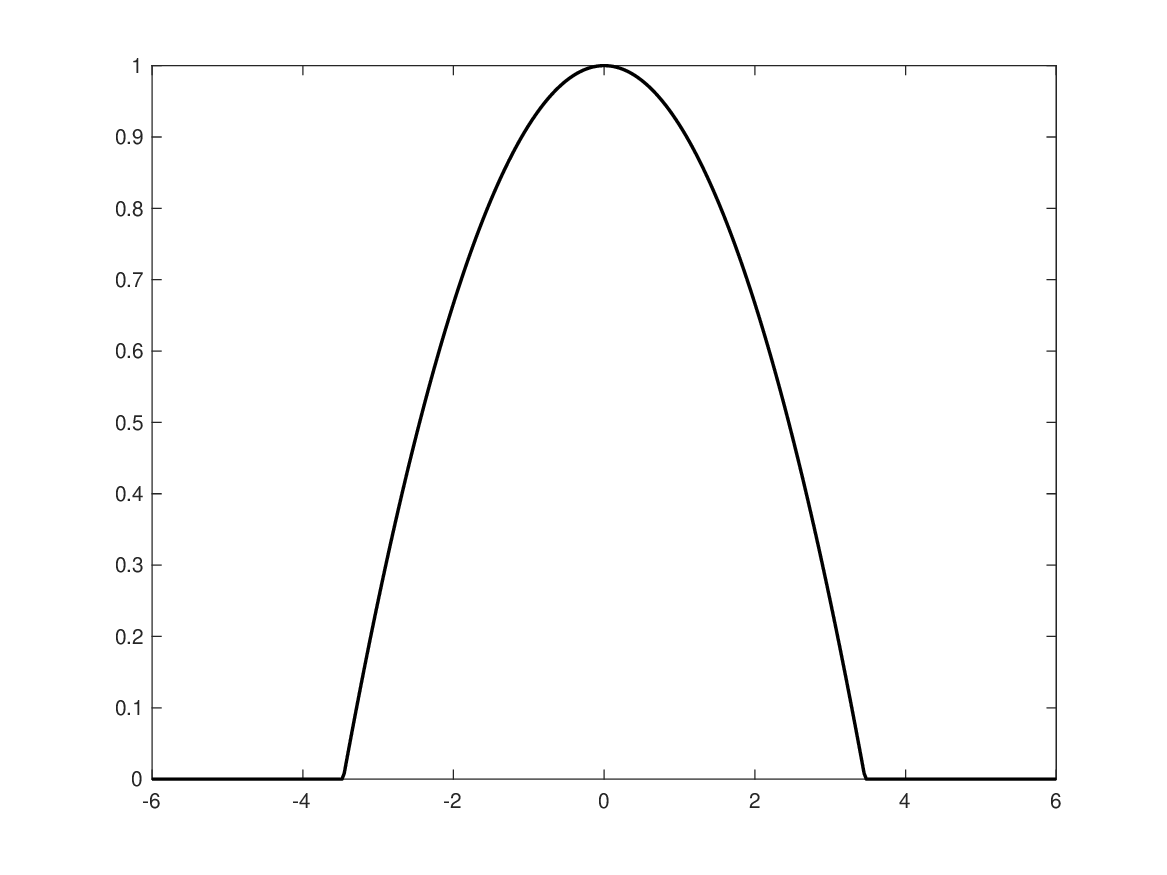}
  \caption{$m=2$}
 \end{subfigure}
 \begin{subfigure}[b]{0.4\textwidth}
  \includegraphics[width=\textwidth]{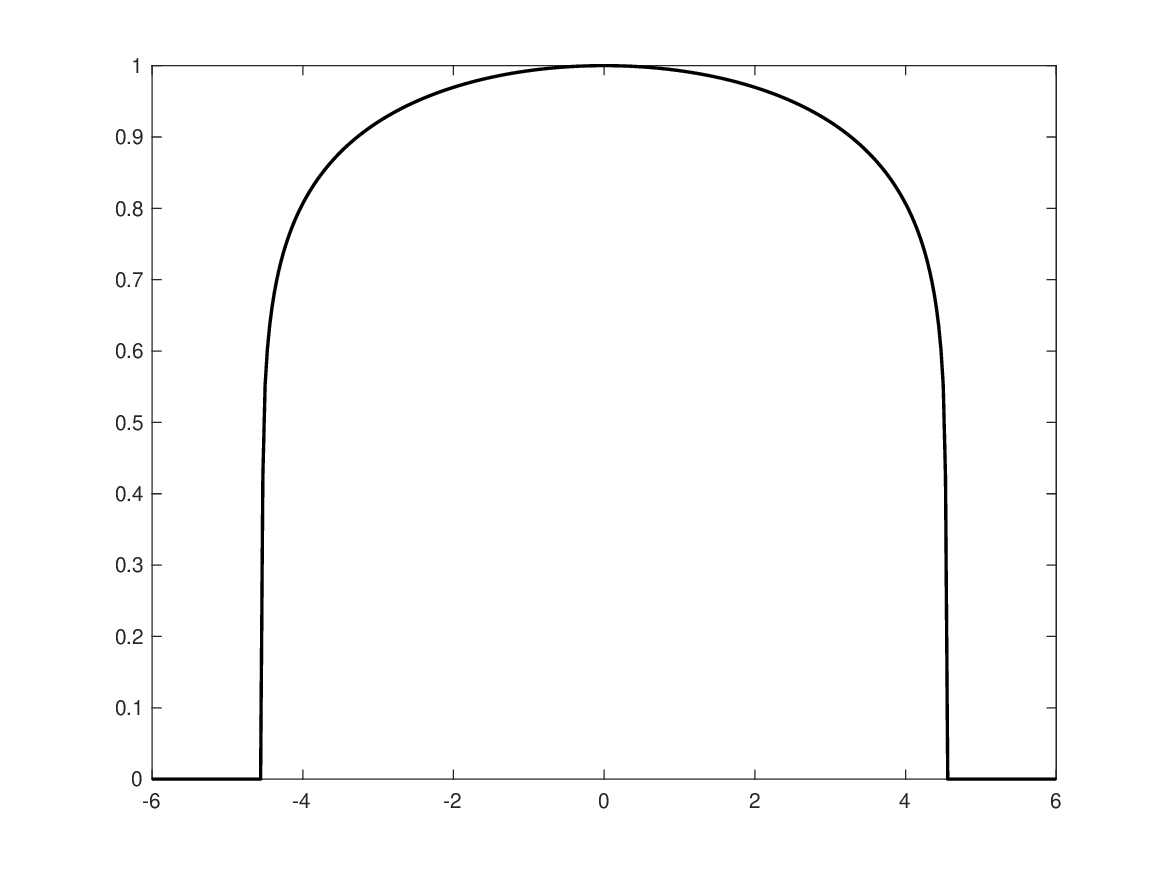}
  \caption{$m=8$}
 \end{subfigure}
 \caption{\footnotesize{The Barenblatt solutions \eqref{eq:Barenblatt} of PME \eqref{eq:PME1D} at $t=1$.}}
 \label{fig:intro}
\end{figure}

A prototypical example of nonlinear degenerate parabolic equations is the porous media equation (PME) \cite{Muskat,Aronson}:
\begin{equation}\label{eq:PME1D}
u_t=(u^m)_{xx},
\end{equation}
used to model gas flow in porous media, where $m>1$ is a constant, and $u\geq0$ denotes the fluid density. 
The PME admits the well-known Barenblatt weak solution:
\begin{equation}\label{eq:Barenblatt}
B_m(x,t)=t^{-p}\left(\left(1-\frac{p(m-1)}{2m}\frac{|x|^2}{t^{2p}}\right)^+\right)^{\frac{1}{m-1}},
\end{equation}
where $z^+:=\max(z,0)$ and $p=\frac{1}{m+1}$. 
The wave fronts $x=\pm t^{p}\sqrt{\frac{2m}{p(m-1)}}$ in the Barenblatt solution propagate at finite speed, as illustrated by the solutions at $t=1$ in Figure \ref{fig:intro}.

Due to the existence of sharp fronts caused by the hyperbolic features of the equations \eqref{eq:PDE1D} and \eqref{eq:PDE2D}, high-order linear schemes suffer from spurious oscillations in the presence of large gradients, i.e. the Gibbs phenomenon, even though they work well in smooth regions. Various numerical methods with nonlinear components have been developed for solving nonlinear degenerate parabolic PDEs, for example, linear approximation scheme based on the nonlinear
Chernoff formula with a relaxtion parameter \cite{Mage,Noch1}, local discontinuous Galerkin
method \cite{ZhWu}, finite volume method \cite{Besse}, kinetic schemes \cite{Aregb}, high-order relaxation schemes \cite{Cavalli}, kernel based integral method \cite{Christ}, moving mesh finite element method \cite{Ngo}, etc.

In the literature (e.g. \cite{liu2011high, Abe1, arbogast2019finite, jiang2021high, Abedian, zhang2022high, Qiu}), it has been shown that the weighted essentially non-oscillatory (WENO) methods, a popular class of high-order accuracy schemes for solving hyperbolic PDEs, are very effective to be adopted in the spatial discretization of degenerate parabolic equations, to achieve high-order accuracy in smooth regions while suppressing oscillations near wave fronts. 
The WENO methods are based on the successful essentially non-oscillatory (ENO) methods with additional advantages.
The ENO methods were first proposed by Harten et al. \cite{ENO0} in the finite volume framework for hyperbolic equations, adopting the smoothest stencil among several candidates to reconstruct the solutions.
To enhance efficiency, Shu and Osher proposed the finite difference framework in the subsequent work \cite{ShuOsher1, ShuOsher2}, which allow for computation in a dimension-splitting fashion. 
The WENO methods use weighted stencils rather than choosing only one stencil, to reconstruct solutions. 
The weighting strategy is a crucial component of WENO methods, designed based on the principle that in smooth regions of the solution, the nonlinear weights are close to the linear weights in linear schemes to yield improvement in accuracy, while being close to zero to minimize the contribution of stencils with large gradient of the solution.
Various types of WENO methods have been proposed since the celebrated WENO methods developed in \cite{WENO0,WENOJS}, e.g.,  very high-order finite difference WENO schemes \cite{Bal}, the mapped WENO schemes \cite{henrick2005mapped} and WENO-Z schemes \cite{borges2008improved, castro2011high} for optimal order near critical points, central WENO schemes \cite{Levy}, energy stable WENO schemes \cite{Yama}, robust WENO schemes \cite{LiuZ}, among others.

Recently in \cite{zhu2018new}, the multi-resolution WENO schemes were developed to solve hyperbolic conservation laws. This class of WENO schemes construct WENO approximations on unequal-sized substencils and exhibit interesting properties. For example, it is flexible to construct linear weights for the multi-resolution WENO schemes, which in general can be
taken as arbitrary positive numbers with the only requirement that their sum equals 1. Such flexibility simplifies the procedure in applications of WENO schemes, e.g., sparse-grid simulations for multidimensional problems \cite{tsybulnik2022efficient}, treatment of negative linear weights for solving degenerate parabolic equations \cite{jiang2021high}, etc.      
In this paper, similar to \cite{jiang2021high}, we adopt the multi-resolution WENO method \cite{zhu2018new} in the finite difference framework developed in \cite{ShuOsher1} for the spatial discretization of nonlinear degenerate parabolic equations. 
The WENO method based on the finite difference framework in \cite{ShuOsher1} is referred to as the alternative finite difference WENO (A-WENO) method in \cite{AlternativeFDWENO} to distinguish it from the more commonly practiced one \cite{WENOJS,ShuOsher2}.
Significant advantages of the A-WENO method include that arbitrary monotone fluxes can be used in this framework \cite{AlternativeFDWENO}, and it allows for a flexible choice of interpolated variables \cite{RIWENO, WBAWENO}.
 
Efficient and high-order temporal numerical schemes are crucial for the performance of high accuracy numerical simulations of time-dependent PDEs. 
For purely hyperbolic equations, explicit methods, e.g., the strong stability preserving Runge-Kutta (SSP-RK) methods and multi-step methods \cite{shu1988total,gottlieb2001strong,gottlieb2009high,gottlieb2011strong}, and Lax-Wendroff type methods \cite{qiu2003finite,li2016two,xu2022third} are widely adopted, as these explicit schemes require less computational effort per time step than implicit schemes and simulations are generally efficient with time-step size $\Delta t $ proportional to spatial grid size $\Delta x$ (i.e., $\Delta t \sim O (\Delta x)$) under the Courant-Friedrichs-Lewy (CFL) condition. 
However, in the context of parabolic equations, time-step sizes of explicit schemes under the CFL condition are much more restrictive which requires $\Delta t \sim O (\Delta x^2)$, due to the stiffness introduced by the diffusion terms. 
In such scenarios, implicit schemes are more desirable to allow for considerably larger time-step sizes, albeit with a significant rise in computational cost per time step.
Computational efficiency of fully implicit schemes can be improved by semi-implicit or implicit-explicit (IMEX) methods, e.g., \cite{Ascher, Kennedy, boscarino2022high,huang2022high}. For solving nonlinear degenerate parabolic equations by high-order WENO schemes,
most of work in the literature uses an explicit approach and small time-step sizes are required, which lead to relatively high computational costs. Recently in  
\cite{arbogast2019finite, zhang2022high}, implicit and semi-implicit approaches for WENO schemes were developed to show much better efficiency than explicit schemes. Another class of efficient temporal schemes for stiff problems are exponential integrators \cite{Oster}. The exponential integrators perform an ``exact'' integration of stiff linear part of the problem to remove the severe time-step size restriction. Among explicit exponential integrators, a popular class of methods are the exponential time differencing (ETD) approaches, e.g., ETD multi-step methods and ETD Runge-Kutta (ETD-RK) methods \cite{Beylk, cox2002exponential, DuQ, kassam2005fourth,  liu2023exponential}. ETD schemes have the advantages such as relatively small numerical errors, good steady-state preservation property, etc. \cite{Beylk, kassam2005fourth}. To deal with stiff nonlinear reaction terms in reaction-diffusion equations and advection-diffusion-reaction equations, a class of implicit
exponential integrators, called ``implicit integration factor'' (IIF) methods \cite{nie2006efficient, ChenZ, jiang2016krylov, LuZ, liu2019krylov}, were developed such that the implicit terms are free of the exponential operation of the stiff linear terms for achieving efficient computations. 

For the fully nonlinear stiff problems such as the degenerate parabolic equations considered in this paper,
the exponential integrators such as the ETD schemes are still very efficient methods to tackle the stiffness explicitly and accurately. Here we adopt the ETD-RK methods in the time stepping due to their larger stability regions than the ETD multi-step methods \cite{cox2002exponential}.  To handle the full nonlinearity, a popular approach is to use the exponential Rosenborg-type method \cite{hochbruck2009exponential} or the exponential propagation iterative method \cite{tokman2006efficient}, which form the linear stiff part for applying the exponential integrators by computing the Jacobian of the nonlinear stiff system in each time step around the numerical solution. However, for highly nonlinear high-order spatial schemes such as the multi-resolution A-WENO schemes used here, it is a challenging problem how to efficiently form the linear stiff part in applying the exponential integrators, since direct computation of
a Jacobian matrix based on high-order A-WENO discretizations of the nonlinear diffusion terms
to obtain the linear stiff part is very complicated and expensive. In this paper, we propose a novel
and effective approach to solve this difficulty by replacing the exact Jacobian of a high-order multi-resolution A-WENO scheme with that of the corresponding high-order linear scheme in the ETD-RK time marching,
based on the fact that in smooth regions the nonlinear weights of WENO schemes closely approximate the corresponding linear weights, while in non-smooth regions the stiff diffusion degenerates.

The remaining of the paper is organized as follows. In Section \ref{sect:2}, we formulate the numerical methods in details for the nonlinear degenerate parabolic equations, incorporating the multi-resolution A-WENO spatial discretization and ETD-RK time integration by the proposed novel approach. 
A comprehensive set of numerical experiments and comparisons with traditional methods are performed in Section \ref{sect:3} to demonstrate the
effectiveness of the proposed approach, and verify high-order accuracy, nonlinear stability and high efficiency of the new methods. The stiffness of the nonlinear degenerate parabolic PDEs is resolved well, and large time-step
size computations of $\Delta t \sim O (\Delta x)$ are achieved. We conclude the paper with remarks in Section \ref{sect:4}.

\section{The numerical methods}\label{sect:2}
%%%%%%%%%%%%%%%%%%%%%%%%%%%%%%%%%%%%%%%%%%%%%%%%%%%%%%%%%%%%%%%%%%%%

In this section, we present the details of the numerical methods that will be employed in the subsequent numerical tests. 
The spatial discretization is based on the multi-resolution WENO interpolations in the alternative formulation of the finite difference method.
The ETD-RK methods are adopted to evolve the nonlinear ODE system resulted from the spatial discretization, with the linear stiff component of the exponential integrators derived from the corresponding high-order linear schemes for the diffusion terms.

\subsection{Spatial discretization}

\subsubsection{Discretization for the diffusion terms}\label{Sect:DiffDisc}

We first consider the parabolic equations with only diffusion terms. In one space dimension, the general form of the equation is
\begin{equation}\label{eq:DiffEq1D}
u_t = g(u)_{xx}.
\end{equation}
A uniform grid $\cdots < x_0 < x_1 < x_2 < \cdots$ is adopted for computation, where $x_i, i = 0, \pm 1, \ldots$ are the grid points, and $\Delta x$ is the grid size. The finite difference method for the diffusion equation \eqref{eq:DiffEq1D} is formulated as follows:
\begin{equation}\label{eq:Diff1DFDM}
\frac{d u_i}{dt} = \frac{\hat{b}_{i+\frac12} - \hat{b}_{i-\frac12}}{\Delta x},\quad i = 0, \pm1, \ldots,
\end{equation}
where the grid function $\{u_i\}_{\forall i}$ approximates the function values $u(x_i, t)$ at grid points, and $\hat{b}_{i+\frac12}$ is the discretization for the flux $b := g(u)_x$ at the interface $x_{i+\frac12} := (x_i+x_{i+1})/2$. 
Here and henceforth, we denote $g_{i} = g(u_i)$ for brevity.

To achieve the $2r$-th order accurate approximation
\begin{equation*}
\frac{\hat{b}_{i+\frac12}-\hat{b}_{i-\frac12}}{\Delta x}=g(u)_{xx}|_{x=x_{i}}+O(\Delta x^{2r}),
\end{equation*}
the numerical flux $\hat{b}_{i+\frac12}$, which depends on $g_{i-r+1},\ldots, g_{i+r}$ in the alternative finite difference formulation \cite{ShuOsher1, jiang2021high}, is taken as
\begin{equation}\label{eq:DiffFlux}
\hat{b}_{i+\frac12}=\sum_{m=0}^{r-1}c_{m}\Delta x^{2m}\widehat{\left(\partial_x^{2m+1}{g}\right)}_{i+\frac12},
\end{equation}
where the constant coefficients $c_0=1, c_1=-\frac{1}{24}, c_2=\frac{7}{5760}, c_3=-\frac{31}{967680},\ldots$ are determined from Taylor expansion to attain the designed accuracy, and $\widehat{(\partial_x^{2m+1}g)}_{i+\frac12}$ are approximations to $\partial_{x}^{2m+1}g$ at the interface $x_{i+\frac12}$, for $m=0,1,\ldots,r-1$.

When the parabolic equation is degenerate, it is desired to control spurious oscillations around sharp fronts and achieve nonlinear stability of simulation.
Following \cite{jiang2021high}, we adopt the multi-resolution WENO interpolations to obtain the lowest order term $\widehat{(\partial_x g)}_{i+\frac12}$, and use central numerical differentiation (linear scheme) to calculate the high order ones $\widehat{(\partial_x^{2m+1}g)}_{i+\frac12}, m\geq 1$.
The $2r$-th order multi-resolution WENO procedure for obtaining $\widehat{(\partial_x g)}_{i+\frac12}$ is outlined as follows.
\begin{itemize}
    \item Step 1: find Lagrange interpolation on nested stencils.
    
    In the multi-resolution WENO, we choose the $r$ nested stencils $S^{(1)}=\{x_{i}, x_{i+1}\}$, $S^{(2)}=\{x_{i-1}, x_{i}, x_{i+1}, x_{i+2}\},$..., $S^{r}=\{x_{i-r+1}, \ldots, x_{i+r}\}$, and perform the Lagrange interpolation on the stencils, i.e., find polynomials $p_{k}(x)$ of degree $2k-1$, for $k=1,2,\ldots, r$, such that
    \begin{equation}
    p_{k}(x_{j})=g_{j},\quad x_{j}\in S^{(k)}.
    \end{equation}
    Consequently, we approximate $g'(x_{i+\frac12})$ by $p_{k}'(x_{i+\frac12})$ with accuracy $O(\Delta x^{2k})$:
    % with the accuracy
    % \begin{equation}
    %     p'_{k}(x_{i+\frac12})=\partial_x g(x_{i+\frac12})+O(\Delta x^{2k}), \quad k=1,2,\ldots, r,
    % \end{equation}
    % for sufficiently smooth grid functions.
    % For example, we can calculate the following expressions,
    \begin{equation}
    \begin{split}
    p_1'(x_{i+\frac12})=&\frac{1}{\Delta x}(-g_i + g_{i+1}),\\
    p_2'(x_{i+\frac12})=&\frac{1}{24\Delta x}(g_{i-1} - 27g_{i} + 27g_{i+1} - g_{i+2}),\\
    p_3'(x_{i+\frac12})=&\frac{1}{1920\Delta x}(-9g_{i-2} + 125g_{i-1} - 2250g_i + 2250g_{i+1} - 125g_{i+2} + 9g_{i+3}),\\
    p_4'(x_{i+\frac12})=&\frac{1}{107520\Delta x}(75g_{i-3} - 1029g_{i-2} + 8575g_{i-1} - 128625g_i + 128625g_{i+1} \\
    & - 8575g_{i+2} + 1029g_{i+3} - 75g_{i+4}),\\
    \cdots
    \end{split}
    \end{equation}
    \item Step 2: form linear weights for interpolation.
    
    We write the highest order interpolation polynomial $p_{r}(x)$ as a convex combination of polynomials $q_1(x), q_{2}(x), \ldots, q_{r}(x)$ with linear coefficients $\{d_{k}\}_{k=1}^{r}$ as follows:
    \begin{equation}\label{eq:LinearInterpDiff}
    p_{r}(x)=\sum_{k=1}^{r}d_{k} q_{k}(x),
    \end{equation}
    where 
    \begin{equation}
        d_k=\theta_{k,r},\quad \theta_{n,m}:=\frac{\tilde{\theta}_{n,m}}{\sum_{\ell=1}^{m}\tilde{\theta}_{\ell,m}},\quad 1\leq n\leq m\leq r,
    \end{equation}
    and 
    \begin{equation}
    \begin{split}
    q_{1}(x)=&p_1(x),\\
    q_{k}(x)=&\frac{1}{\theta_{k,k}}p_{k}(x)-\sum_{\ell=1}^{k-1}\frac{\theta_{\ell,k}}{\theta_{k,k}}q_{\ell}(x),\quad k=2,3,\ldots, r.
    \end{split}
    \end{equation}
    Following the practice in \cite{zhu2018new, jiang2021high}, we take $\tilde{\theta}_{n,m}=10^{n-1}$.
    Consequently, we have the linear weights
    \begin{equation}
    \begin{split}
    d_1=&1, \quad \text{for}~r=1,\\
    d_1=&\frac{1}{11}, d_2=\frac{10}{11},\quad \text{for}~r=2,\\
    d_1=&\frac{1}{111}, d_2=\frac{10}{111}, d_3=\frac{100}{111}, \quad \text{for}~r=3,\\
    d_1=&\frac{1}{1111}, d_2=\frac{10}{1111}, d_3=\frac{100}{1111}, d_4=\frac{1000}{1111}, \quad \text{for}~ r=4,\\
    \cdots&
    \end{split}
    \end{equation}
    and 
    \begin{equation}
    \begin{split}
    q_1(x)=&p_1(x),\\
    q_2(x)=&\frac{11}{10} p_2(x)-\frac{1}{10}q_1(x),\\
    q_3(x)=&\frac{111}{100}p_3(x)-\frac{1}{10}q_2(x)-\frac{1}{100}q_1(x),\\
    q_4(x)=&\frac{1111}{1000}p_4(x)-\frac{1}{10}q_3(x)-\frac{1}{100}q_2(x)-\frac{1}{1000}q_1(x),\\
    \cdots&
    \end{split}
    \end{equation}
    \item Step 3: compute nonlinear weights for interpolation.

    The nonlinear weights $\omega_{k}$ are obtained based on the smoothness indicators $\beta_k$ and the linear weights $d_{k}$ for the polynomials $q_k(x)$, $k=1,2,\ldots, r$.
    Following the practice in \cite{jiang2021high}, we use the smoothness indicator measuring the smoothness of $p_k(x)$ instead of $q_k(x)$ on the interval $[x_{i}, x_{i+1}]$, as $p_k(x)$ is the major component of $q_k(x)$. It is defined as 
    \begin{equation}\label{eq:smoothindDiff}
    \beta_k=\sum_{m=1}^{2k-1}\Delta x^{2m-1}\int_{x_{i}}^{x_{i+1}}\left(\frac{d^m p_{k}(x)}{dx^m}\right)^{2}dx.
    \end{equation}
    The detailed expressions of $\beta_1,\beta_2,\beta_3$ and $\beta_4$ are given in the Appendix \ref{appd:diffusion}.
    Consequently, we calculate the nonlinear weights $\omega_k$ as
    \begin{equation}
    \omega_k=\frac{\tilde{\omega}_{k}}{\sum_{m=1}^{r} {\tilde{\omega}_{m}}},\quad \tilde{\omega}_k=d_k\left(1+\frac{\tau_r}{\beta_k+\epsilon}\right),\quad k=1,2,\ldots,r,
    \end{equation}
    where $\tau_r=\left(\sum_{m=1}^{r-1}|\beta_{m}-\beta_{r}|\right)^{\frac{r+1}{2}}$, and $\epsilon$ is a small positive number to avoid dividing by zero, e.g., $\epsilon=10^{-10}$.
    
    \item Step 4: compute multi-resolution WENO interpolation.
    
    Finally, we obtain the multi-resolution WENO interpolation polynomial by replacing the linear weights with the nonlinear weights in \eqref{eq:LinearInterpDiff}:
    \begin{equation*}
        Q(x)=\sum_{k=1}^{r}\omega_{k} q_{k}(x),
    \end{equation*}
    and consequently, 
    \begin{equation}\label{eq:Diff1D1stTerm}
        \widehat{(\partial_x g)}_{i+\frac12}=\sum_{k=1}^{r}\omega_{k} q_{k}'(x_{i+\frac12}).
    \end{equation}
\end{itemize}

It remains to calculate the high-order derivative terms in \eqref{eq:DiffFlux}, which is based on linear central numerical differentiation with much less computational costs than the first-order derivative term. 
For different orders of approximation, one can calculate:
\begin{itemize}
    \item $r=2$:
    \begin{equation}\label{eq:Diff1DOtherTermsr2}
    c_{1}\Delta x^{2}\widehat{\left(\partial_x^{3}{g}\right)}_{i+\frac12}=\frac{1}{24\Delta x}(g_{i-1} - 3g_{i} + 3g_{i+1} - g_{i+2}),
    \end{equation}
    \item $r=3$:
    \begin{equation}\label{eq:Diff1DOtherTermsr3}
    \begin{split}
    \sum_{m=1}^{r-1}c_{m}\Delta x^{2m}\widehat{\left(\partial_x^{2m+1}{g}\right)}_{i+\frac12}=&\frac{1}{5760\Delta x}(-37g_{i-2} + 425g_{i-1} - 1090g_i \\
    &+ 1090g_{i+1} - 425g_{i+2} + 37g_{i+3}),
    \end{split}    
    \end{equation}
    \item $r=4$:
    \begin{equation}\label{eq:Diff1DOtherTermsr4}
    \begin{split}
    \sum_{m=1}^{r-1}c_{m}\Delta x^{2m}\widehat{\left(\partial_x^{2m+1}{g}\right)}_{i+\frac12}=&\frac{1}{322560\Delta x}(351 g_{i-3} - 4529g_{i-2} + 31171g_{i-1} - 73325g_i \\
    & + 73325g_{i+1} - 31171g_{i+2} + 4529g_{i+3} - 351g_{i+4}).\\
    \end{split}
    \end{equation}
\end{itemize}

The $2r$-th order alternative formulation of finite difference multi-resolution WENO scheme for the one-dimensional diffusion equation \eqref{eq:DiffEq1D} is then obtained by integrating \eqref{eq:Diff1DFDM}, \eqref{eq:DiffFlux}, \eqref{eq:Diff1D1stTerm}, \eqref{eq:Diff1DOtherTermsr2} (or \eqref{eq:Diff1DOtherTermsr3}, \eqref{eq:Diff1DOtherTermsr4}).

Likewise, for the two-dimensional diffusion equation
\begin{equation}\label{eq:DiffEq2D}
    u_t=g_1(u)_{xx}+g_2(u)_{yy}, 
\end{equation}
we adopt the uniform grid $\cdots<x_0<x_1<x_2<\cdots$ and $\cdots<y_0<y_1<y_2<\cdots$ for the $x$ and $y$ directions, respectively, where $(x_i,y_j), i,j=0,\pm1,\ldots$ are grid points, and $\Delta x, \Delta y$ are the grid sizes. 

The alternative formulation of finite difference method for the two-dimensional diffusion equation \eqref{eq:DiffEq2D} is formulated as follows:
\begin{equation}\label{eq:Diff2DFDM}
\frac{d u_{i,j}}{dt}=\frac{\hat{b}^{1}_{i+\frac12,j}-\hat{b}^{1}_{i-\frac12,j}}{\Delta x}+\frac{\hat{b}^{2}_{i,j+\frac12}-\hat{b}^{2}_{i,j-\frac12}}{\Delta y},\quad i,j=0,\pm1,\ldots,
\end{equation}
where the grid function $\{u_{i,j}\}_{\forall i,j}$ approximates the function values $u(x_i,y_j,t)$ at grid points, and $\hat{b}^{1}_{i+\frac12,j}$ and $\hat{b}^{2}_{i, j+\frac12}$ are the discretizations for the fluxes $b^1:=g_1(u)_x$ and $b^2:=g_2(u)_y$ at $(x_{i+\frac12},y_{j})$ and $(x_{i}, y_{j+\frac12})$, respectively.
Since we are using finite difference method, the discretization in \eqref{eq:Diff2DFDM} is performed in a dimension by dimension fashion, and the computation for each dimension follows exactly the same procedure as we have established for the one-dimensional case. 

The resulting ODE system from the scheme \eqref{eq:Diff1DFDM} for the one-dimensional problems or the scheme \eqref{eq:Diff2DFDM} for the two-dimensional problems is generally denoted by
\begin{equation}\label{eq:DiffODEs}
\mathbf{u}_t=G(\mathbf{u}),
\end{equation}
where $\mathbf{u}$ is the vector $\{u_{i}\}_{\forall i}$ in the one-dimensional case, or $\{u_{i,j}\}_{\forall i,j}$ in the two-dimensional case, at grid points.

To prepare for the ETD time evolution method in the later subsections, we describe the corresponding linear scheme of \eqref{eq:DiffODEs} as follows. This scheme is obtained by replacing the nonlinear weights $\omega_k$'s with the linear weights $d_k$'s in the WENO interpolation \eqref{eq:Diff1D1stTerm}:
\begin{equation}\label{eq:DiffODEsLinear}
\mathbf{u}_t = G_L(\mathbf{u}).
\end{equation}
In one space dimension, the $2r$-th order linear scheme \eqref{eq:DiffODEsLinear} is given as follows:
\begin{itemize}
    \item $r=1$:
    \begin{equation}
        \frac{du_i}{dt}=\frac{1}{\Delta x^2}(g_{i-1}-2g_{i}+g_{i+1}),
    \end{equation}
    \item $r=2$:
        \begin{equation}
        \frac{du_i}{dt}=\frac{1}{12\Delta x^2}(-g_{i-2}+16g_{i-1}-30g_{i}+16g_{i+1}-g_{i+2}),
    \end{equation}
    \item $r=3$:
        \begin{equation}
        \frac{du_i}{dt}=\frac{1}{180\Delta x^2}(2g_{i-3}-27g_{i-2}+270g_{i-1}-490g_{i}+270g_{i+1}-27g_{i+2}+2g_{i+3}),
    \end{equation}
    \item $r=4$:
        \begin{equation}
        \begin{split}
            \frac{du_i}{dt}=&\frac{1}{5040\Delta x^2}(-9g_{i-4}+128g_{i-3}-1008g_{i-2}+8064g_{i-1}-14350g_{i}\\
            &+8064g_{i+1}-1008g_{i+2}+128g_{i+3}-9g_{i+4}).
        \end{split}
    \end{equation}
\end{itemize}
The formulation of the $2r$-th order linear scheme \eqref{eq:DiffODEsLinear} in two space dimensions is obtained in a dimension by dimension manner.

\subsubsection{Discretization for the convection terms}\label{Sect:ConvDisc}
In this section, we describe the multi-resolution WENO discretization for the convection terms in the degenerate
parabolic equations, using the alternative formulation of finite difference method. Consider the one-dimensional convection equation
\begin{equation}
u_t+f(u)_x=0.
\end{equation}
The semi-discrete conservative finite difference scheme is formulated as 
\begin{equation}
\frac{d u_i}{dt}=-\frac{\hat{f}_{i+\frac12}-\hat{f}_{i-\frac12}}{\Delta x},\quad i=0,\pm 1,\ldots,
\end{equation}
where $\{u_{i}\}_{\forall i}$ is the grid function to approximate $u(x_i,t)$ at grid points, and $\hat{f}_{i+\frac12}$ is the numerical flux at the interface $x_{i+\frac12}$.
%For the time being, we denote $f_i=f(u_{i})$.

To achieve the $2r$-th order accurate approximation, the numerical flux $\hat{f}_{i+\frac12}$, which depends on $u_{i-r},\ldots, u_{i+r+1}$, is defined in a similar fashion as for the diffusion equations:
\begin{equation}\label{eq:ConvFlux}
\hat{f}_{i+\frac12}=h(u^{-}_{i+\frac12}, u^{+}_{i+\frac12})+\sum_{m=1}^{r-1}c_{m}\Delta x^{2m}\widehat{\left(\partial_x^{2m}{f}\right)}_{i+\frac12}.
\end{equation}
Here the coefficients $c_1, c_2, c_3,\ldots$ are the same as in \eqref{eq:DiffFlux}. $h(\cdot,\cdot)$ is the numerical flux function of $f(u)$ obtained from an exact or approximate Riemann solver, e.g., the Lax-Friedrichs flux. $u_{i+\frac12}^{\pm}$ is computed using the multi-resolution WENO interpolation for the grid function $\{u_{j}\}_{\forall j}$ at the interface $x_{i+\frac12}$ with left/right bias, and $\widehat{(\partial_x^{2m}f)}_{i+\frac12}$ are approximations to $\partial_{x}^{2m}f$ at $x_{i+\frac12}$, for $m=1,\ldots,r-1$.
In the following, we first outline
the $(2r+1)$-th order multi-resolution WENO procedure for computing $u_{i+\frac12}^{-}$.
The procedure of the multi-resolution WENO interpolation for $u_{i+\frac12}^{+}$ follows similarly in that it
is mirror-symmetric with respect to $x_{i+\frac12}$.

\begin{itemize}
    \item Step 1: find Lagrange interpolation on nested stencils.

    In the multi-resolution WENO interpolation, we choose the nested stencils $S^{(0)}=\{x_{i}\}$, $S^{(1)}=\{x_{i-1}, x_{i}, x_{i+1}\}, \ldots$, $S^{(r)}=\{x_{i-r}, \ldots, x_{i+r}\}$, and perform the Lagrange interpolation on these stencils, i.e., find polynomials $p_{k}(x)$ of degree $2k$ for $k=0,1,\ldots, r$, such that,
    \begin{equation}
        p_k(x_j)=u_{j},\quad x_{j}\in S^{(k)}.
    \end{equation}
    Consequently, we use $p_{k}(x_{i+\frac12})$ to approximate $u(x_{i+\frac12})$ with accuracy $\Delta x^{2k+1}$ as follows:
    \begin{equation}
    \begin{split}
        p_{0}(x_{i+\frac12})=&u_{i},\\
        p_{1}(x_{i+\frac12})=&\frac{1}{8} (-u_{i-1} + 6 u_{i} + 3 u_{i+1}),\\
        p_{2}(x_{i+\frac12})=&\frac{1}{128}(3 u_{i-2} - 20 u_{i-1} + 90 u_{i} + 60 u_{i+1} - 5 u_{i+2}),\\
        p_{3}(x_{i+\frac12})=&\frac{1}{1024}(-5 u_{i-3} + 42 u_{i-2} - 175 u_{i-1} + 700 u_{i} + 525 u_{i+1} - 70 u_{i+2} + 7 u_{i+3}),\\
        p_{4}(x_{i+\frac12})=&\frac{1}{32768}(35 u_{i-4} -360 u_{i-3} + 1764 u_{i-2} - 5880 u_{i-1} + 22050 u_{i} + 17640 u_{i+1} \\
        & - 2940 u_{i+2} + 504 u_{i+3} - 45 u_{i+4}),\\
        \cdots&
    \end{split}
    \end{equation}
    \item Step 2: form linear weights for interpolation. 
        
        We write the highest order interpolation polynomial $p_r(x)$ as a convex combination of polynomials $q_0(x), q_1(x),\ldots, q_r(x)$ with linear coefficients $\{d_{k}\}_{k=0}^{r}$ as follows:
        \begin{equation}\label{eq:LinearInterpConv}
            p_r(x)=\sum_{k=0}^{r}d_{k} q_{k}(x),
        \end{equation}
    where
    \begin{equation}
        d_k=\theta_{k,r},\quad \theta_{n,m}:=\frac{\tilde{\theta}_{n,m}}{\sum_{\ell=0}^{m}\tilde{\theta}_{\ell,m}},\quad 0\leq n\leq m\leq r,
    \end{equation}
    and 
    \begin{equation}
    \begin{split}
    q_{0}(x)=&p_0(x),\\
    q_{k}(x)=&\frac{1}{\theta_{k,k}}p_{k}(x)-\sum_{\ell=0}^{k-1}\frac{\theta_{\ell,k}}{\theta_{k,k}}q_{\ell}(x),\quad k=1,2,\ldots, r.
    \end{split}
    \end{equation}
    In practice, we take $\tilde{\theta}_{n,m}=10^{n}$.
    Consequently, we have the linear weights
    \begin{equation}
    \begin{split}
    d_0=&1, \quad \text{for}~r=0,\\
    d_0=&\frac{1}{11}, d_1=\frac{10}{11},\quad \text{for}~r=1,\\
    d_0=&\frac{1}{111}, d_1=\frac{10}{111}, d_2=\frac{100}{111}, \quad \text{for}~r=2,\\
    d_0=&\frac{1}{1111}, d_1=\frac{10}{1111}, d_2=\frac{100}{1111}, d_3=\frac{1000}{1111}, \quad \text{for}~ r=3,\\
    d_0=&\frac{1}{11111}, d_1=\frac{10}{11111}, d_2=\frac{100}{11111}, d_3=\frac{1000}{11111}, d_4=\frac{10000}{11111}\quad \text{for}~ r=4,\\
    \cdots&
    \end{split}
    \end{equation}
    and 
    \begin{equation}
    \begin{split}
    q_0(x)=&p_0(x),\\
    q_1(x)=&\frac{11}{10} p_1(x)-\frac{1}{10}q_0(x),\\
    q_2(x)=&\frac{111}{100}p_2(x)-\frac{1}{10}q_1(x)-\frac{1}{100}q_0(x),\\
    q_3(x)=&\frac{1111}{1000}p_3(x)-\frac{1}{10}q_2(x)-\frac{1}{100}q_1(x)-\frac{1}{1000}q_0(x),\\
    q_4(x)=&\frac{11111}{10000}p_4(x)-\frac{1}{10}q_3(x)-\frac{1}{100}q_2(x)-\frac{1}{1000}q_1(x)-\frac{1}{10000}q_0(x),\\
    \cdots&
    \end{split}
    \end{equation}    
    \item Step 3: compute nonlinear weights for interpolation.
    
    For $k=1, \ldots, r$,
    we compute the nonlinear weight $\omega_{k}$ for the multi-resolution WENO interpolation, based on the smoothness indicator $\beta_k$ which measures the smoothness of $q_{k}(x)$ on the interval $[x_{i-\frac12},x_{i+\frac12}]$:
    \begin{equation}\label{eq:WENOweightsconv}
    \beta_k=\sum_{m=1}^{2k}\Delta x^{2m-1}\int_{x_{i-\frac12}}^{x_{i+\frac12}}\left(\frac{d^m q_{k}(x)}{dx^m}\right)^{2}dx,\quad k\geq 1.
    \end{equation}
    The detailed expressions of $\beta_1,\beta_2,\beta_3$ and $\beta_4$ are given in the Appendix B. 
    The only exception for \eqref{eq:WENOweightsconv} is $\beta_0$, which is magnified from zero to a tiny value. See \cite{zhu2018new} for details.
    Consequently, we calculate the nonlinear weights $\omega_k$ as
    \begin{equation}
    \omega_k=\frac{\tilde{\omega}_{k}}{\sum_{m=0}^{r} {\tilde{\omega}_{m}}},\quad \tilde{\omega}_k=d_k\left(1+\frac{\tau_r}{\beta_k+\epsilon}\right),\quad k=0,1,2,\ldots,r,
    \end{equation}
    where $\tau_r=\left(\frac{1}{r}\sum_{m=0}^{r-1}|\beta_{m}-\beta_{r}|\right)^{r}$, and $\epsilon$ is a small positive number to avoid dividing by zero, e.g. $\epsilon=10^{-10}$.
    
    \item Step 4: compute multi-resolution WENO interpolation.
    
    Finally, the multi-resolution WENO interpolation $u_{i+\frac12}^{-}$ is obtained by replacing the linear weights $d_k$'s with the nonlinear weights $\omega_k$'s in \eqref{eq:LinearInterpConv} and evaluating at $x_{i+\frac12}$:
    \begin{equation*}
        u_{i+\frac12}^{-}=\sum_{k=0}^{r}\omega_{k} q_{k}(x_{i+\frac12}).
    \end{equation*}
\end{itemize}
    % # To calculate β₁
    % κ = 3 # κ = 1,2,3,4 for the third-order, fifth-order, seventh-order, ninth-order approximations, respectively.
    % ζ₀ = (v4-v3)^2
    % ζ₁ = (v5-v4)^2
    % if ζ₀≥ζ₁
    %     γ01=1/11; γ11=10/11;
    % else
    %     γ01=10/11; γ11=1/11;
    % end
    % σ₀ = γ01*(1+abs(ζ₀-ζ₁)^κ/(ζ₀+ϵ))
    % σ₁ = γ11*(1+abs(ζ₀-ζ₁)^κ/(ζ₁+ϵ))
    % σ = σ₀ + σ₁
    % β₁ = ( σ₀*(v4-v3) + σ₁*(v5-v4) )^2/σ^2

It remains to calculate the high-order derivative terms in \eqref{eq:ConvFlux}, which is based on linear central numerical differentiation with much less computational costs than the first term. 
Here and henceforth, we denote $f(u_i)$ by $f_i$ for brevity.
The high-order derivative terms are computed as follows:
\begin{itemize}
    \item $r=2$:
    \begin{equation}
        c_1\Delta x^2(\widehat{\partial_x^{2}f})_{i+\frac12} = \frac{1}{48}(-f_{i-1} + f_{i} + f_{i+1} - f_{i+2}),
    \end{equation}
    \item $r=3$:
    \begin{equation}
    \begin{split}
            \sum_{m=1}^{r-1}c_{m}\Delta x^{2m}\widehat{\left(\partial_x^{2m}{f}\right)}_{i+\frac12} =& \frac{1}{3840}(19f_{i-2} - 137f_{i-1} + 118f_{i} \\
            &+118f_{i+1} - 137f_{i+2} + 19f_{i+3}),
    \end{split}
    \end{equation}
    \item $r=4$:
    \begin{equation}
    \begin{split}
        \sum_{m=1}^{r-1}c_{m}\Delta x^{2m}\widehat{\left(\partial_x^{2m}{f}\right)}_{i+\frac12} =& \frac{1}{215040}(-243f_{i-3} + 2279f_{i-2} - 9859f_{i-1} + 7823f_{i} \\
        &+ 7823f_{i+1} - 9859 f_{i+2} + 2279f_{i+3} - 243f_{i+4}).
    \end{split}
    \end{equation}
\end{itemize}

The discretizations for the two-dimensional convection equation
\begin{equation}
u_t+f_1(x)_x+f_2(x)_y=0
\end{equation}
are performed in a dimension by dimension manner, with each spatial dimension being exactly the same as in the one-dimensional case.

To this end, we denote the resulting ODE system from these finite difference spatial discretizations of the convection equations by
\begin{equation}
\mathbf{u}_t=F(\mathbf{u}),
\end{equation}
for both the one-dimensional and two-dimensional cases.

\subsection{ETD temporal discretization}

\subsubsection{A new semilinearization approach}
The overall nonlinear ODE system resulting from the multi-resolution A-WENO spatial discretization for the nonlinear degenerate parabolic equation \eqref{eq:PDE1D} (or \eqref{eq:PDE2D})  is denoted by
\begin{equation}\label{eq:ODEs}
    \mathbf{u}_t=G(\mathbf{u})+F(\mathbf{u}):=L(\mathbf{u}),
\end{equation}
where $G(\mathbf{u})$ and $F(\mathbf{u})$ correspond to the discretizations for the diffusion and convection terms of the equations, respectively.
Suitable time discretization approaches for the ODE system \eqref{eq:ODEs}, which is stiff due to $G(\mathbf{u})$ from the diffusion terms, are desired to obtain efficient and accurate approximation to the solution.
Implicit schemes are often used to solve such kind of stiff systems and achieve large time-step size computations. 
Challenges in designing implicit methods include developing efficient iterative solvers for nonlinear algebraic systems, computing the exact Jacobian matrices of nonlinear terms which are often highly complex for nonlinear schemes (e.g., WENO methods \cite{gottlieb2006fifth}), etc. 

As an alternative approach for solving stiff systems efficiently, the exponential integrator methods are originally designed to solve the semilinear ODEs with stiff linear  part and non-stiff nonlinear part:
\begin{equation}\label{eq:ETD_ODEs0}
\mathbf{u}_t=C\mathbf{u}+N(\mathbf{u}),
\end{equation}
where $C$ is a constant matrix and $N(\mathbf{u})$ is a nonlinear vector function of $\mathbf{u}$.
The idea behind the exponential integrator methods is to multiply \eqref{eq:ETD_ODEs0} by the integrating factor $e^{-Ct}$ to absorb the linear stiff term:
\begin{equation}\label{eq:ETD_ODEs1}
    \frac{d }{dt} (e^{-Ct}\mathbf{u})=e^{-Ct} N(\mathbf{u}).
\end{equation}
Consequently, the system can be transformed into the exact integral formula,
\begin{equation}\label{eq:ETD_ODEs2}
\mathbf{u}^{n+1}=e^{C\Delta t}\mathbf{u}^n+\int_{0}^{\Delta t}e^{(\Delta t-\tau)C} N(\mathbf{u}(t^n+\tau)) d\tau,
\end{equation}
after being integrated over one time interval $[t^{n}, t^{n+1}]$ where $t^n$ and $t^{n+1}$ are the time steps, and $\Delta t=t^{n+1}-t^{n}$. $\mathbf{u}^{n+1}$ and $\mathbf{u}^{n}$ are the numerical solutions of $\mathbf{u}$ at $t^{n+1}$ and $t^{n}$ respectively. 
Different numerical discretization strategies for the integral in \eqref{eq:ETD_ODEs2} give rise to various methods in the family of the exponential integrators, e.g., the implicit integration factor methods, ETD multi-step methods, ETD-RK methods, etc. The severe time-step size restriction imposed by the stiff linear part is removed, as this stiff component of the system is integrated exactly here.

For a fully nonlinear stiff ODE system with the general form $\mathbf{u}_t=L(\mathbf{u})$, the exponential Rosenborg-type method takes a linearization at every time step to yield the equivalent semilinear reformulation:
\begin{equation}\label{eq:ETD_ODEs3}
\mathbf{u}_t=C^n\mathbf{u} + N^n(\mathbf{u}), \quad t\in[t^n, t^{n+1}],
\end{equation}
where $C^n:=L'(\mathbf{u}^n)$ is the Jacobian matrix of $L(\mathbf{u})$ and $N^{n}(u):=L(\mathbf{u})-C^n\mathbf{u}$ is the nonlinear remainder.
Then the exponential integrators are ready to be applied in solving \eqref{eq:ETD_ODEs3} and remove the stiffness from the linear part of the system.  

Note that in the stiff ODE system \eqref{eq:ODEs}, 
since the stiffness arises from the parabolic term, it is reasonable to compute the Jacobian matrix $C^n$ only based on $G(\mathbf{u})$ in the semilinear reformulation \eqref{eq:ETD_ODEs3}. However, even without the complex computation of the Jacobian matrix from the hyperbolic term $F(\mathbf{u})$, due to the high nonlinearity of high-order WENO discretizations for the diffusion term, it is still very difficult and expensive to compute the exact Jacobian matrix of $G(\mathbf{u})$ to obtain $C^n$ in \eqref{eq:ETD_ODEs3}. The feasible approach we propose here to resolve this difficulty is based on 
the key observation that the stiffness of the system indeed comes from the non-degenerate region of the diffusion terms, where the solution is smooth and approximations by the nonlinear WENO discretizations are very close to these by the corresponding linear schemes.
Therefore, we adopt the Jacobian matrix of the spatial discretization of the linear scheme \eqref{eq:DiffODEsLinear} to formulate the semilinear system \eqref{eq:ETD_ODEs3}, i.e., $C^{n}=G_L'(\mathbf{u}^n)$ and $N^{n}(\mathbf{u})=G(\mathbf{u})+F(\mathbf{u})-C^n\mathbf{u}$, where $G_{L}(\mathbf{u})$ is the spatial discretization for the diffusion terms in the linear scheme \eqref{eq:DiffODEsLinear}.
More specifically, for the problems with one space dimension, we have 
\begin{itemize}
    \item $r=1$, 
    \begin{equation}
        C^{n}_{i\cdot}=\frac{1}{\Delta x^2}(0,\ldots,0, g_{i-1}', -2g_{i}', g_{i+1}', 0, \ldots,0),
    \end{equation}
    \item $r=2$,
        \begin{equation}
        C^{n}_{i\cdot}=\frac{1}{12\Delta x^2}(0, \ldots, 0, -g_{i-2}', 16g_{i-1}', -30g_{i}', 16g_{i+1}', -g_{i+2}', 0, \ldots, 0),
    \end{equation}
    \item $r=3$,
        \begin{equation}
        \begin{split}
        C^{n}_{i\cdot}=&\frac{1}{180\Delta x^2}(0, \ldots, 0, 2g_{i-3}', -27g_{i-2}', 270g_{i-1}', -490g_{i}',\\
        & 270g_{i+1}', -27g_{i+2}', 2g_{i+3}', 0, \ldots, 0),            
        \end{split}
    \end{equation}
    \item $r=4$,
        \begin{equation}
        \begin{split}
            C^{n}_{i\cdot}=&\frac{1}{5040\Delta x^2}(0, \ldots, 0, -9g_{i-4}', 128g_{i-3}', -1008g_{i-2}', 8064g_{i-1}', -14350g_{i}',\\
            &8064g_{i+1}', -1008g_{i+2}', 128g_{i+3}', -9g_{i+4}', 0, \ldots, 0),
        \end{split}
    \end{equation}
\end{itemize}
where $C^{n}_{i\cdot}$ denotes the $i$-th row of the matrix $C^{n}$, $g_{k}'=g'(u_{k})$, and the lower ellipsis stands for zero components of the vectors.
Likewise, the Jacobian matrix of the spatial discretization of the linear scheme for the problems with two space dimensions is computed in a similar manner.

In the consequent subsections, we present an efficient and accurate class of explicit exponential integrators, the ETD Runge-Kutta schemes \cite{cox2002exponential}, for solving \eqref{eq:ETD_ODEs3} using the new semilinearization approach proposed above.
For conciseness, we drop the superscript $n$.

\subsubsection{ETD-RK methods}

It is important to note that the matrices involved in the algorithms are highly sparse.
Therefore, the data structure used for storing and computing these matrices is specialized for sparse matrices in our implementation.

The first-order ETD scheme is derived through approximation to the nonlinear part $N(\mathbf{u})$ in the integrand \eqref{eq:ETD_ODEs2} by the constant quantity $N(\mathbf{u}^n)$, and a direct evaluation of the integral gives
\begin{equation}
    \begin{split}
    \mathbf{u}^{n+1}=&e^{\Delta t C}\mathbf{u}^{n}+\Delta t \varphi_{1}(\Delta t C) N(\mathbf{u}^n)\\
    =&\mathbf{u}^n+\Delta t\varphi_{1}(\Delta tC)(C\mathbf{u}^n+N(\mathbf{u}^n)),        
    \end{split}
\end{equation}
where $\varphi_1(z):=\frac{e^{z}-1}{z}$ is a member of the family of $\varphi$-functions (see e.g. \cite{Oster}).
Analogous to the classic RK methods, the high-order ETD-RK schemes were derived in  \cite{cox2002exponential}. Also see e.g. \cite{liu2023exponential} for their formulations in terms of the $\varphi$-functions. In the following, we show the third-order and the fourth-order ETD-RK methods used for solving the system \eqref{eq:ETD_ODEs3} in this paper:
\begin{itemize}
%    \item ETD-RK2:
%    \begin{equation}\label{eq:ETDRK2}
%    \begin{split}
%        \mathbf{a}^n=&\mathbf{u}^n+\Delta t \varphi_1 (\Delta t C) (C\mathbf{u}^{n}+N(\mathbf{u}^n)),\\
%        \mathbf{u}^{n+1}=&\mathbf{u}^n+\Delta t\varphi_{1}(\Delta t C)(C\mathbf{u}^n+N(\mathbf{u}^n))+\Delta t\varphi_{2}(\Delta t C)(-N(\mathbf{u}^n)+N(\mathbf{a}^n)),
%    \end{split}
%    \end{equation}
    \item ETD-RK3
    \begin{equation}\label{eq:ETDRK3}
    \begin{split}
        \mathbf{a}^n=&\mathbf{u}^n+\frac{\Delta t}{2}\varphi_1(\frac{\Delta t}{2}C)(C\mathbf{u}^n+N(\mathbf{u}^n)),\\
        \mathbf{b}^n=&\mathbf{u}^n+\Delta t\varphi_{1}(\Delta t C)(C\mathbf{u}^n-N(\mathbf{u}^n)+2N(\mathbf{a}^n)),\\
        \mathbf{u}^{n+1}=&\mathbf{u}^n+\Delta t\varphi_{1}(\Delta t C)(C\mathbf{u}^n+N(\mathbf{u}^n))\\
        &+\Delta t\varphi_{2}(\Delta t C)(-3N(\mathbf{u}^n)+4N(\mathbf{a}^n)-N(\mathbf{b}^n))\\
        &+\Delta t\varphi_{3}(\Delta t C)(4N(\mathbf{u}^n)-8N(\mathbf{a}^n)+4N(\mathbf{b}^n)),
    \end{split}
    \end{equation}
    \item ETD-RK4
    \begin{equation}\label{eq:ETDRK4}
    \begin{split}
        \mathbf{a}^n=&\mathbf{u}^n+\frac{\Delta t}{2}\varphi_{1}(\frac{\Delta t}{2}C)(C\mathbf{u}^n+N(\mathbf{u}^n)),\\
        \mathbf{b}^n=&\mathbf{u}^{n}+\frac{\Delta t}{2}\varphi_{1}(\frac{\Delta t}{2} C)(C\mathbf{u}^n+N(\mathbf{a}^n)),\\
        \mathbf{c}^n=&\mathbf{a}^n+\frac{\Delta t}{2}\varphi_{1}(\frac{\Delta t}{2} C)(C\mathbf{a}^n-N(\mathbf{u}^n)+2N(\mathbf{b}^n)),\\
        \mathbf{u}^{n+1}=&\mathbf{u}^n+\Delta t\varphi_{1}(\Delta t C)(C\mathbf{u}^n+N(\mathbf{u}^n))\\
        &+\Delta t\varphi_{2}(\Delta t C)(-3N(\mathbf{u}^{n})+2N(\mathbf{a}^n)+2N(\mathbf{b}^n)-N(\mathbf{c}^n))\\
        &+\Delta t\varphi_{3}(\Delta t C)(4N(\mathbf{u}^n)-4N(\mathbf{a}^n)-4N(\mathbf{b}^n)+4N(\mathbf{c}^n)),
    \end{split}
    \end{equation}
\end{itemize}
where the $\varphi$-functions in \eqref{eq:ETDRK3}-\eqref{eq:ETDRK4} are defined by the recurrence relation $\varphi_{\ell}(z)=z\varphi_{\ell+1}(z)+\frac{1}{\ell !}, \ell=0,1,\ldots$ \cite{Oster}.
For example, for small values of $\ell$,
\begin{equation}
\varphi_0=e^{z},\quad\varphi_1(z)=\frac{e^z-1}{z},\quad \varphi_2(z)=\frac{e^z-1-z}{z^2},\quad \varphi_3(z)=\frac{e^{z}-1-z-\frac12 z^2}{z^3}.
\end{equation}

\subsubsection{Fast computation of $\varphi$-functions in the implementation}

Notice that every stage in the ETD-RK schemes \eqref{eq:ETDRK3}-\eqref{eq:ETDRK4} is a linear combination of $\varphi$-functions of matrices acting on a set of vectors,
\begin{equation}\label{eq:LinCombPhis}
    \varphi_0(A) \mathbf{v}_0+\varphi_1(A)\mathbf{v}_1+\cdots+\varphi_{p}(A)\mathbf{v}_{p},
\end{equation}
where $A$ is a $\mathcal{N}$-by-$\mathcal{N}$ sparse matrix, $\mathbf{v}_k\in\mathbb{R}^\mathcal{N}, k=0,\ldots,p$, and $\mathcal{N}$ is the size of the system \eqref{eq:ETD_ODEs3}. Hence fast computation of $\varphi$-functions is crucial for the efficient implementation of these ETD methods.
In this paper, the algorithm \textit{phipm} developed in \cite{niesen2012algorithm}, which is based on the Krylov subspace method for matrix exponentials (see e.g. \cite{moler2003nineteen}), is adopted for the computation of \eqref{eq:LinCombPhis}.
One may also conduct the computation based on its modified versions, e.g. \textit{phipm-simul-iom2} in \cite{luan2019further}.
For completeness, a brief description of the algorithm is given below.

The Krylov subspace with a dimension $m$ ($m \ll \mathcal{N}$) for the matrix-vector pair $A\in\mathbb{R}^{\mathcal{N}\times \mathcal{N}},\mathbf{v}\in\mathbb{R}^\mathcal{N}$ is defined as 
\begin{equation*}
    K_{m}=\text{span}\{\mathbf{v}, A\mathbf{v},\ldots,A^{m-1}\mathbf{v}\}.
\end{equation*}
The Arnoldi iteration \cite{TrefBau}, a stabilized Gram-Schmidt process, is employed to obtain an orthonormal basis of $K_m$.
We denote the basis by $\{\mathbf{\tilde v}_1,\mathbf{\tilde v}_2,\ldots, \mathbf{\tilde v}_m\}$ and let $V_{m}=[\mathbf{\tilde v}_1, \mathbf{\tilde v}_2,\ldots, \mathbf{\tilde v}_m]\in\mathbb{R}^{\mathcal{N}\times m}$, then the Hessenberg matrix $H_m=V_{m}^{T}AV_{m}\in\mathbb{R}^{m\times m}$ obtained as one of the products of the Arnoldi iteration is the matrix representation of $A$ in the Krylov subspace $K_m$ with the basis $\{\mathbf{\tilde v}_{1}, \ldots,\mathbf{\tilde v}_m\}$.
Note that $V_{m}H_{m}V_{m}^{T}$ is the matrix representation of $A$ in the Krylov subspace $K_m$ with the standard basis of $\mathbb{R}^\mathcal{N}$, so one can approximate $\varphi_{k}(A)\mathbf{v}$ in the Krylov subspace $K_m$ as
\begin{equation*}
    \varphi_{k}(A)\mathbf{v}\approx\varphi_{k}(V_{m}H_{m}V_{m}^{T})\mathbf{v}=V_{m}\varphi_{k}(H_{m})V_{m}^{T}\mathbf{v}=||\mathbf{v}||V_{m}\varphi_{k}(H_{m})\mathbf{e}_{1},
\end{equation*}
where $||\cdot||$ is the standard Euclidean norm and $\mathbf{e}_1$ is the first standard basis in $\mathbb{R}^m$.
The Krylov subspace method reduces the computation of $\varphi_{k}(A)\mathbf{v}$ down to that of $\varphi_{k}(H_{m})\mathbf{e}_1$, whose computation is discussed in \cite{niesen2012algorithm} and the references therein.

As that pointed out in \cite{niesen2012algorithm}, the linear combination \eqref{eq:LinCombPhis} is actually the solution of the initial value problem of the ordinary differential equation
\begin{equation}\label{eq:ODEphifunc}
\mathbf{y}_t=A\mathbf{y}+\mathbf{v}_1+t\mathbf{v}_2+\cdots+\frac{t^{p-1}}{(p-1)!}\mathbf{v}_p,\quad \mathbf{y}(0)=\mathbf{v}_0,
\end{equation}
at $t=1$. 
A time-stepping method built upon efficient computation of $\varphi$-function of matrices by the Krylov subspace method is used to compute $\mathbf{y}(1)$.
As the two critical parameters affecting the computational cost, the dimension $m$ of the Krylov subspace $K_m$ and the time-step size $\tau$ in the time-stepping method for \eqref{eq:ODEphifunc} are determined adaptively during computations to optimize the performance. See \cite{niesen2012algorithm} for details of the method.

\subsubsection{Other time marching approaches for comparison}
To show the efficiency of the ETD schemes with the proposed techniques for solving the complex ODE systems resulting from the multi-resolution A-WENO spatial discretization, we compare their computational costs and accuracy with some commonly used explicit and implicit strong stability preserving Runge-Kutta (SSP-RK) methods. 
Specifically, the following three-stage third order explicit Runge-Kutta (SSP-ERK3) and three-stage third order diagonally implicit Runge-Kutta (SSP-IRK3) schemes are employed for the comparison. Their forms for solving the ODE system $\mathbf{u}_t=L(\mathbf{u})$ are

\begin{itemize}
    \item 
SSP-ERK3:
\begin{equation}
\begin{split}
\mathbf{u}^{(1)}=&\mathbf{u}^{n}+\Delta t L(\mathbf{u}^{n}),\\
\mathbf{u}^{(2)}=&\frac{3}{4} \mathbf{u}^{n}+\frac14 (\mathbf{u}^{(1)}+\Delta t L(\mathbf{u}^{(1)})),\\
\mathbf{u}^{n+1}=&\frac{1}{3} \mathbf{u}^{n}+\frac{2}{3}(\mathbf{u}^{(2)}+\Delta t L(\mathbf{u}^{(2)})),
\end{split}
\end{equation}
    \item SSP-IRK3:
\begin{equation}
\begin{split}
    \mathbf{u}^{(1)} = & \mathbf{u}^{n} + \Delta t(a_{11}L(\mathbf{u}^{(1)})),\\
    \mathbf{u}^{(2)} = & \mathbf{u}^{n} + \Delta t(a_{21}L(\mathbf{u}^{(1)})+a_{22}L(\mathbf{u}^{(2)})),\\
    \mathbf{u}^{(3)} = & \mathbf{u}^{n} + \Delta t(a_{31}L(\mathbf{u}^{(1)})+a_{32}L(\mathbf{u}^{(2)})+a_{33}L(\mathbf{u}^{(3)})),\\
    \mathbf{u}^{n+1} = & \mathbf{u}^{n} + \Delta t(a_{41}L(\mathbf{u}^{(1)})+a_{42}L(\mathbf{u}^{(2)})+a_{43}L(\mathbf{u}^{(3)})),\\
\end{split}
\end{equation}
% \begin{equation}
% \begin{split}
%     \mathbf{u}^{(1)} = & \mathbf{u}^{n} + \Delta t(a_{11}L(\mathbf{u}^{(1)})),\\
%     \mathbf{u}^{(2)} = & \mathbf{u}^{n} + \Delta t(a_{21}L(\mathbf{u}^{(1)})+a_{22}L(\mathbf{u}^{(2)})),\\
%     \mathbf{u}^{(3)} = & \mathbf{u}^{n} + \Delta t(a_{31}L(\mathbf{u}^{(1)})+a_{32}L(\mathbf{u}^{(2)})+a_{33}L(\mathbf{u}^{(3)})),\\
%     \mathbf{u}^{(4)} = & \mathbf{u}^{n} + \Delta t(a_{41}L(\mathbf{u}^{(1)})+a_{42}L(\mathbf{u}^{(2)})+a_{43}L(\mathbf{u}^{(3)})+a_{44}L(\mathbf{u}^{(4)})),\\
%     \mathbf{u}^{n+1} = & \mathbf{u}^{n} + \Delta t(a_{51}L(\mathbf{u}^{(1)})+a_{52}L(\mathbf{u}^{(2)})+a_{53}L(\mathbf{u}^{(3)})+a_{54}L(\mathbf{u}^{(4)})),\\
% \end{split}
% \end{equation}
where the coefficients are given as
\begin{equation*}
\begin{split}
&a_{11}=0.1464466094067262,\\
&a_{21}=0.3535533905932738, a_{22}=0.1464466094067262,\\
&a_{31}=0.3535533905932738, a_{32}=0.3535533905932738, a_{33}=0.1464466094067262,\\
&a_{41}=0.33333333333333333, a_{42}=0.33333333333333333, a_{43}=0.33333333333333333.
\end{split}
\end{equation*}
% \begin{equation*}
% \begin{split}
% &a_{11}=0.112701665379258,\\
% &a_{21}=0.258198889747161, a_{22}=0.112701665379258,\\
% &a_{31}=0.258198889747161, a_{32}=0.258198889747161, a_{33}=0.112701665379258,\\
% &a_{41}=0.258198889747161, a_{42}=0.258198889747161, a_{43}=0.258198889747161, a_{44}=0.112701665379258,\\
% &a_{51}=0.250000000000000, a_{52}=0.250000000000000, a_{53}=0.250000000000000, a_{54}=0.250000000000000.
% \end{split}
% \end{equation*}
\end{itemize}
As discussed in the previous sections, the high complexity and costs in the calculation of Jacobian matrices for high-order WENO spatial discretizations lead to significant challenges in applying a fully implicit time-stepping method, which is also a motivation for us to design a new semilinearization technique for the exponential integrators in this paper. To simplify the implementation in applying the SSP-IRK3 method for the numerical comparison studies here, we couple the SSP-IRK3 method with the linear spatial discretization scheme \eqref{eq:DiffODEsLinear}.
The Newton iteration is used to solve the nonlinear algebraic systems at every time step of the SSP-IRK3 scheme.
Note that appropriate preconditioning techniques could speed up the computation of an implicit method such as the SSP-IRK3 method. Here, we do not further explore more sophisticated nonlinear system solvers for the implicit methods, but leave this interesting topic in our future work.
For more detailed introduction of these SSP-RK schemes, we refer the readers to the monograph \cite{gottlieb2011strong}.
In the following section, The numerical results obtained from these different time-stepping approaches are presented to compare their efficiency and accuracy.

\section{Numerical experiments}\label{sect:3}

This section is devoted to the numerical tests of the ETD-RK multi-resolution A-WENO schemes developed in the previous sections. 
All numerical examples are adopted from the published benchmarks that are widely used in the literature, e.g. \cite{KurganovT,liu2011high,arbogast2019finite,jiang2021high,zhang2022high,vijaywargiya2023two}. 

We present the examples in order of increasing complexity and compare the results of the ETD-RK methods with those obtained using the explicit SSP-RK (SSP-ERK3) and the diagonally implicit SSP-RK (SSP-IRK3) methods.
For the ETD-RK and the implicit SSP-RK methods, we adopt the time-step size $\Delta t = \text{CFL} \times \Delta x$ for one-dimensional problems and $\Delta t = \text{CFL} \times \min{\{\Delta x, \Delta y\}}$ for two-dimensional problems, respectively, where CFL is a constant that may take different values in different examples.
For the SSP-ERK3 method, we follow the CFL stability condition of explicit methods and take $\Delta t = \frac{0.4}{c/\Delta x + b/\Delta x^2}$ for one-dimensional problems, where $c = \max_{u}|f'(u)|$ and $b = \max_{u}|g'(u)|$; and $\Delta t = \frac{0.4}{c_x/\Delta x + b_x/\Delta x^2 + c_y/\Delta y + b_y/\Delta y^2}$ for two-dimensional problems, where $c_x = \max_{u}|f_1'(u)|$, $b_x = \max_{u}|g_1'(u)|$, $c_y = \max_{u}|f_2'(u)|$, and $b_y = \max_{u}|g_2'(u)|$.
All computations in this paper are conducted using Matlab R2023b on an Apple M2 Pro chip 12-core CPU with 16 GB of RAM.
In the following discussions of numerical results, $N$ and $M$ denote the number of grid points in the $x-$ and the $y-$ direction of a computational grid respectively. 
%For a straightforward comparison, we present the ratio $\frac{\Delta t}{\Delta x}$ for all methods.
%This ratio will be fixed for ETD and implicit methods, but it will change for the explicit method with time and we present its average throughout the simulation.

%%%%%%%%%%%%%%%%%%%%%%%%%%%%%%%%%%%%%%%%%%%%%%%%%%%%%%%%%%%%%%%%%%%%
\subsection{One-dimensional (1D) examples}

\begin{exmp}\label{ex:ex1}
\textbf{1D heat equation.}

We solve the one-dimensional heat equation
\begin{equation}\label{eq:heatEq1D}
u_t=u_{xx},
\end{equation}
on the domain $\Omega=[-\pi, \pi]$ with the periodic boundary condition and the initial condition $u(x,0)=\sin(x)$.
The problem has the exact solution $u(x,t)=e^{-t}\sin(x)$.

The solution is computed up to $T=1$ using the ETD-RK multi-resolution A-WENO methods of different temporal and spatial orders with the time-step size $\Delta t=\Delta x$.
The numerical errors and orders of convergence are presented in Table \ref{tab:ex1_1}.
The desired accuracy orders of the spatial discretizations, instead of the temporal discretizations, are observed in the table. The results verify that for a pure diffusion problem with a smooth solution, the numerical errors of the temporal discretizations using these exponential integrators for the diffusion term are very small. 
In addition, we compare the computational efficiency of the ETD-RK3 scheme with
the SSP-ERK3 scheme and the SSP-IRK3 scheme, coupled with the fourth-order spatial discretizations (i.e., the $r=2$ case in the section 2.1.1). Specifically, the SSP-ERK3 scheme is coupled with the fourth-order multi-resolution A-WENO discretization, while the fourth-order linear spatial discretization is used for the SSP-IRK3 scheme as that discussed in the section 2.2.4. For the ETD-RK3 scheme, both the A-WENO and the linear spatial discretizations are applied for the comparison.   
The computation is conducted on the grids with $N=20, 40, 60,\ldots, 140$.
On all grids, we take $\Delta t=\Delta x$ for the SSP-IRK3 scheme as that for the ETD-RK3 scheme.
The $L^1$ errors versus the CPU times for different methods are shown in Figure \ref{fig:efficiency_ex1}.
In the figure, we observe the superiority of the efficiency of the ETD-RK method compared to both the explicit and the implicit SSP-RK methods. It takes less CPU time costs for the ETD-RK3 scheme than the other two schemes to achieve a similar level of small numerical errors on refined meshes.  

\begin{table}[!htbp]
\centering
\begin{tabular}{ccccccccc}
\toprule[1.5pt]
\multicolumn{8}{c}{ETD-RK3} \\
\midrule
& MRWENO4 & & MRWENO6 & & MRWENO8\\
\cline{2-5} \cline{6-9}
$N$ & $L^1$ Error & Order & $L^1$ Error & Order & $L^1$ Error & Order\\
\midrule
20 & $1.58\times10^{-4}$ & - & $2.48\times10^{-6}$ & - & $4.31\times10^{-8}$ & -\\
40 & $9.92\times10^{-6}$ & $3.99$ & $3.93\times10^{-8}$ & $5.98$ & $1.72\times10^{-10}$ & $7.97$\\
60 & $1.96\times10^{-6}$ & $4.00$ & $3.46\times10^{-9}$ & $5.99$ & $6.74\times10^{-12}$ & $7.99$ \\
80 & $6.22\times10^{-7}$ & $4.00$ & $6.16\times10^{-10}$ & $6.00$ & $6.60\times10^{-13}$ & $8.07$\\
100 & $2.55\times10^{-7}$ & $4.00$ & $1.62\times10^{-10}$ & $6.00$ & $1.16\times10^{-13}$ & $7.80$\\
\midrule
\multicolumn{8}{c}{ETD-RK4} \\
\midrule
& MRWENO4 & & MRWENO6 & & MRWENO8\\
\cline{2-5} \cline{6-9}
$N$ & $L^1$ Error & Order & $L^1$ Error & Order & $L^1$ Error & Order & & \\
\midrule
20 & $1.58\times10^{-4}$ & - & $2.48\times10^{-6}$ & - & $4.31\times10^{-8}$ & -\\
40 & $9.92\times10^{-6}$ & $3.99$ & $3.93\times10^{-8}$ & $5.98$ & $1.72\times10^{-10}$ & $7.97$\\
60 & $1.96\times10^{-6}$ & $4.00$ & $3.46\times10^{-9}$ & $5.99$ & $6.74\times10^{-12}$ & $7.99$ \\
80 & $6.22\times10^{-7}$ & $4.00$ & $6.16\times10^{-10}$ & $6.00$ & $6.60\times10^{-13}$ & $8.08$\\
100 & $2.55\times10^{-7}$ & $4.00$ & $1.62\times10^{-10}$ & $6.00$ & $1.15\times10^{-13}$ & $7.84$\\
\bottomrule[1.5pt]
\end{tabular}
\caption{\textbf{Example \ref{ex:ex1}.} Numerical errors of the ETD-RK multi-resolution A-WENO methods with the time-step size $\Delta t=\Delta x$. 
MRWENO$2r$ stands for the $2r$-th order multi-resolution A-WENO discretization in space.}
\label{tab:ex1_1}
\end{table}

\begin{figure}[!htbp]
 \centering
  \includegraphics[width=0.5\textwidth]{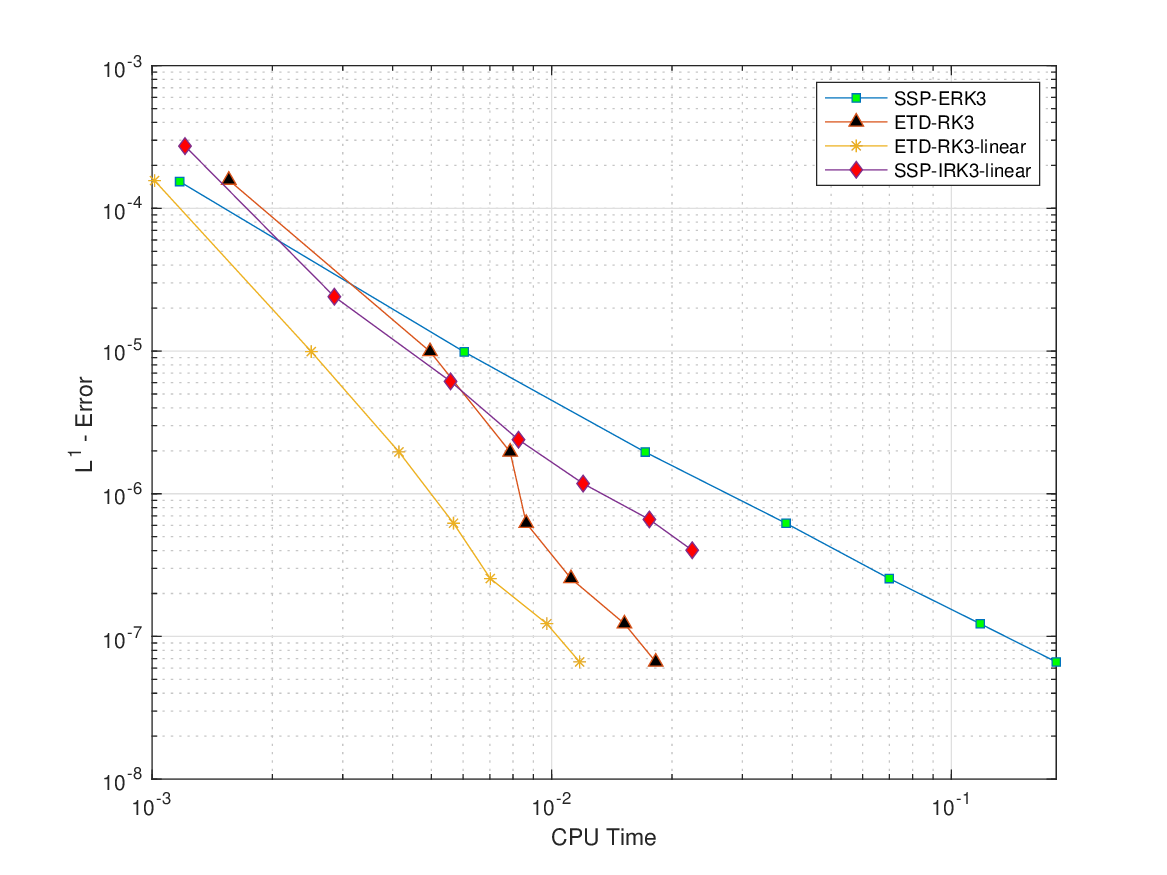}
 \caption{\textbf{Example \ref{ex:ex1}.} 
 Comparison of efficiency for different time-stepping methods. 
 SSP-IRK3-linear and ETD-RK3-linear indicate that the linear spatial discretization is used for them. CPU time unit: second.}
 \label{fig:efficiency_ex1}
\end{figure}

\end{exmp}

\begin{exmp}\label{ex:ex2}
\textbf{1D nonlinear stiff reaction-diffusion equation.}

We solve the fully nonlinear stiff reaction-diffusion equation 
\begin{equation}\label{eq:ex2}
    u_t=128(u^8)_{xx}+R(u),
\end{equation}
where the reaction term $R(u)=\frac{1}{1024u^7}-\frac{u}{8}+128 u^8-1$, on the domain $\Omega=[-\pi,\pi]$ with the periodic boundary condition and the initial condition $u(x,0)=\frac12(\sin(x)+2)^{\frac18}$. 
The problem has the exact solution $u(x,t)=\frac12(e^{-t}\sin(x)+2)^{\frac18}$.

The solution is computed up to $T=1$ using the ETD-RK multi-resolution A-WENO methods of different temporal and spatial orders with the time-step size $\Delta t=0.01\Delta x$.
Note that due to these complex nonlinear stiff diffusion and reaction terms, a smaller CFL number is required in this example than Example 1. However, time-step
size $\Delta t \sim O (\Delta x)$ can still be preserved in the mesh refinement study.    
The numerical errors and orders of convergence are presented in Table \ref{tab:ex2_1}. We observe that for the fourth-order multi-resolution A-WENO scheme coupled with either the ETD-RK3 or the ETD-RK4 temporal discretizations, the numerical errors of the spatial discretization dominate and a fourth-order accuracy is obtained.  
However, for the sixth-order and the eighth-order multi-resolution A-WENO schemes coupled with either the ETD-RK3 or the ETD-RK4 temporal discretizations, the numerical errors of the temporal discretizations dominate along with the refinement of the meshes, hence the desired accuracy orders of the temporal discretizations, instead of the spatial discretizations, are observed in the table.
Comparing with Example 1 which is a linear problem with a pure diffusion term, in this example the equation has a highly nonlinear and complex reaction term in addition to a stiff nonlinear diffusion term, and numerical errors of the spatial discretizations and the temporal discretizations from the different terms have a richer structure.

In addition, we compare the computational efficiency of the ETD-RK3 method with the SSP-ERK3 and the SSP-IRK3 methods, coupled with the fourth-order spatial discretization (i.e., the $r=2$ case in the section 2.1.1) on the grids with $N=50,100,150,\ldots, 350$. Similar to Example 1, the SSP-ERK3 scheme is coupled with
the fourth-order multi-resolution A-WENO discretization, while the fourth-order linear spatial
discretization is used for the SSP-IRK3 scheme as that discussed in the section 2.2.4. For the
ETD-RK3 scheme, both the A-WENO and the linear spatial discretizations are applied in the
comparison.
The time-step size $\Delta t=0.02\Delta x$ is taken for the simulations using the ETD-RK3 method on all grids, while for the SSP-IRK3 method, we take $\Delta t=0.02\Delta x$ for the simulations on the grids with  $N=50, 100$, $\Delta t=0.01\Delta x$ for the simulations on the grids with $N=150,200$, and $\Delta t=0.005\Delta x$ for the simulations on the other finer grids to ensure the stability in solving this stiff nonlinear problem.
The $L^1$ errors versus the CPU times for different methods are shown in Figure \ref{fig:efficiency_ex2}, which demonstrates that the ETD-RK3 method is more efficient than both the explicit and the implicit SSP-RK methods here. It takes less CPU time costs for the ETD-RK3 method than the other two methods to achieve a similar level of small numerical errors.

\begin{table}[!htbp]
\centering
\begin{tabular}{ccccccccc}
\toprule[1.5pt]
\multicolumn{8}{c}{ETD-RK3} \\
\midrule
& MRWENO4 & & MRWENO6 & & MRWENO8\\
\cline{2-5} \cline{6-9}
$N$ & $L^1$ Error & Order & $L^1$ Error & Order & $L^1$ Error & Order\\
\midrule
50 & $1.19\times10^{-1}$ & - & $2.78\times10^{-3}$ & - & $2.52\times10^{-3}$ & -\\
100 & $6.74\times10^{-3}$ & $4.15$ & $3.22\times10^{-4}$ & $3.11$ & $3.18\times10^{-4}$ & $2.99$ \\
150 & $1.36\times10^{-3}$ & $3.96$ & $9.47\times10^{-5}$ & $3.01$ & $9.44\times10^{-5}$ & $2.99$ \\
200 & $4.38\times10^{-4}$ & $3.92$ & $3.99\times10^{-5}$ & $3.00$ & $3.99\times10^{-5}$ & $2.99$ \\
250 & $1.84\times10^{-4}$ & $3.90$ & $2.05\times10^{-5}$ & $3.00$ & $2.04\times10^{-5}$ & $3.00$ \\
\midrule
\multicolumn{8}{c}{ETD-RK4} \\
\midrule
& MRWENO4 & & MRWENO6 & & MRWENO8\\
\cline{2-5} \cline{6-9}
$N$ & $L^1$ Error & Order & $L^1$ Error & Order & $L^1$ Error & Order & & \\
\midrule
50 & $1.16\times10^{-1}$ & - & $2.70\times10^{-4}$ & - & $1.25\times10^{-5}$ & -\\
100 & $6.42\times10^{-3}$ & $4.18$ & $4.78\times10^{-6}$ & $5.82$ & $7.43\times10^{-7}$ & $4.07$ \\
150 & $1.26\times10^{-3}$ & $4.01$ & $5.03\times10^{-7}$ & $5.56$ & $1.44\times10^{-7}$ & $4.05$ \\
200 & $3.98\times10^{-4}$ & $4.00$ & $1.01\times10^{-7}$ & $5.56$ & $4.48\times10^{-8}$ & $4.05$ \\
250 & $1.63\times10^{-4}$ & $4.00$ & $4.19\times10^{-8}$ & $3.96$ & $1.91\times10^{-8}$ & $3.83$ \\
\bottomrule[1.5pt]
\end{tabular}
\caption{\textbf{Example \ref{ex:ex2}.}
Numerical errors of the ETD-RK multi-resolution A-WENO methods with the time-step size $\Delta t=0.01\Delta x$.
MRWENO$2r$ stands for the $2r$-th order multi-resolution A-WENO discretization in space.}
\label{tab:ex2_1}
\end{table}

\begin{figure}[!htbp]
 \centering
  \includegraphics[width=0.5\textwidth]{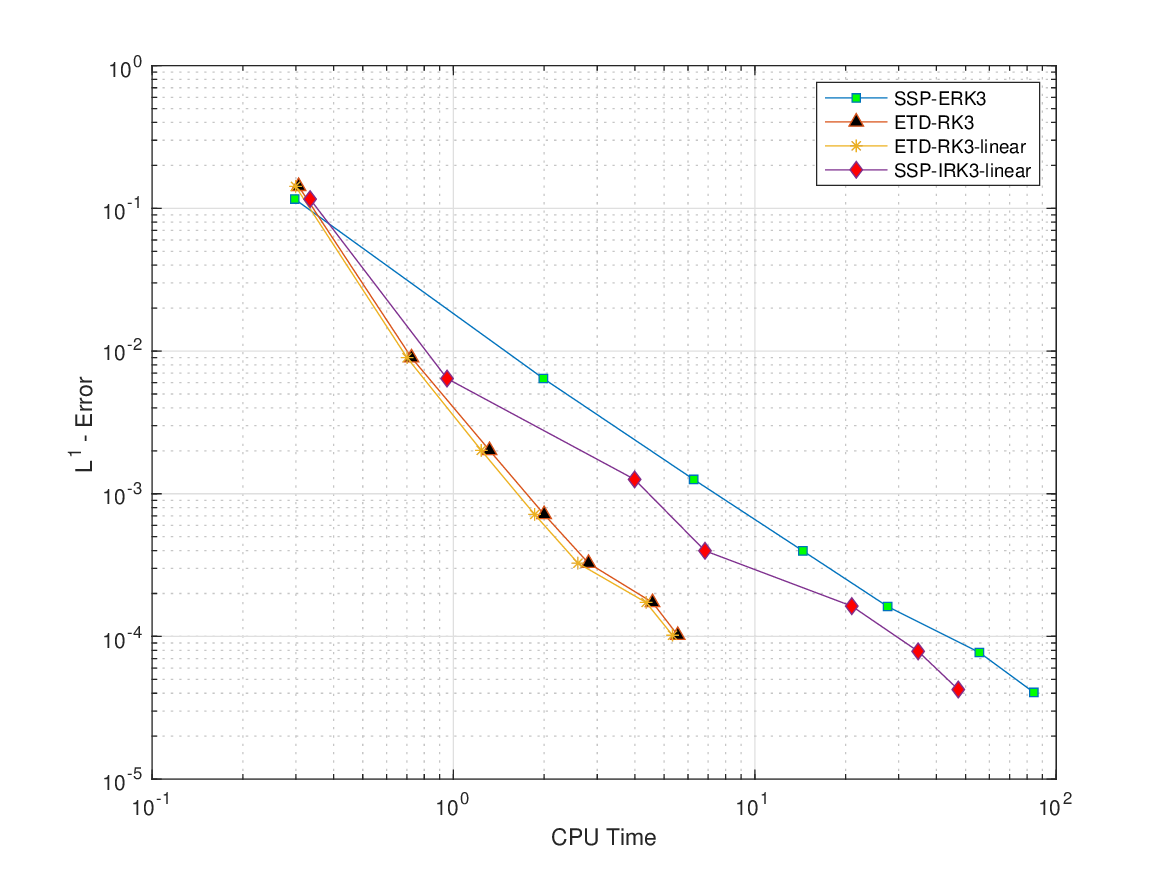}
 \caption{\textbf{Example \ref{ex:ex2}.}
 Comparison of efficiency for different time-stepping methods. 
 SSP-IRK3-linear and ETD-RK3-linear indicate that the linear
 spatial discretization is used for them. CPU time unit: second.}
 \label{fig:efficiency_ex2}
\end{figure}

\end{exmp}

\begin{exmp}\label{ex:ex3}
\textbf{1D PME with the Barenblatt solution.}

We solve the 1D PME \eqref{eq:PME1D} with the Barenblatt solution \eqref{eq:Barenblatt} and the homogeneous Dirichlet boundary condition.
The solution at $t_0=1$ serves as the initial condition.

The solutions with different values of $m$ are computed up to $T=2$ on the domain $\Omega=[-6, 6]$. The computational grid with $N=200$ is used. The time-step size is taken as $\Delta t = \Delta x$.
The numerical solutions obtained from the simulations using the ETD-RK4 scheme coupled with the sixth-order multi-resolution A-WENO (A-WENO6) spatial discretization are shown in Figure \ref{fig:ex3}.
Numerical results with other spatial and temporal order accuracy exhibit a similar pattern and are therefore not shown here to save space.
From the numerical results, we observe that the sharp wave fronts are captured stably and the numerical solutions match the exact solutions very well.

Then, we solve the equation with different values of $m$  up to $T=10$, applying different third-order time-marching approaches and comparing their computational efficiency. The computational grid with $N=300$ is used. The simulations are performed on the domain $\Omega=[-9, 9]$ to include the entire non-zero profile of the solutions. 
The time-step sizes used in different methods are chosen to be at their maximum values to achieve stable computations and  numerical solutions that approximate well the exact solutions. 
The comparison of the ratios of the time-step sizes to the spatial grid size, $\frac{\Delta t}{\Delta x}$, and the corresponding CPU times of different methods are shown in Table \ref{tab:ex3}. 
From the numerical results in the table, we observe that the ETD-RK3 method, which is coupled with either the multi-resolution A-WENO spatial discretizations or the corresponding linear spatial discretizations, allows for much larger time-step sizes and takes much less CPU time costs than both the explicit and the implicit SSP-RK methods here.   
Moreover, the permitted maximum time-step sizes have very small changes as the complexity of the equation and the accuracy order of spatial discretization increase, which shows the robustness of the ETD-RK method. Again, similar to the previous examples, the high computational efficiency of the ETD-RK method is verified here.  

\begin{figure}[!htbp]
 \centering
 \begin{subfigure}[b]{0.45\textwidth}
  \includegraphics[width=\textwidth]{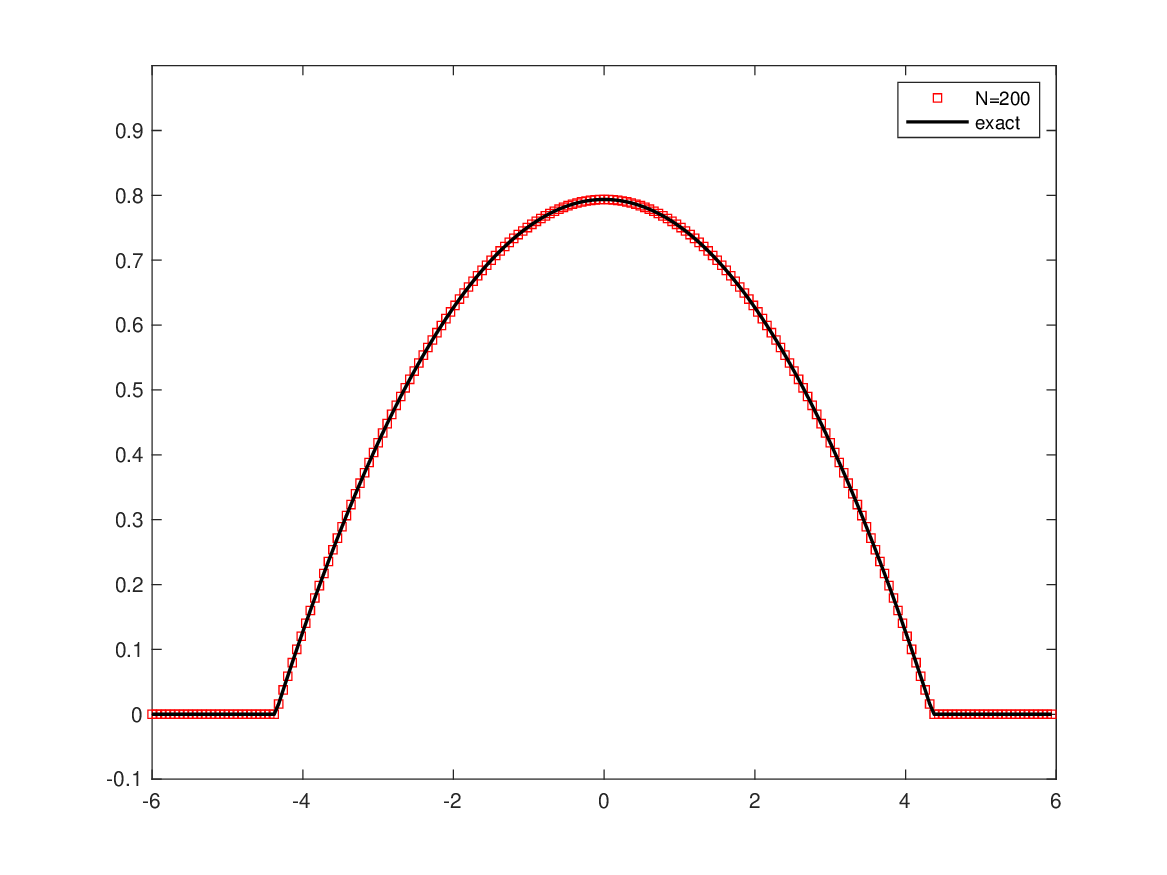}
  \caption{$m=2$.}
 \end{subfigure}
 \begin{subfigure}[b]{0.45\textwidth}
  \includegraphics[width=\textwidth]{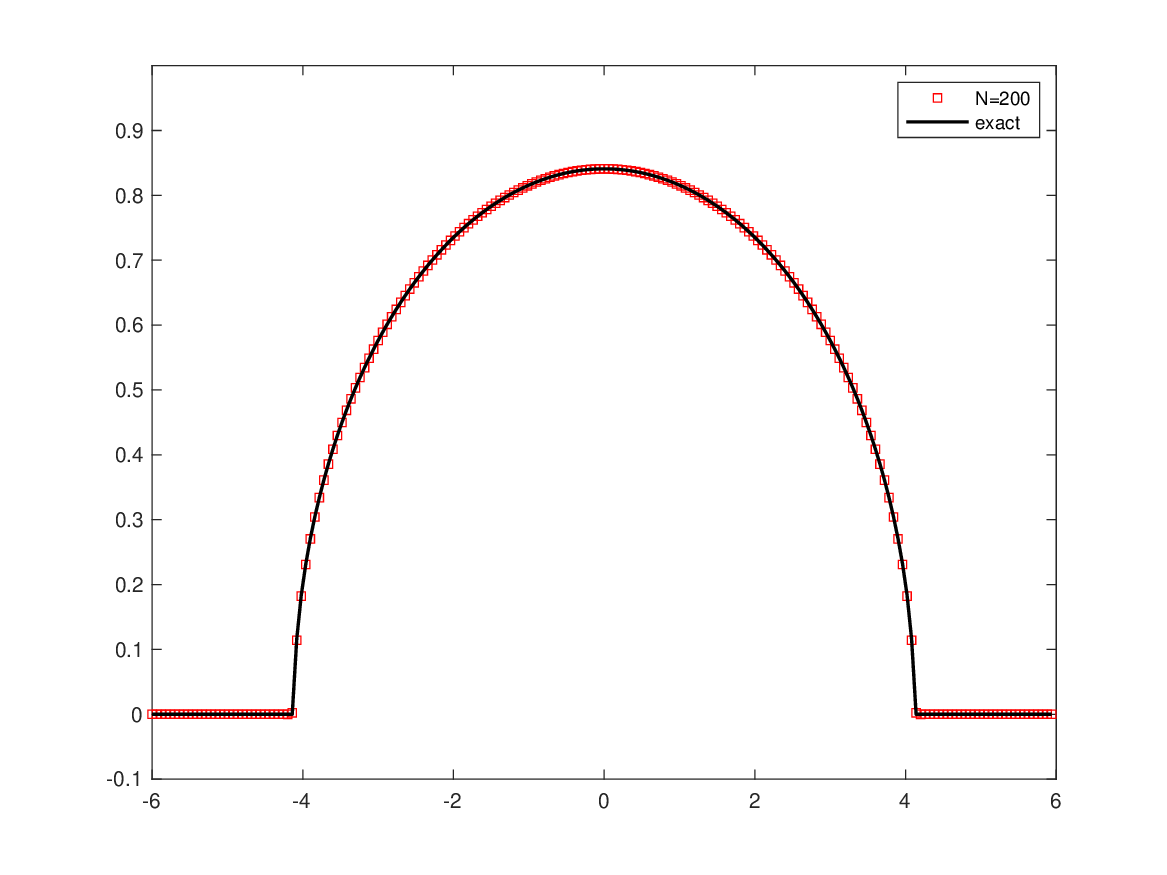}
  \caption{$m=3$.}
 \end{subfigure}
 \begin{subfigure}[b]{0.45\textwidth}
  \includegraphics[width=\textwidth]{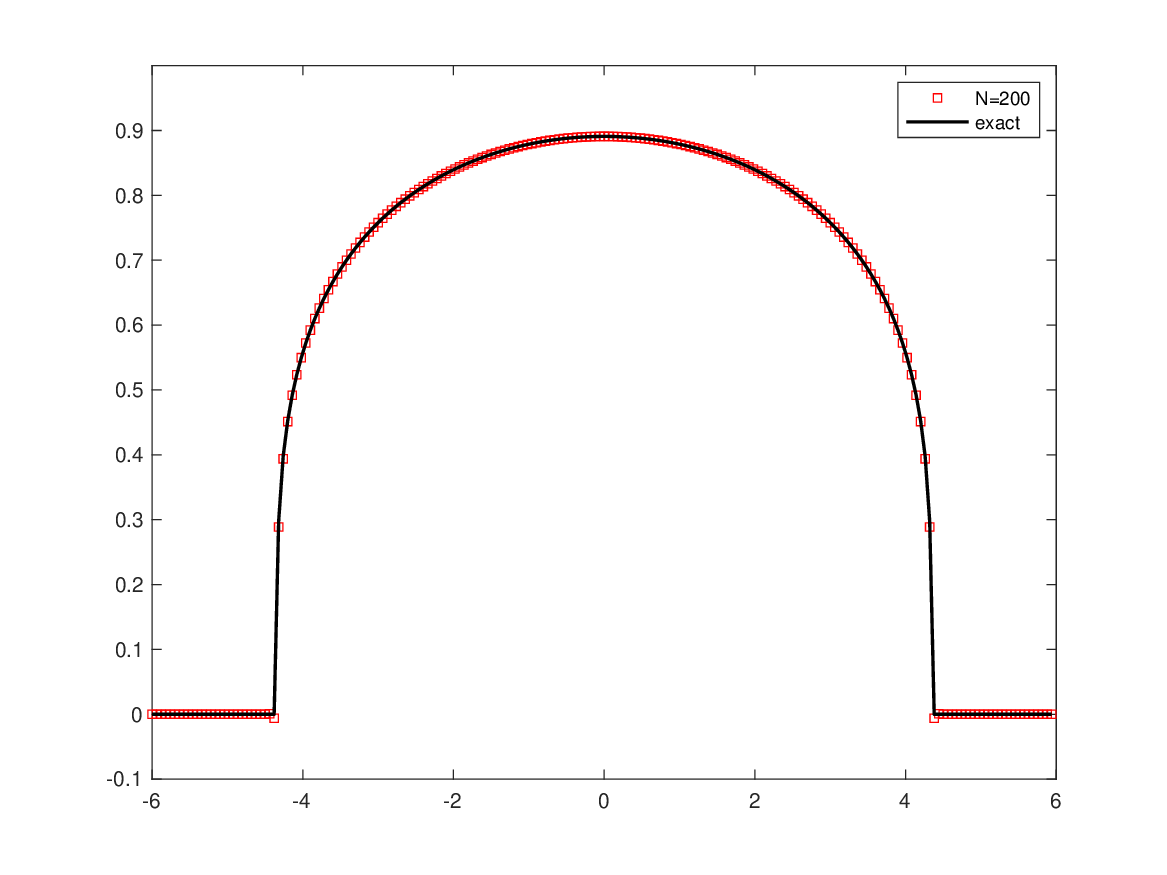}
  \caption{$m=5$.}
 \end{subfigure}
 \begin{subfigure}[b]{0.45\textwidth}
  \includegraphics[width=\textwidth]{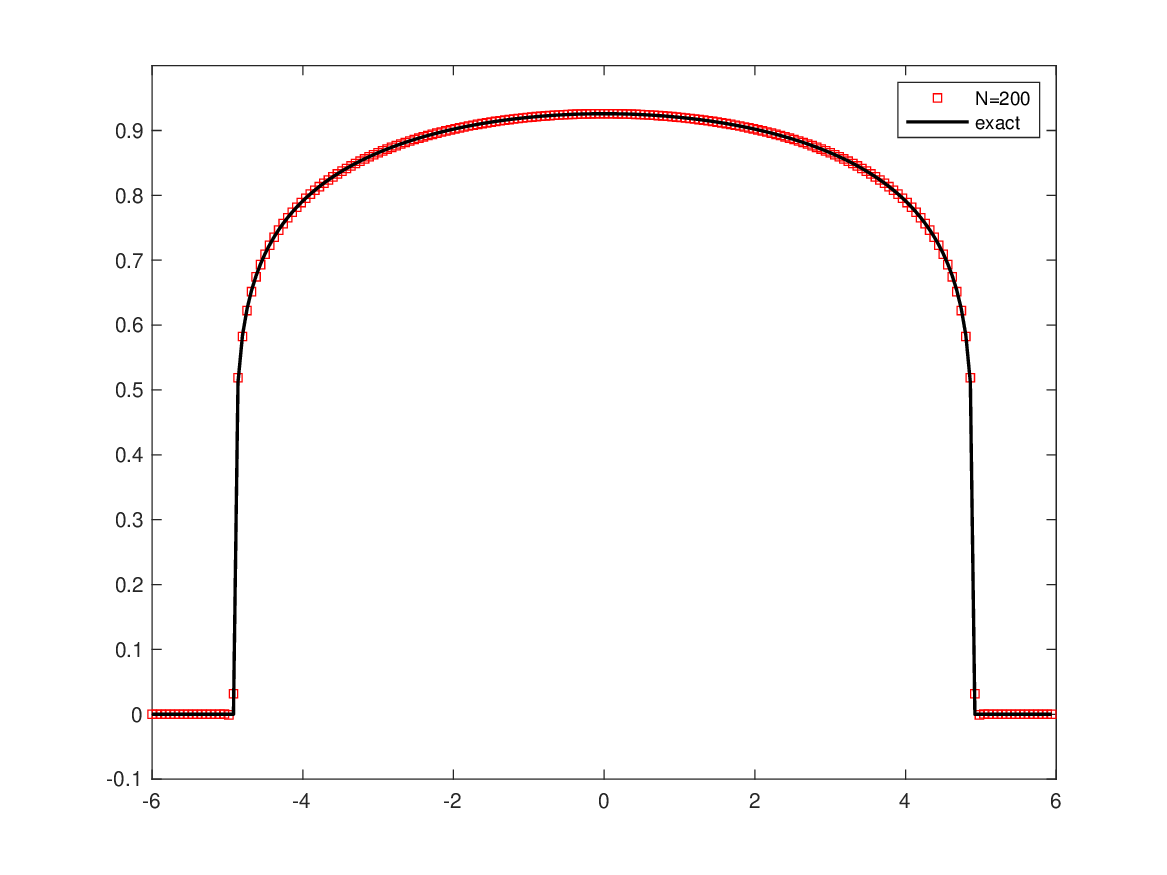}
  \caption{$m=8$.}
 \end{subfigure}
 \caption{\textbf{Example \ref{ex:ex3}.} 
 Numerical solutions of the ETD-RK4 multi-resolution A-WENO6 method on the grid with $N=200$, for the 1D PME with different values of $m$ for the Barenblatt solution at $T=2$. The time-step size is $\Delta t=\Delta x$.}
 \label{fig:ex3}
\end{figure}

\begin{table}[!htbp]
\centering
\begin{tabular}{cccccccccccc}
\toprule[1.5pt]
\multicolumn{7}{c}{$m=2$} \\
\midrule
Time & MRWENO4 & & MRWENO6 & & MRWENO8\\
\cline{2-3} \cline{4-5} \cline{6-7}
marching & $\Delta t/\Delta x$  & CPU (s) & $\Delta t/\Delta x$  & CPU (s) & $\Delta t/\Delta x$  & CPU (s)\\
\midrule
ETD-RK3 & $1.6$ & $0.27$ & $1.6$ & $0.31$ & $1.5$ & $0.35$ \\
SSP-ERK3 & $0.0198$ & $2.13$ & $0.0198$ & $1.93$ & $0.0198$ & $2.82$ \\
ETD-RK3-linear & $1.6$ & $0.26$ & $1.5$ & $0.30$ & $1.4$ & $0.33$ \\
SSP-IRK3-linear & $0.3$ & $1.28$ & $0.3$ & $1.41$ & $0.2$ & $2.31$\\
\midrule
\multicolumn{7}{c}{$m=3$}\\
\midrule
Time & MRWENO4 & & MRWENO6 & & MRWENO8\\
\cline{2-3} \cline{4-5} \cline{6-7}
marching & $\Delta t/\Delta x$  & CPU (s) & $\Delta t/\Delta x$  & CPU (s) & $\Delta t/\Delta x$  & CPU (s)\\
\midrule
ETD-RK3 & $1.5$ & $0.34$ & $1.4$ & $0.38$ & $1.5$ & $0.42$ \\
SSP-ERK3 & $0.0166$ & $2.72$ & $0.0166$ & $2.54$ & $0.0166$ & $3.53$ \\
ETD-RK3-linear & $1.5$ & $0.32$ & $1.4$ & $0.38$ & $1.5$ & $0.41$ \\
SSP-IRK3-linear & $0.2$ & $2.00$ & $0.2$ & $2.18$ & $0.2$ & $2.38$ \\
\midrule
\multicolumn{7}{c}{$m=5$} \\
\midrule
Time & MRWENO4 & & MRWENO6 & & MRWENO8\\
\cline{2-3} \cline{4-5} \cline{6-7}
marching & $\Delta t/\Delta x$  & CPU (s) & $\Delta t/\Delta x$  & CPU (s) & $\Delta t/\Delta x$  & CPU (s)\\
\midrule
ETD-RK3 & $1.5$ & $0.42$ & $1.5$ & $0.47$ & $1.5$ & $0.51$ \\
SSP-ERK3 & $0.0125$ & $3.73$ & $0.0125$ & $3.47$ & $0.0125$ & $4.74$ \\
ETD-RK3-linear & $1.4$ & $0.40$ & $1.4$ & $0.47$ & $1.4$ & $0.51$ \\
SSP-IRK3-linear & $0.1$ & $3.83$ & $0.1$ & $4.14$ & $0.1$ & $4.60$ \\
\midrule
\multicolumn{7}{c}{$m=8$} \\
\midrule
Time & MRWENO4 & & MRWENO6 & & MRWENO8\\
\cline{2-3} \cline{4-5} \cline{6-7}
marching & $\Delta t/\Delta x$  & CPU (s) & $\Delta t/\Delta x$  & CPU (s) & $\Delta t/\Delta x$  & CPU (s)\\
\midrule
ETD-RK3 & $1.4$ & $0.51$ & $1.4$ & $0.57$ & $1.4$ & $0.62$ \\
SSP-ERK3 & $0.009$ & $5.25$ & $0.009$ & $4.90$ & $0.009$ & $6.65$ \\
ETD-RK3-linear & $1.4$ & $0.50$ & $1.3$ & $0.58$ & $1.3$ & $0.62$ \\
SSP-IRK3-linear & $0.09$ & $4.52$ & $0.07$ & $5.94$ & $0.06$ & $7.59$\\
\bottomrule[1.5pt]
\end{tabular}
\caption{
\textbf{Example \ref{ex:ex3}.}
Maximum ratios of the time-step sizes to the spatial grid size and the corresponding CPU times of different methods in the computation of the one-dimensional PME with the Barenblatt solutions. MRWENO$2r$ stands for the $2r$-th order multi-resolution A-WENO discretization in space.
SSP-IRK3-linear and ETD-RK3-linear indicate that the corresponding linear spatial discretizations of the MRWENO schemes are used.}
\label{tab:ex3}
\end{table}

\end{exmp}

\begin{exmp}\label{ex:ex4}
\textbf{Interacting boxes.}

We solve the PME \eqref{eq:PME1D} with $m=6$ for  the interaction of two boxes depicted by the initial condition 
\begin{equation}
u(x,0)=
\begin{cases}
1,& -4<x<-1,\\ 
2,& 0<x<3,\\
0,& \text{otherwise},
\end{cases}
\end{equation}
on the domain $\Omega=[-6,6]$ with the homogeneous Dirichlet boundary condition. This model
describes how temperature changes when two hot spots are suddenly placed in the domain.

The solution is computed up to $T=0.8$ using the ETD-RK4 scheme coupled with the sixth-order multi-resolution A-WENO spatial discretization. The computation is carried out on the grid with $N=160$ and the time-step size is taken as $\Delta t=0.01\Delta x$. The numerical results at different times are shown in Figure \ref{fig:ex4_ETDRK4}, from which we observe the boxes merging as time progresses. Similar to Example 3, the sharp
wave fronts of the solution are captured stably with high resolution in the simulation.

\begin{figure}[!htbp]
 \centering
 \begin{subfigure}[b]{0.3\textwidth}
  \includegraphics[width=\textwidth]{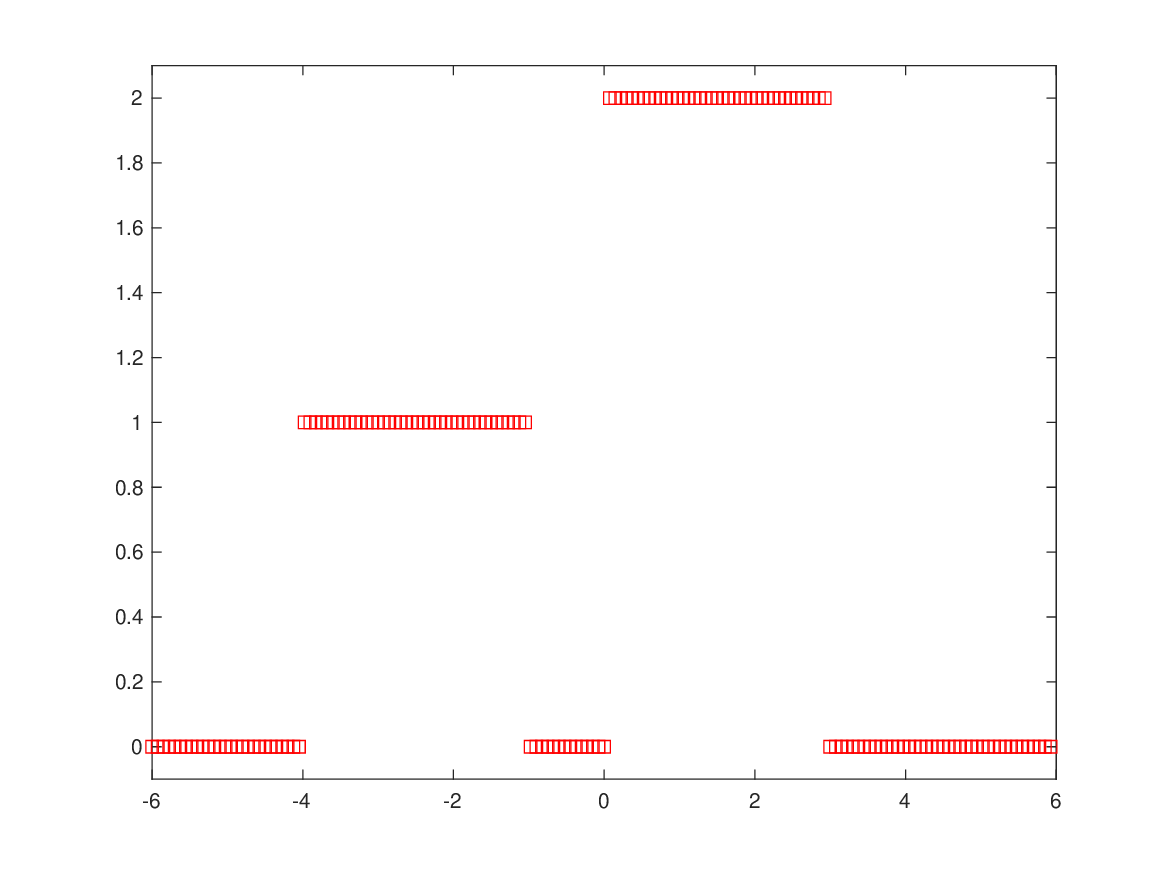}
  \caption{$t=0$}
 \end{subfigure}
 \begin{subfigure}[b]{0.3\textwidth}
  \includegraphics[width=\textwidth]{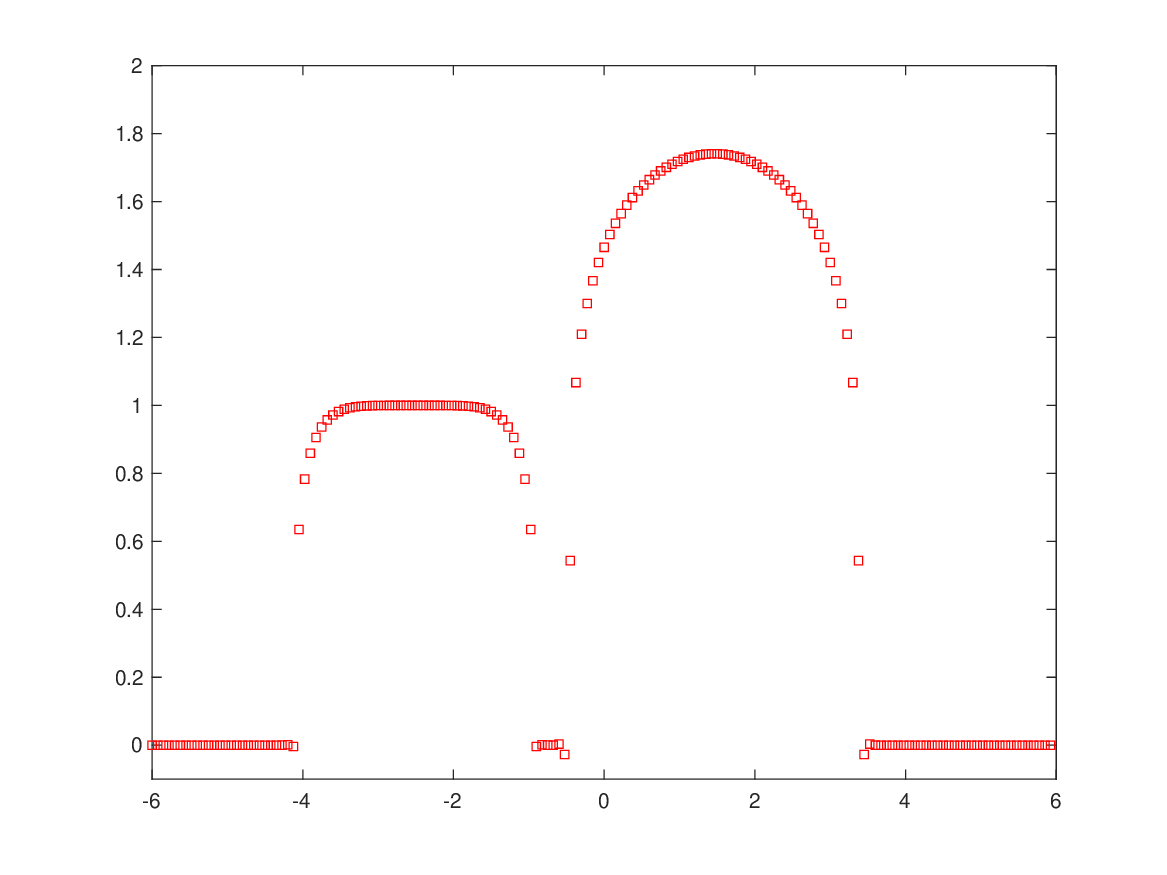}
  \caption{$t = 0.01$}
 \end{subfigure}

 \begin{subfigure}[b]{0.3\textwidth}
  \includegraphics[width=\textwidth]{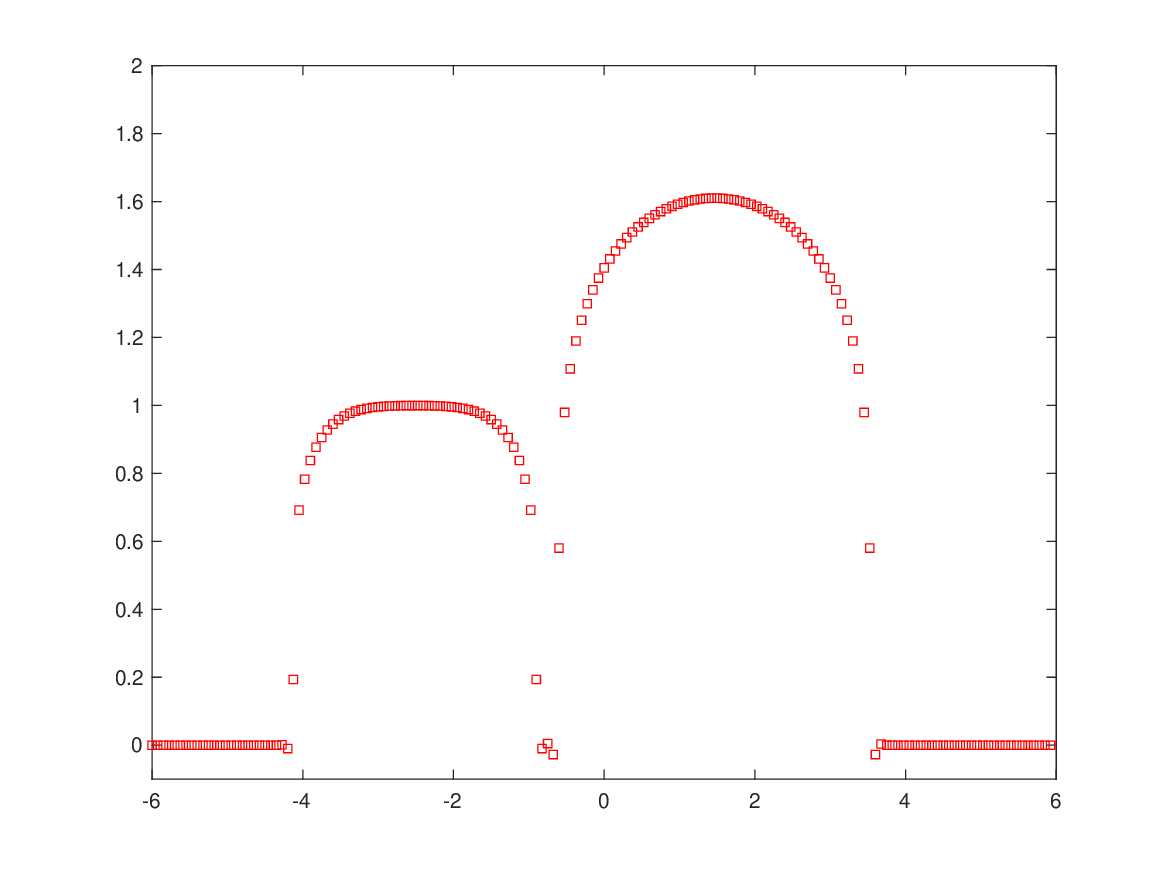}
  \caption{$t = 0.02$}
 \end{subfigure}
 \begin{subfigure}[b]{0.3\textwidth}
  \includegraphics[width=\textwidth]{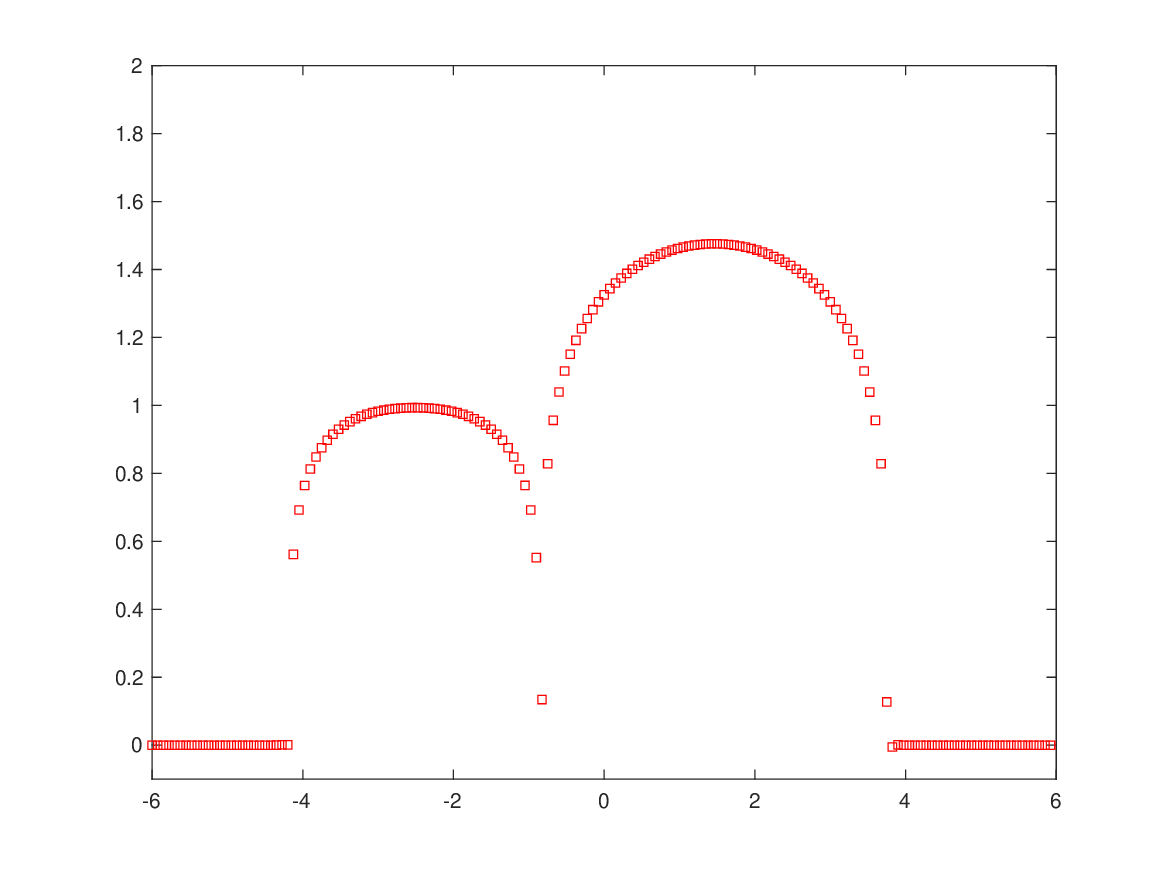}
  \caption{$t = 0.04$}
 \end{subfigure}

 \begin{subfigure}[b]{0.3\textwidth}
  \includegraphics[width=\textwidth]{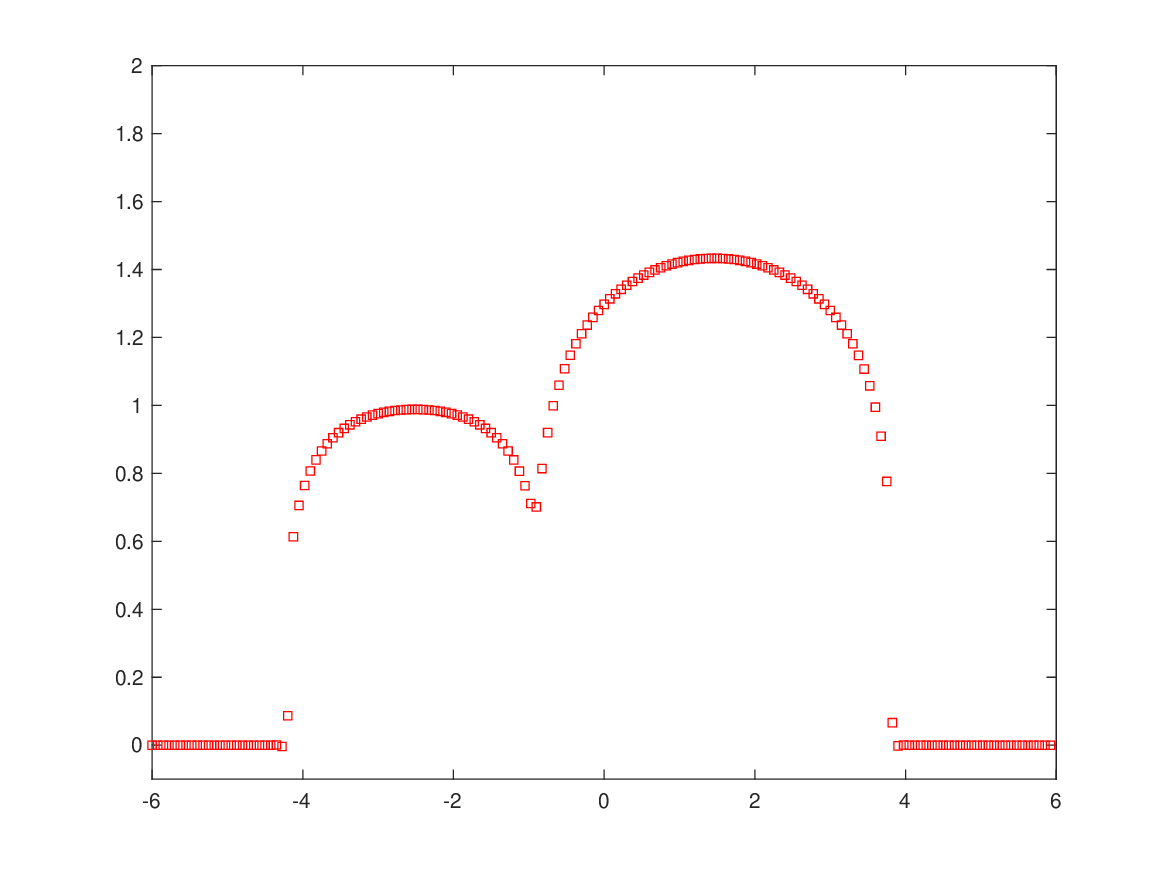}
  \caption{$t = 0.05$}
 \end{subfigure}
 \begin{subfigure}[b]{0.3\textwidth}
  \includegraphics[width=\textwidth]{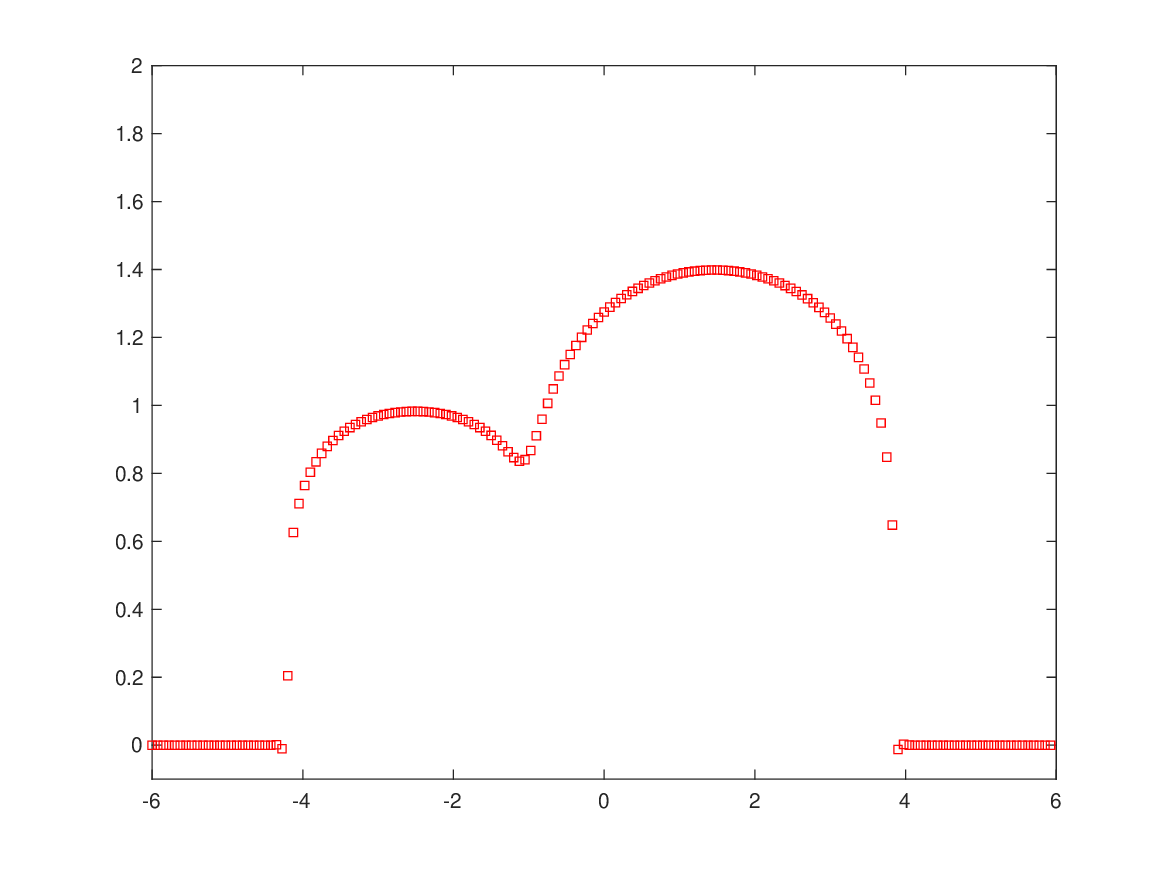}
  \caption{$t = 0.06$}
 \end{subfigure}

 \begin{subfigure}[b]{0.3\textwidth}
  \includegraphics[width=\textwidth]{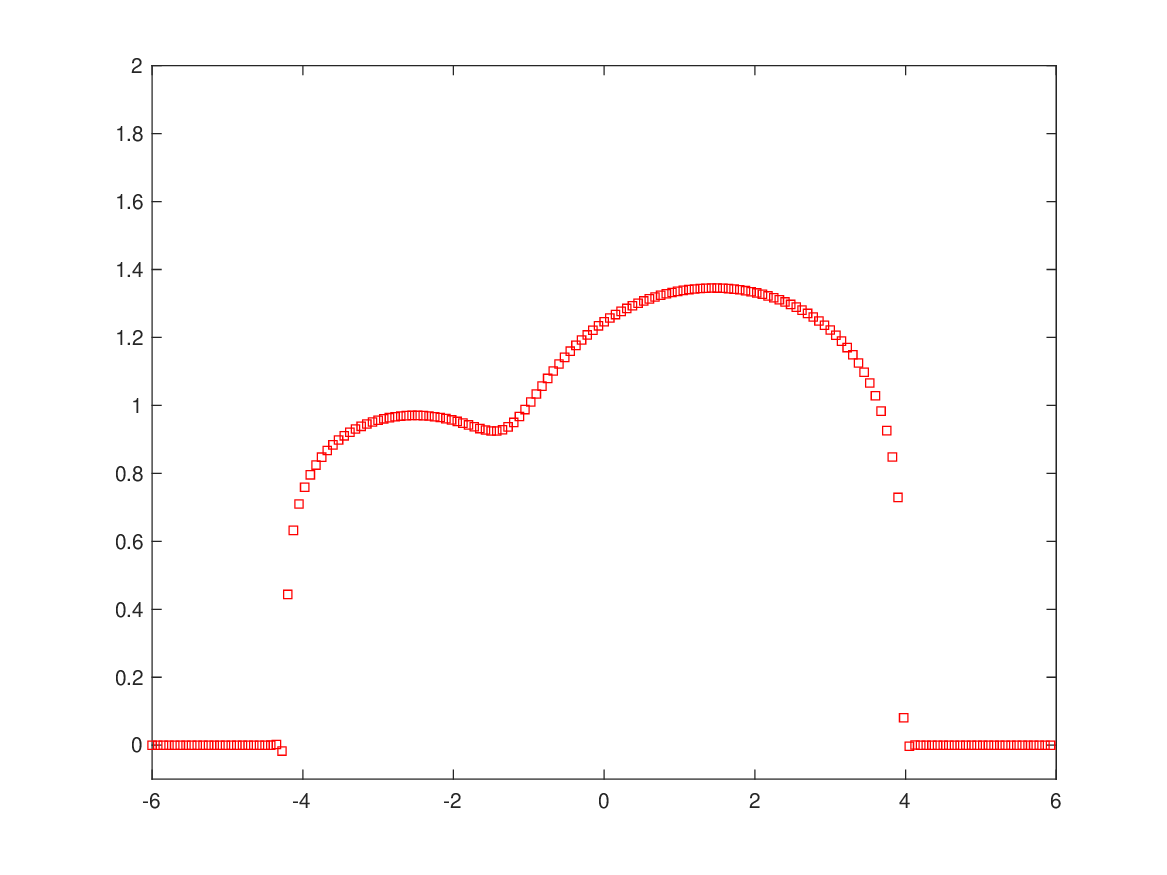}
  \caption{$t=0.08$}
 \end{subfigure}
 \begin{subfigure}[b]{0.3\textwidth}
  \includegraphics[width=\textwidth]{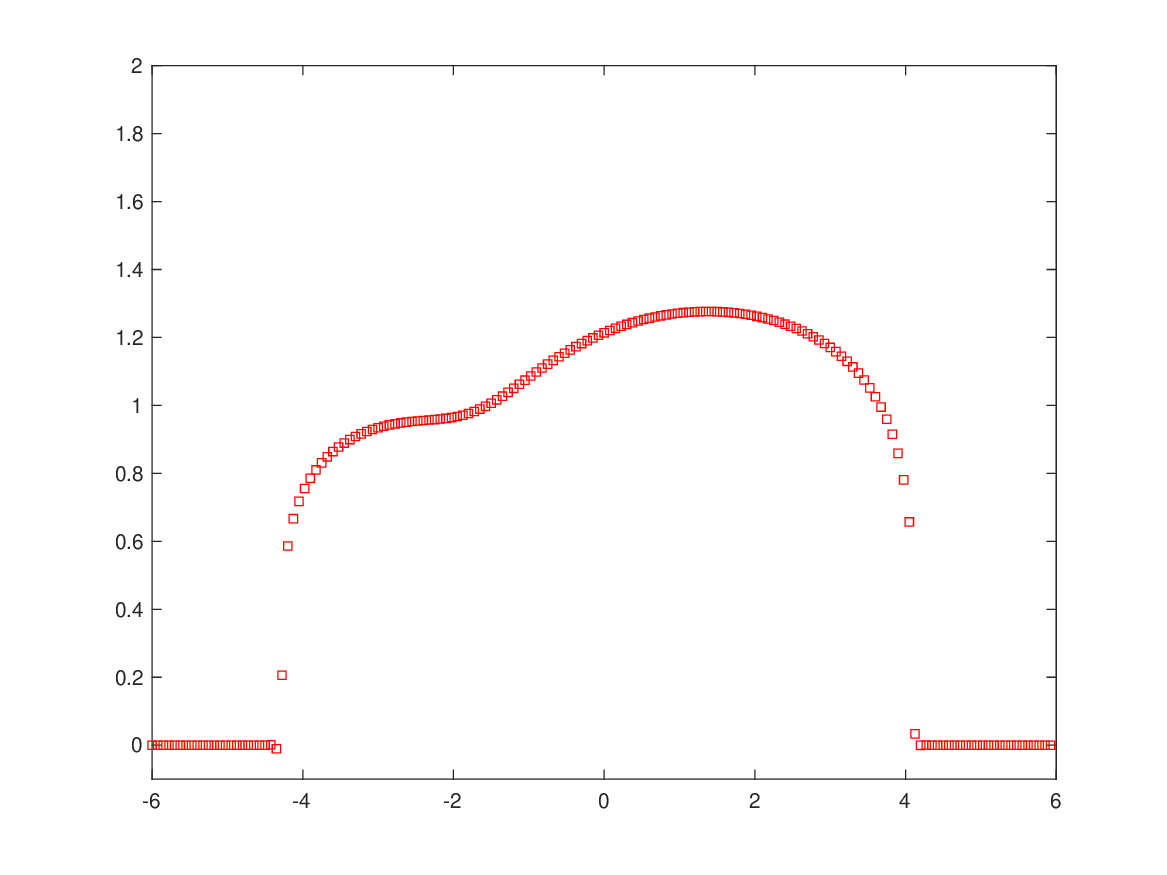}
  \caption{$t=0.12$}
 \end{subfigure}

 \begin{subfigure}[b]{0.3\textwidth}
  \includegraphics[width=\textwidth]{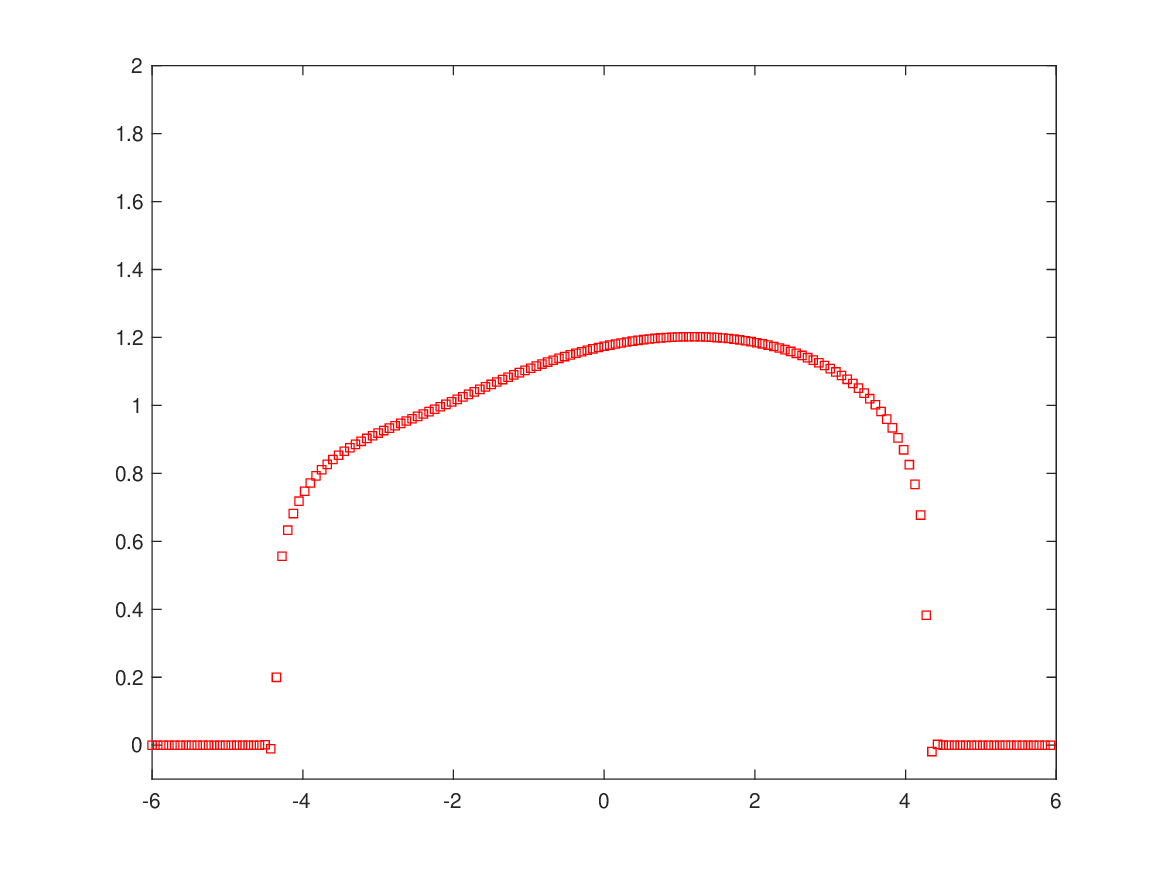}
  \caption{$t=0.2$}
 \end{subfigure}
 \begin{subfigure}[b]{0.3\textwidth}
  \includegraphics[width=\textwidth]{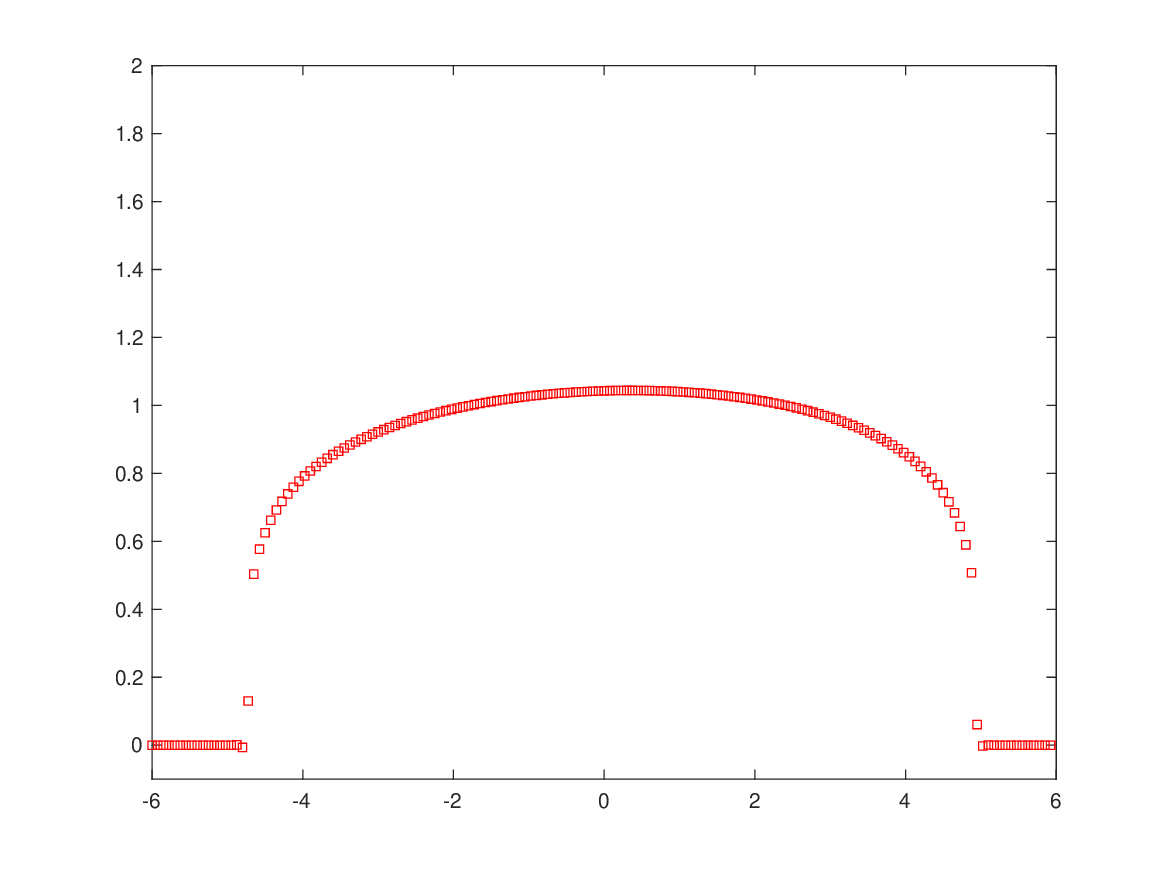}
  \caption{$t=0.8$}
 \end{subfigure}
 \caption{\textbf{Example \ref{ex:ex4}.} 
 Numerical solutions of the ETD-RK4 multi-resolution A-WENO6 method on the grid with $N=160$, for the problem of interaction of two boxes. The time-step size is $\Delta t=0.01\Delta x$.}
 \label{fig:ex4_ETDRK4}
\end{figure}

\end{exmp}

\begin{exmp}\label{ex:ex5}
\textbf{1D Buckley-Leverett equation.}

We consider the viscous Buckley-Leverett equation
\begin{equation}\label{eq:VBL}
u_t+f(u)_x=\epsilon(\nu(u)u_x)_{x},
\end{equation}
which is used to model the two-phase flow in porous media.
The nonlinear diffusion coefficient is taken as
\begin{equation}
\nu(u)=\begin{cases}
4u(1-u), &0\leq u\leq1,\\
0, &\text{otherwise,}
\end{cases}
\end{equation}
so that the function of the parabolic term $g(u)_{xx}=\epsilon(\nu(u)u_x)_x$ is given by
\begin{equation}
g(u)=
\begin{cases}
0, & u<0,\\
\epsilon(2u^2-\frac43 u^3),& 0\leq u\leq 1,\\
\frac23 \epsilon,& u>1.\\
\end{cases}
\end{equation}
Two different convection fluxes are considered, namely, the flux without gravitational effects given by
\begin{equation}\label{eq:BL1DNG}
f(u)=\frac{u^2}{u^2+(1-u)^2},
\end{equation}
and the one with gravitational effects given by
\begin{equation}\label{eq:BL1DGR}
f(u)=
\frac{u^2}{u^2+(1-u)^2}(1-5(1-u)^2).
\end{equation}

We take $\epsilon=0.01$, and solve two initial-boundary value problems. The first initial-boundary value problem defined on $\Omega=[0,1]$ has the initial condition
\begin{equation}\label{eq:ex5_IC1}
u(x,0)=
\begin{cases}
1-3x,& 0\leq x\leq \frac13,\\
0,& \frac{1}{3}<x\leq 1,
\end{cases}
\end{equation}
and the boundary condition $u(0,t)=1, u(1,t)=0$. The second problem is a Riemann problem defined on $\Omega=[0,1]$, which has the initial condition
\begin{equation}\label{eq:ex5_IC2}
u(x,0)=
\begin{cases}
0,& 0\leq x<1-\frac{1}{\sqrt{2}},\\
1, & 1-\frac{1}{\sqrt{2}}\leq x\leq1,
\end{cases}
\end{equation}
and the boundary condition $u(0,t)=0$ and $u(1,t)=1$.

The ETD-RK4 multi-resolution A-WENO6 scheme is used to solve these problems. The simulations are performed on the computational grids with $N=100$ and $N=800$. The time-step sizes are chosen to be the maximum permitted values for achieving a stable computation in each problem. The obtained numerical solutions at $T=0.2$ are shown in Figure \ref{fig:ex5}. We observe that the numerical solutions on different grids match very well, which verifies the convergence of the numerical solutions. Similar to the previous examples, the large gradients of the solutions are captured stably with high resolution in the simulations, which verifies the nonlinear stability of the method.
Moreover, the desired time-step
size $\Delta t \sim O (\Delta x)$ for the ETD-RK method is still preserved on different grids in this example.

\begin{figure}[!htbp]
 \centering
 \begin{subfigure}[b]{0.45\textwidth}
  \includegraphics[width=\textwidth]{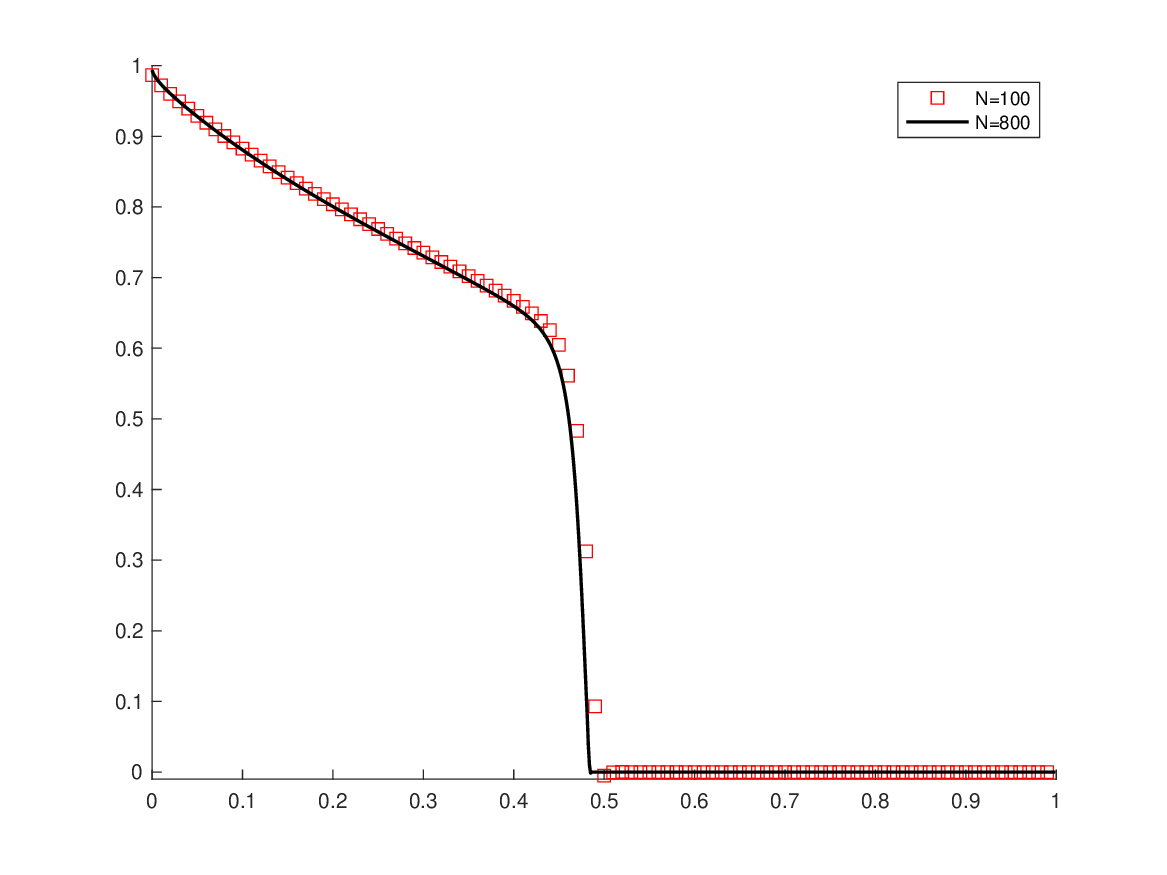}
  \caption{Initial-boundary value problem \eqref{eq:ex5_IC1} without gravity.
  Time-step size satisfies $\frac{\Delta t}{\Delta x}=0.9$.}
 \end{subfigure}
 \begin{subfigure}[b]{0.45\textwidth}
  \includegraphics[width=\textwidth]{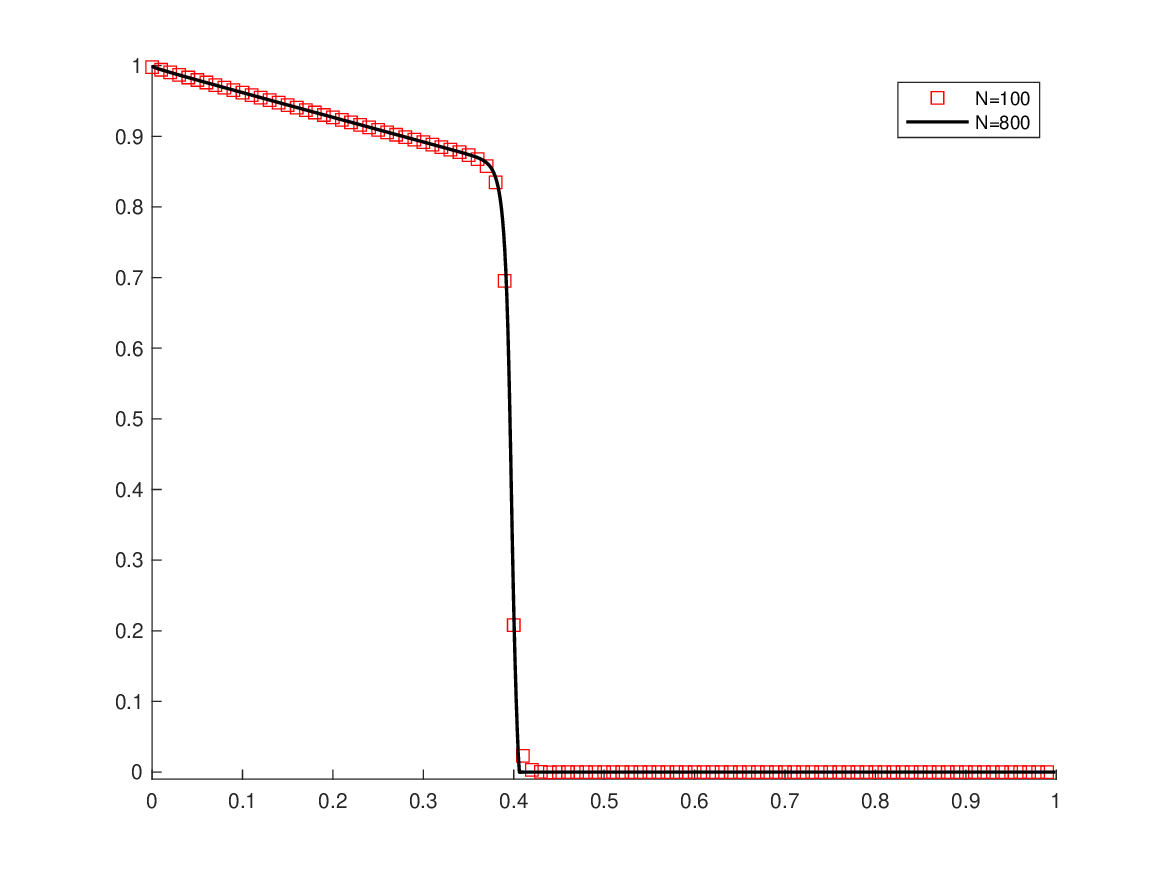}
  \caption{Initial-boundary value problem \eqref{eq:ex5_IC1} with gravity.
  Time-step size satisfies $\frac{\Delta t}{\Delta x}=0.5$.}
 \end{subfigure}
 \begin{subfigure}[b]{0.45\textwidth}
  \includegraphics[width=\textwidth]{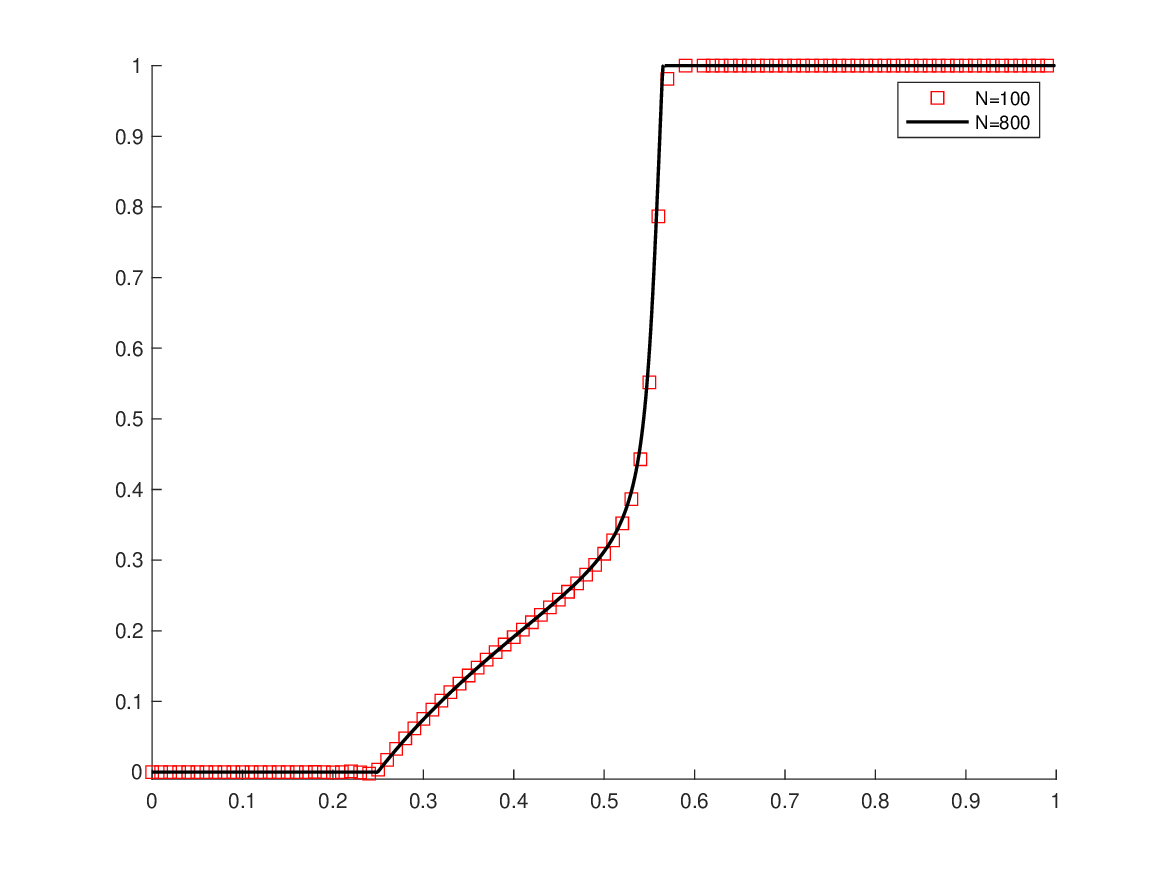}
  \caption{Riemann problem \eqref{eq:ex5_IC2} without gravity. 
  Time-step size satisfies $\frac{\Delta t}{\Delta x}=0.4$.}
 \end{subfigure}
 \begin{subfigure}[b]{0.45\textwidth}
  \includegraphics[width=\textwidth]{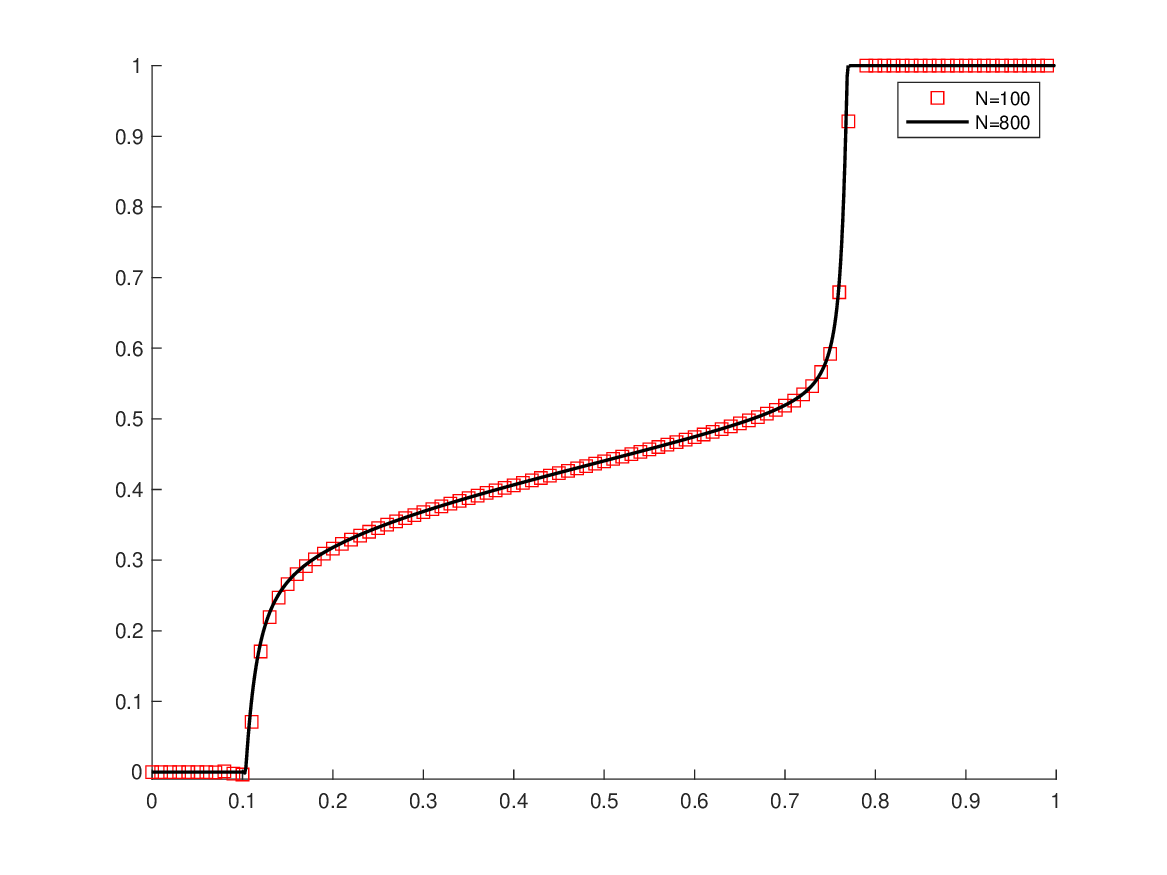}
  \caption{Riemann problem \eqref{eq:ex5_IC2} with gravity.
  Time-step size satisfies $\frac{\Delta t}{\Delta x}=0.1$.}
 \end{subfigure}
 \caption{
 \textbf{Example \ref{ex:ex5}.} 
 Numerical solutions of the Buckley-Leverett equation at $T=0.2$. The ETD-RK4 multi-resolution A-WENO6 scheme is used.
 The time-step sizes are chosen to be the
maximum permitted values for achieving a stable computation.}
 \label{fig:ex5}
\end{figure}

\end{exmp}

\begin{exmp}\label{ex:ex6}
\textbf{A 1D strongly degenerate convection-diffusion equation.}

We consider the convection-diffusion equation
\begin{equation}
u_t+(u^2)_x=\epsilon(\nu(u)u_x)_x,
\end{equation}
with the strongly degenerate viscosity coefficient 
\begin{equation}\label{eq:degevis}
\nu(u)=
\begin{cases}
0,& |u|\leq\frac14,\\
1,& |u|>\frac14.
\end{cases}
\end{equation}
So the equation is hyperbolic when $u\in[-\frac14,\frac14]$ and parabolic otherwise.
The corresponding function of the parabolic term is
\begin{equation}
g(u)=
\begin{cases}
\epsilon(u+\frac14), &u<-\frac14,\\
\epsilon(u-\frac14), &u>\frac14,\\
0,&|u|\leq\frac14.
\end{cases}
\end{equation}

We solve the problem on the domain $\Omega=[-2, 2]$ with $\epsilon=0.1$, the initial condition
\begin{equation}
u(x,0)=
\begin{cases}
1,&-\frac{1}{\sqrt{2}}-\frac25<x<-\frac{1}{\sqrt{2}}+\frac25,\\
-1,& \frac{1}{\sqrt{2}}-\frac25<x<\frac{1}{\sqrt{2}}+\frac25,\\
0,&\text{otherwise},
\end{cases}
\end{equation}
and the homogeneous Dirichlet boundary condition.
The ETD-RK4 multi-resolution A-WENO6 scheme is used to carry out the simulations on the computational grids with $N=100$ and $N=800$. The time-step size is taken as $\Delta t=0.08\Delta x$ such that the desired time-step
size condition, $\Delta t \sim O (\Delta x)$, for the ETD-RK method is satisfied on different grids. 
The obtained numerical results at $T=0.7$ are reported in Figure \ref{fig:ex6}. It is observed that the numerical solutions on different grids match very well, which verifies the convergence of the numerical solutions. The sharp interfaces of the solutions and the kinks where the type of the equation changes are captured stably with high resolution in the simulations, which indicates good nonlinear stability of the method.

\begin{figure}[!htbp]
 \centering
  \includegraphics[width=0.5\textwidth]{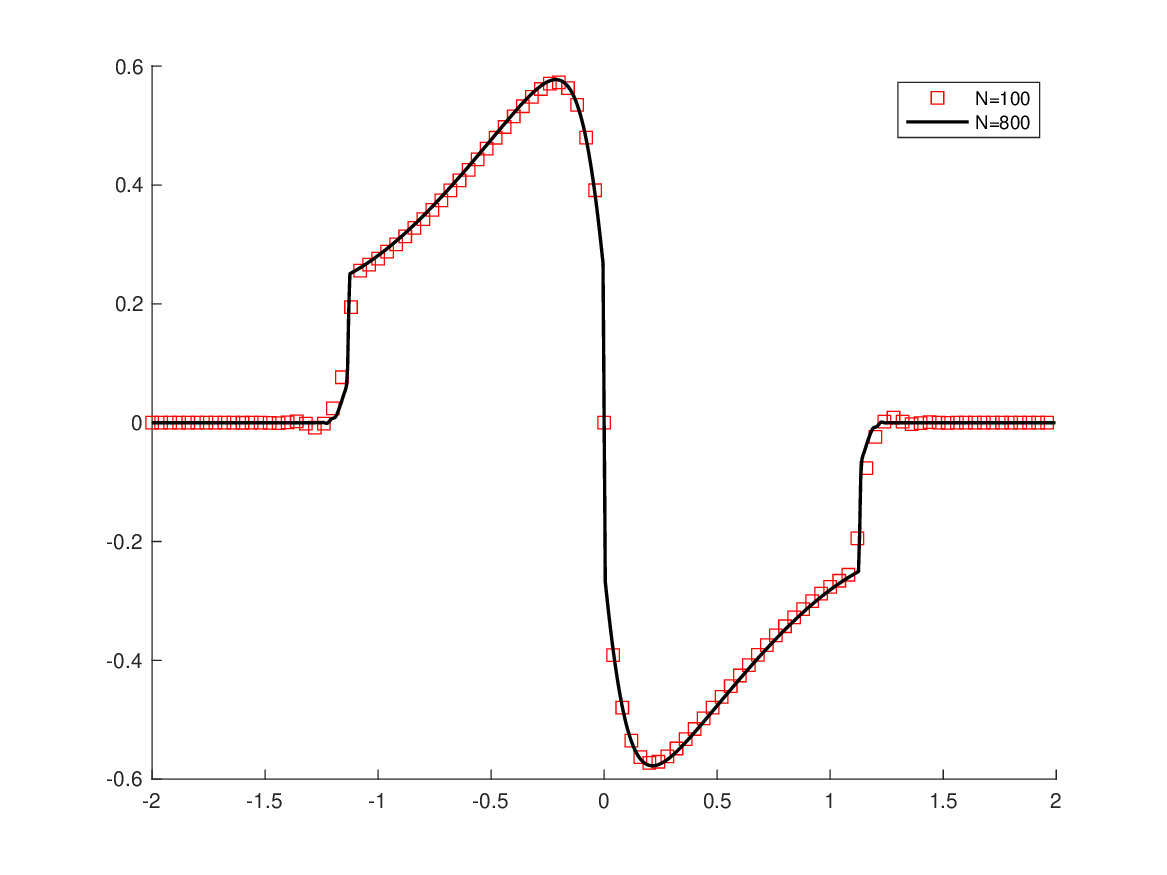}
 \caption{\textbf{Example \ref{ex:ex6}.} 
 Numerical solutions of the 1D strongly degenerate convection-diffusion equation at $T=0.7$. The ETD-RK4 multi-resolution A-WENO6 scheme is used. The time-step size is $\Delta t=0.08\Delta x$.}
 \label{fig:ex6}
\end{figure}

\end{exmp}

\subsection{Two-dimensional (2D) examples}

\begin{exmp}\label{ex:ex7}
\textbf{2D heat equation.}

We consider the two-dimensional heat equation
\begin{equation}\label{eq:heat2D}
u_t=u_{xx}+u_{yy},
\end{equation}
on the domain $\Omega=[-\pi,\pi]^2$ with the periodic boundary condition and the initial condition $u(x,y,0)=\sin(x+y)$.
The exact solution of the problem is $u(x,y,t)=e^{-2t}\sin(x+y)$.

The solution is computed up to $T=1$ using the ETD-RK multi-resolution  A-WENO methods of different temporal and spatial orders with the time-step size $\Delta t=\Delta x$.
The numerical errors and orders of convergence of these methods are reported in Table \ref{tab:ex7_1}.
The desired accuracy orders of the spatial discretizations, instead of the temporal discretizations, are observed in the table. Similar to the example of 1D heat equation, the results verify that for a pure diffusion problem with a smooth solution, the numerical errors of the temporal discretizations using these exponential integrators for the diffusion terms are very small.

In addition, we compare the computational efficiency of the 
 ETD-RK3 method with the 
SSP-ERK3 and the SSP-IRK3 methods, coupled with the fourth-order spatial discretizations (i.e., the $r=2$ case in the section 2.1.1). Specifically, the SSP-ERK3 scheme is coupled with the fourth-order multi-resolution A-WENO discretization, while the fourth-order linear spatial discretization is used for the SSP-IRK3 scheme as that discussed in the section 2.2.4. For the ETD-RK3 scheme, both the A-WENO and the linear spatial discretizations are used in the comparison.   
The computation is conducted on the grids with $N=M=20, 40, 60, \ldots, 140$. 
On all grids, we take $\Delta t=\Delta x$ in the ETD-RK3 and the SSP-IRK3 methods.
The $L^1$ errors versus the $CPU$ times for different methods are shown in Figure \ref{fig:efficiency_ex7}. The numerical results reported in the figure again verify the superiority of the efficiency of the ETD-RK method compared to both the explicit and the implicit SSP-RK methods. It takes much less CPU time costs for the ETD-RK3 scheme than the other two schemes to achieve a similar level of numerical errors on all meshes. 

\begin{table}[!htbp]
\centering
\begin{tabular}{ccccccccc}
\toprule[1.5pt]
\multicolumn{8}{c}{ETD-RK3} \\
\midrule
& MRWENO4 & & MRWENO6 & & MRWENO8\\
\cline{2-5} \cline{6-9}
$N\times M$ & $L^1$ Error & Order & $L^1$ Error & Order & $L^1$ Error & Order\\
\midrule
$20\times20$ & $7.28\times10^{-4}$ & - & $1.14\times10^{-5}$ & - & $1.99\times10^{-7}$ & -\\
$40\times40$ & $4.59\times10^{-5}$ & $3.99$ & $1.81\times10^{-7}$ & $5.98$ & $7.94\times10^{-10}$ & $7.97$ \\
$60\times60$ & $9.08\times10^{-6}$ & $4.00$ & $1.60\times10^{-8}$ & $5.99$ & $3.11\times10^{-11}$ & $7.99$ \\
$80\times80$ & $2.87\times10^{-6}$ & $4.00$ & $2.85\times10^{-9}$ & $6.00$ & $3.05\times10^{-12}$ & $8.08$ \\
$100\times100$ & $1.18\times10^{-6}$ & $4.00$ & $7.47\times10^{-10}$ & $6.00$ & $5.27\times10^{-13}$ & $7.86$ \\
\midrule
\multicolumn{8}{c}{ETD-RK4} \\
\midrule
& MRWENO4 & & MRWENO6 & & MRWENO8\\
\cline{2-5} \cline{6-9}
$N\times M$ & $L^1$ Error & Order & $L^1$ Error & Order & $L^1$ Error & Order & & \\
\midrule
$20\times20$ & $7.28\times10^{-4}$ & - & $1.14\times10^{-5}$ & - & $1.99\times10^{-7}$ & -\\
$40\times40$ & $4.59\times10^{-5}$ & $3.99$ & $1.81\times10^{-7}$ & $5.98$ & $7.94\times10^{-10}$ & $7.97$ \\
$60\times60$ & $9.08\times10^{-6}$ & $4.00$ & $1.60\times10^{-8}$ & $5.99$ & $3.11\times10^{-11}$ & $7.99$ \\
$80\times80$ & $2.87\times10^{-6}$ & $4.00$ & $2.85\times10^{-9}$ & $6.00$ & $3.05\times10^{-12}$ & $8.08$ \\
$100\times100$ & $1.18\times10^{-6}$ & $4.00$ & $7.47\times10^{-10}$ & $6.00$ & $5.28\times10^{-13}$ & $7.86$ \\
\bottomrule[1.5pt]
\end{tabular}
\caption{
\textbf{Example \ref{ex:ex7}.} 
Numerical errors of the ETD-RK multi-resolution A-WENO methods with the time-step size $\Delta t=\Delta x$. MRWENO$2r$ stands for the $2r$-th order multi-resolution A-WENO discretization in space.}
\label{tab:ex7_1}
\end{table}

\begin{figure}[!htbp]
 \centering
  \includegraphics[width=0.5\textwidth]{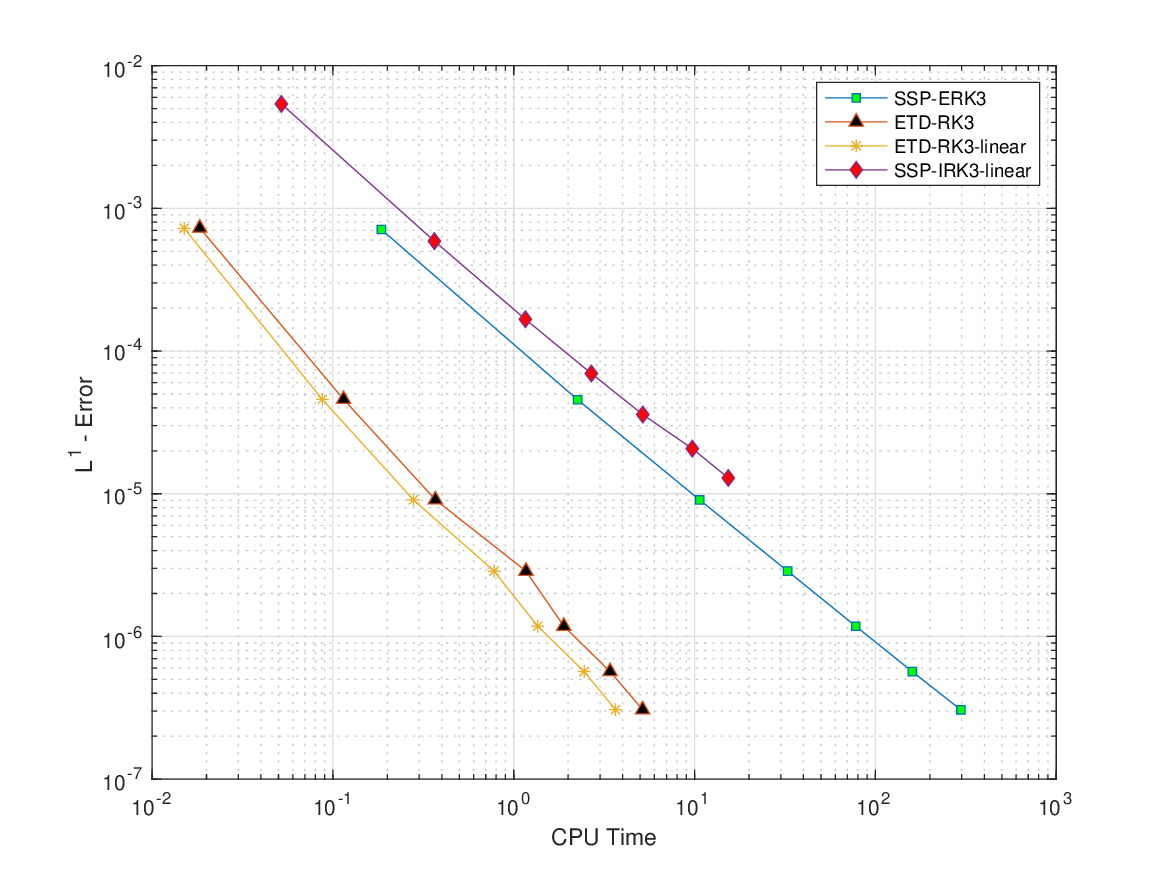}
 \caption{\textbf{Example \ref{ex:ex7}.} 
 Comparison of efficiency for different time-stepping methods.
 SSP-IRK3-linear and ETD-RK3-linear indicate that the linear spatial discretization is used for them. CPU time unit: second.}
 \label{fig:efficiency_ex7}
\end{figure}
\end{exmp}

\begin{exmp}\label{ex:ex8}
\textbf{2D nonlinear stiff reaction-diffusion equation.}

We solve the two-dimensional fully nonlinear stiff reaction-diffusion equation
\begin{equation}\label{eq:PME2DNR}
u_t=128(u^8)_{xx}+128(u^8)_{yy}+R(u),
\end{equation}
where the reaction term $R(u)=-\frac{u}{4}+256 u^{8}+\frac{1}{512 u^7}-2$, on the domain $\Omega=[-\pi, \pi]^2$ with the periodic boundary and the initial condition $u(x,y,0)=\frac12(\sin(x+y)+2)^{\frac18}$.
The problem has the exact solution $u(x,y,t)=\frac12(e^{-2t}\sin(x+y)+2)^{\frac18}$.

The solution is computed up to $T=0.2$ using the ETD-RK multi-resolution  A-WENO methods of different temporal and spatial orders with the time-step size $\Delta t=0.01\Delta x$. Similar to the 1D problem (Example 2), a smaller CFL number is required in this example than Example 7, due to these complex nonlinear stiff diffusion and reaction terms in this 2D example. However, the desired ratio of time-step
size and spatial grid size $\Delta t \sim O (\Delta x)$ can still be preserved in the mesh refinement study. 
The numerical errors and orders of convergence are presented in Table \ref{tab:ex8_1}. We observe that for the fourth-order multi-resolution A-WENO scheme coupled with either the ETD-RK3 or the ETD-RK4 temporal discretizations, the numerical errors of the spatial discretization show strong influence and a fourth-order / close to fourth-order convergence rate is obtained. However, for the sixth-order and the eighth-order multi-resolution A-WENO schemes coupled with either the ETD-RK3 or the ETD-RK4 temporal discretizations, the numerical errors of the temporal discretizations dominate along with the refinement of the meshes, hence the desired accuracy orders of the temporal discretizations, instead of the spatial discretizations, are observed in the table. Again, similar to the numerical experiments of 1D problems, in this 2D example the equation has a highly nonlinear and complex reaction term in addition to the stiff nonlinear diffusion terms, so numerical errors of the spatial discretizations and the temporal discretizations from the different terms have a richer structure than Example 7 which is a linear problem with only the diffusion terms.

In addition, we compare the computational efficiency of the ETD-RK3 method with the SSP-ERK3 and the SSP-IRK3 methods, coupled with the fourth-order spatial discretization (i.e., the $r=2$ case in the section 2.1.1). 
The computation is conducted on the grids with $N=M=50, 100, 150,\ldots, 250, 300$ for the ETD-RK3 method, while for the explicit and diagonally implicit SSP-RK methods, the one on the most refined grid with $N=M=300$ is omitted to save the simulation time, without loss of generality. Similar to the previous examples, the SSP-ERK3 scheme is coupled with
the fourth-order multi-resolution A-WENO discretization, while the fourth-order linear spatial
discretization is used for the SSP-IRK3 scheme as that discussed in the section 2.2.4. For the
ETD-RK3 scheme, both the A-WENO and the linear spatial discretizations are applied in the comparison.
We take the time-step size $\Delta t=0.01\Delta x$ in the computations using the ETD-RK3 method on all grids, while for the SSP-IRK3 method, we take $\Delta t=0.01\Delta x$ in the simulations on the grids with $N=M=50, 100$, $\Delta t=0.005\Delta x$ in the simulations on the grids with $N=M=150,200$, and $\Delta t=0.002\Delta x$ in the simulation on the grid with $N=M=250$ to ensure the stability in solving this 2D stiff nonlinear problem.
The $L^1$ errors versus the CPU times for different methods are shown in Figure \ref{fig:efficiency_ex8}, which verifies that the ETD-RK3 method is much more efficient than both the explicit and the implicit SSP-RK methods here. It takes much less CPU time costs for the ETD-RK3 method than the other two methods to reach a similar level of numerical errors. A consistent conclusion with that in the 1D problem (Example 2) is drawn for this 2D problem. 

\begin{table}[!htbp]
\centering
\begin{tabular}{ccccccccc}
\toprule[1.5pt]
\multicolumn{8}{c}{ETD-RK3} \\
\midrule
& MRWENO4 & & MRWENO6 & & MRWENO8\\
\cline{2-5} \cline{6-9}
$N\times M$ & $L^1$ Error & Order & $L^1$ Error & Order & $L^1$ Error & Order\\
\midrule
$50\times50$ & $2.10\times10^{-4}$ & - & $3.48\times10^{-5}$ & - & $3.43\times10^{-5}$ & -\\
$100\times100$ & $1.54\times10^{-5}$ & $3.77$ & $4.38\times10^{-6}$ & $2.99$ & $4.38\times10^{-6}$ & $2.97$ \\
$150\times150$ & $3.47\times10^{-6}$ & $3.67$ & $1.31\times10^{-6}$ & $2.99$ & $1.30\times10^{-6}$ & $2.98$ \\
$200\times200$ & $1.24\times10^{-6}$ & $3.59$ & $5.52\times10^{-7}$ & $2.99$ & $5.52\times10^{-7}$ & $2.99$ \\
$250\times250$ & $5.64\times10^{-7}$ & $3.52$ & $2.83\times10^{-7}$ & $2.99$ & $2.83\times10^{-7}$ & $2.99$ \\
\midrule
\multicolumn{8}{c}{ETD-RK4} \\
\midrule
& MRWENO4 & & MRWENO6 & & MRWENO8\\
\cline{2-5} \cline{6-9}
$N\times M$ & $L^1$ Error & Order & $L^1$ Error & Order & $L^1$ Error & Order & & \\
\midrule
$50\times50$ & $1.76\times10^{-4}$ & - & $7.68\times10^{-7}$ & - & $3.25\times10^{-7}$ & -\\
$100\times100$ & $1.10\times10^{-5}$ & $4.00$ & $2.75\times10^{-8}$ & $4.80$ & $2.06\times10^{-8}$ & $3.98$ \\
$150\times150$ & $2.17\times10^{-6}$ & $4.00$ & $4.70\times10^{-9}$ & $4.36$ & $4.09\times10^{-9}$ & $3.98$ \\
$200\times200$ & $6.88\times10^{-7}$ & $4.00$ & $1.41\times10^{-9}$ & $4.19$ & $1.30\times10^{-9}$ & $4.00$ \\
$250\times250$ & $2.82\times10^{-7}$ & $4.00$ & $5.62\times10^{-10}$ & $4.12$ & $5.35\times10^{-10}$ & $3.96$ \\
\bottomrule[1.5pt]
\end{tabular}
\caption{
\textbf{Example \ref{ex:ex8}.} 
Numerical errors of the ETD-RK multi-resolution A-WENO methods with the time-step size
$\Delta t=0.01\Delta x$. MRWENO$2r$ stands for the $2r$-th order multi-resolution A-WENO discretization in space.}
\label{tab:ex8_1}
\end{table}

\begin{figure}[!htbp]
 \centering
  \includegraphics[width=0.5\textwidth]{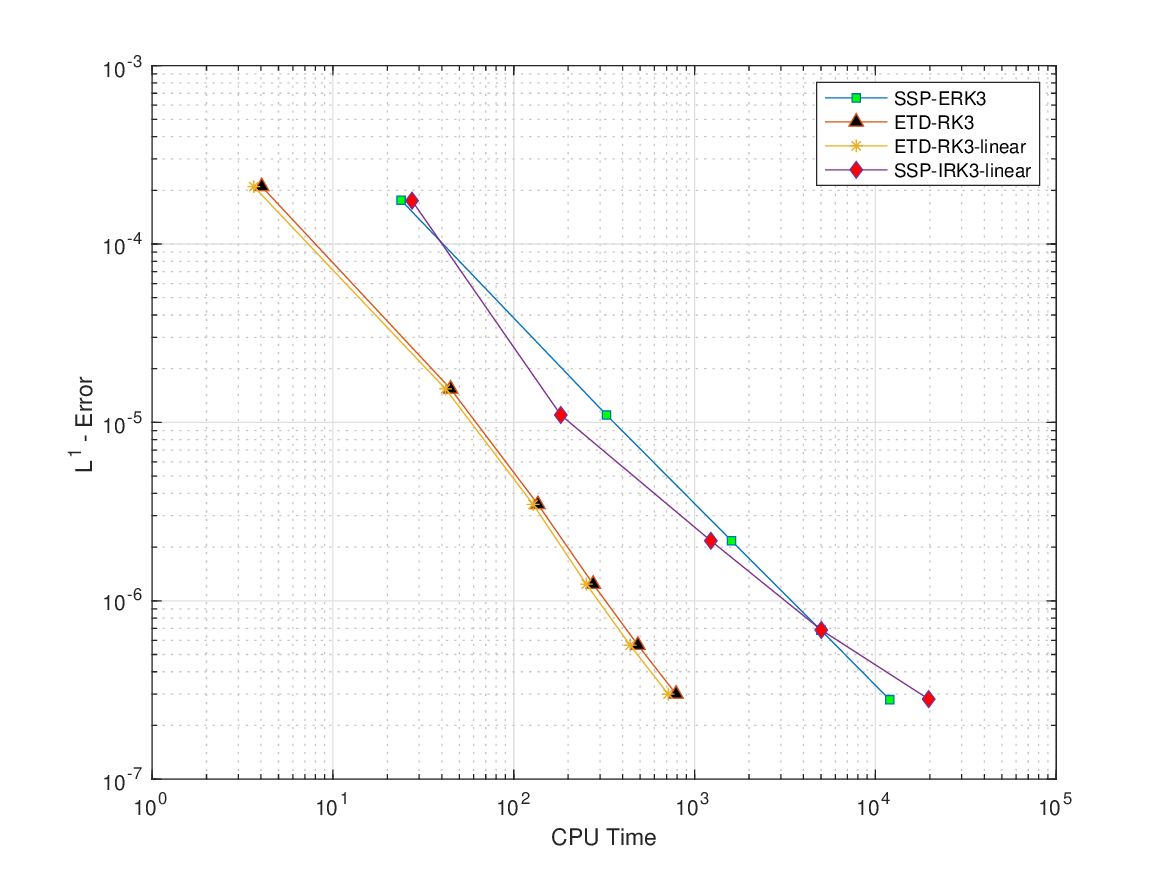}
 \caption{\textbf{Example \ref{ex:ex8}.} 
 Comparison of efficiency for different time-stepping methods.
 SSP-
IRK3-linear and ETD-RK3-linear indicate that the linear
 spatial discretization is used for them. CPU time unit: second.}
 \label{fig:efficiency_ex8}
\end{figure}

\end{exmp}

\begin{exmp}\label{ex:Barenblatt2D}
\textbf{2D PME with the Barenblatt solution.}

We solve the two-dimensional PME
\begin{equation}\label{eq:PME2D}
u_t=(u^m)_{xx}+(u^{m})_{yy}
\end{equation}
with the  Barenblatt solution
\begin{equation}\label{eq:Barenblatt2D}
B_m(x,y,t)=t^{-p}\left(\left(1-\frac{p(m-1)}{4m}\frac{(|x|^2+|y|^2)}{t^{p}}\right)^+\right)^{\frac{1}{m-1}},\quad p=\frac{1}{m}.
\end{equation}
The homogeneous Dirichlet boundary condition is applied.
The solution at $t_0=1$  serves as the initial condition.

The solutions with different values of $m$ are computed up to $T=2$, using the ETD-RK4 scheme coupled with the multi-resolution A-WENO6 spatial discretization. 
To capture the entire non-zero profile of the solutions, we perform the computations on the domain $\Omega=[-6, 6]^2$ with a grid of $N\times M=200\times200$ for the $m=2, 3, 5$ cases, and on the domain $\Omega=[-7, 7]^2$ with a grid of $N\times M=233\times233$ for the $m=8$ case.
The time-step size is taken as  $\Delta t=0.5\Delta x$.
The numerical results for all cases with different values of $m$ are shown in Figure \ref{fig:2DBarenblatt}.
From the numerical results presented in the figures, we observe that the numerical solutions exhibit excellent agreement with the exact solutions and the non-oscillatory performance of the proposed high-order scheme to capture the sharp wave fronts.

In addition, we solve the equation with different values of $m$  up to $T=5$, using different third-order time-marching approaches and comparing their computational efficiency.
The simulations are performed on the domain $\Omega=[-8, 8]$ to include the entire support of a solution, and a computational grid with $N\times M=200\times200$ is used.
The time-step sizes used in different methods are chosen to be at their maximum values to achieve stable computations and  numerical solutions that approximate well the exact solutions.
The comparison of the ratios of the time-step sizes to the spatial grid size, $\Delta t/\Delta x$, and the corresponding CPU times of different methods are shown in Table \ref{tab:2Dbarenblatt}. Again, similar to the 1D PME example (Example 3), the high computational efficiency of the ETD-RK method is verified in solving this 2D problem. From the numerical results in the table, we observe that the ETD-RK3 method, which is coupled with either the multi-resolution A-WENO spatial discretizations or the corresponding linear spatial discretizations, allows for much larger time-step sizes and takes much less CPU time costs than both the explicit and the implicit SSP-RK methods here.
Furthermore, the permitted maximum time-step sizes have very small changes as the stiffness of the equation and the accuracy order of spatial discretization increase, which shows the robustness of the ETD-RK method. 

\begin{figure}[!htbp]
 \centering
 \begin{subfigure}[b]{0.3\textwidth}
  \includegraphics[width=\textwidth]{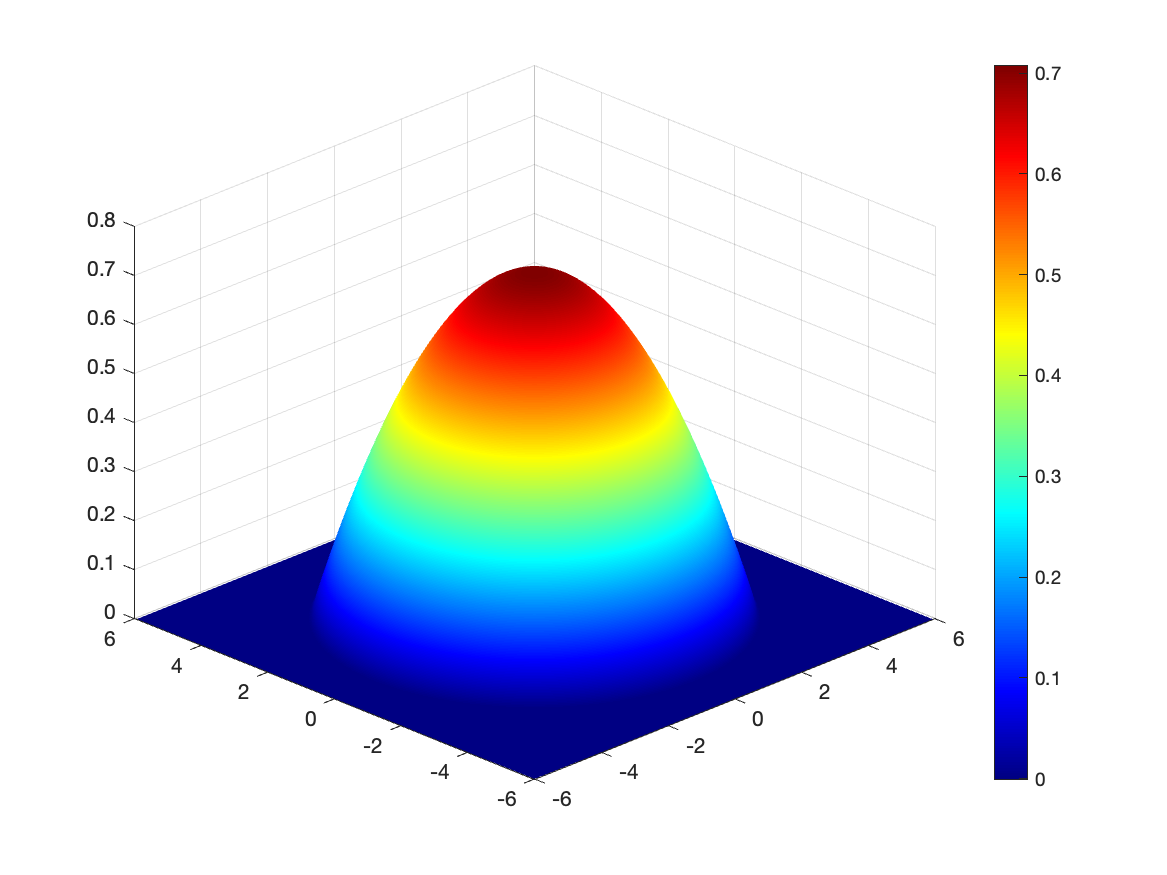}
  \caption{Surface, $m=2$}
 \end{subfigure}
 \begin{subfigure}[b]{0.3\textwidth}
  \includegraphics[width=\textwidth]{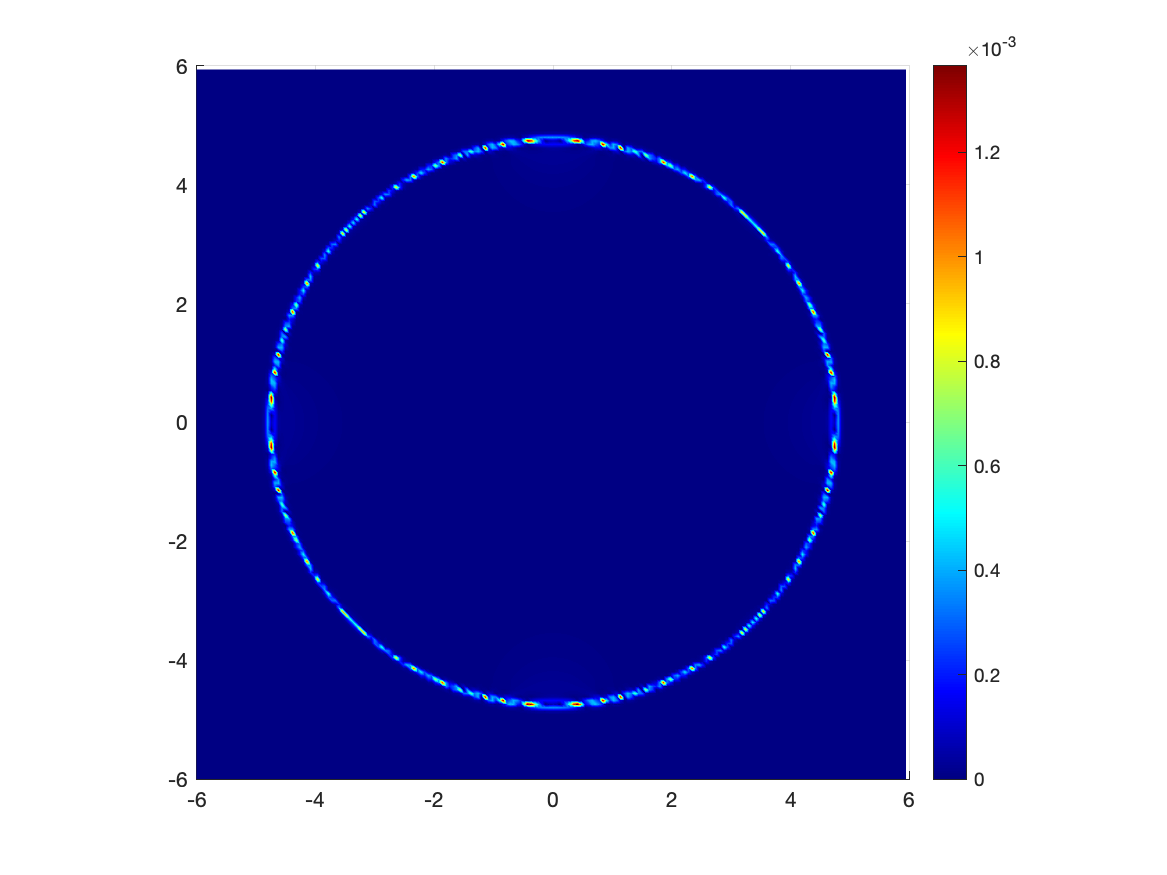}
  \caption{Error, $m=2$}
 \end{subfigure}
 \begin{subfigure}[b]{0.3\textwidth}
  \includegraphics[width=\textwidth]{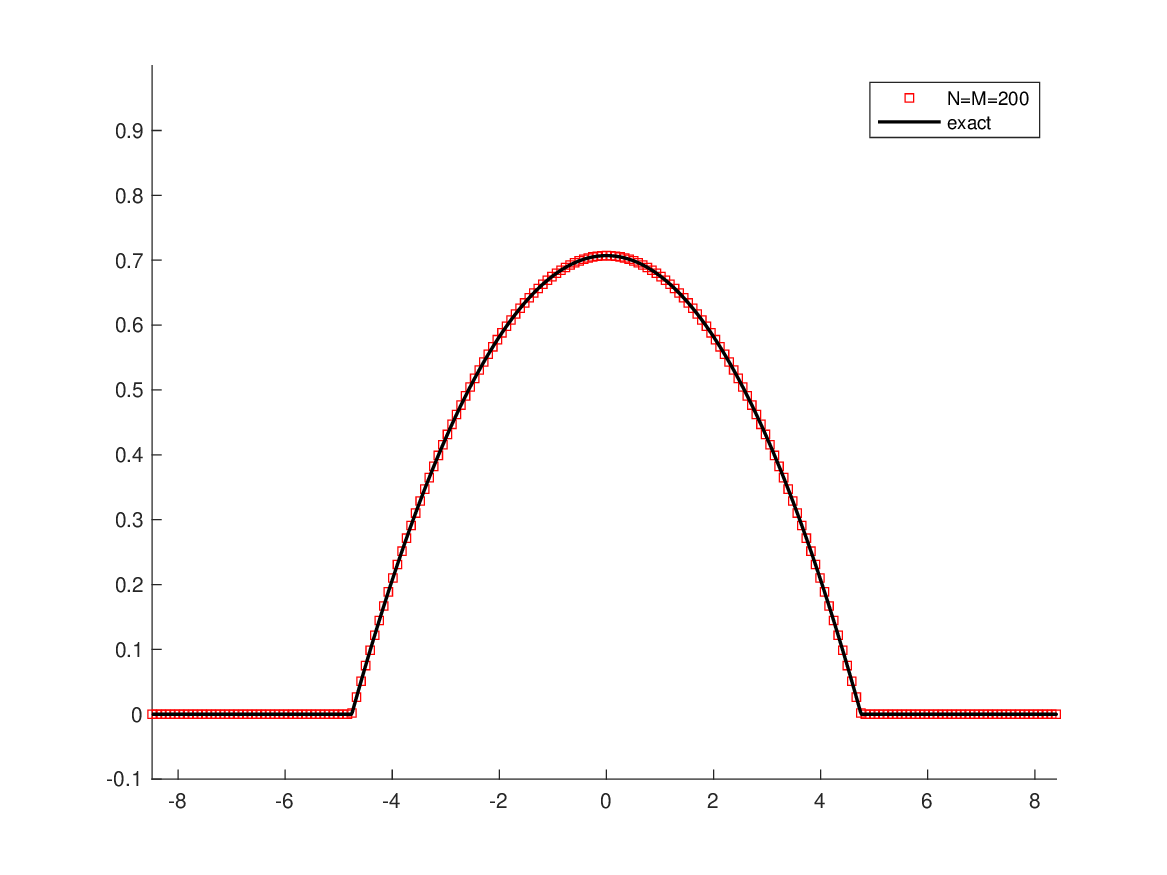}
  \caption{Slice, $m=2$}
 \end{subfigure}

 \begin{subfigure}[b]{0.3\textwidth}
  \includegraphics[width=\textwidth]{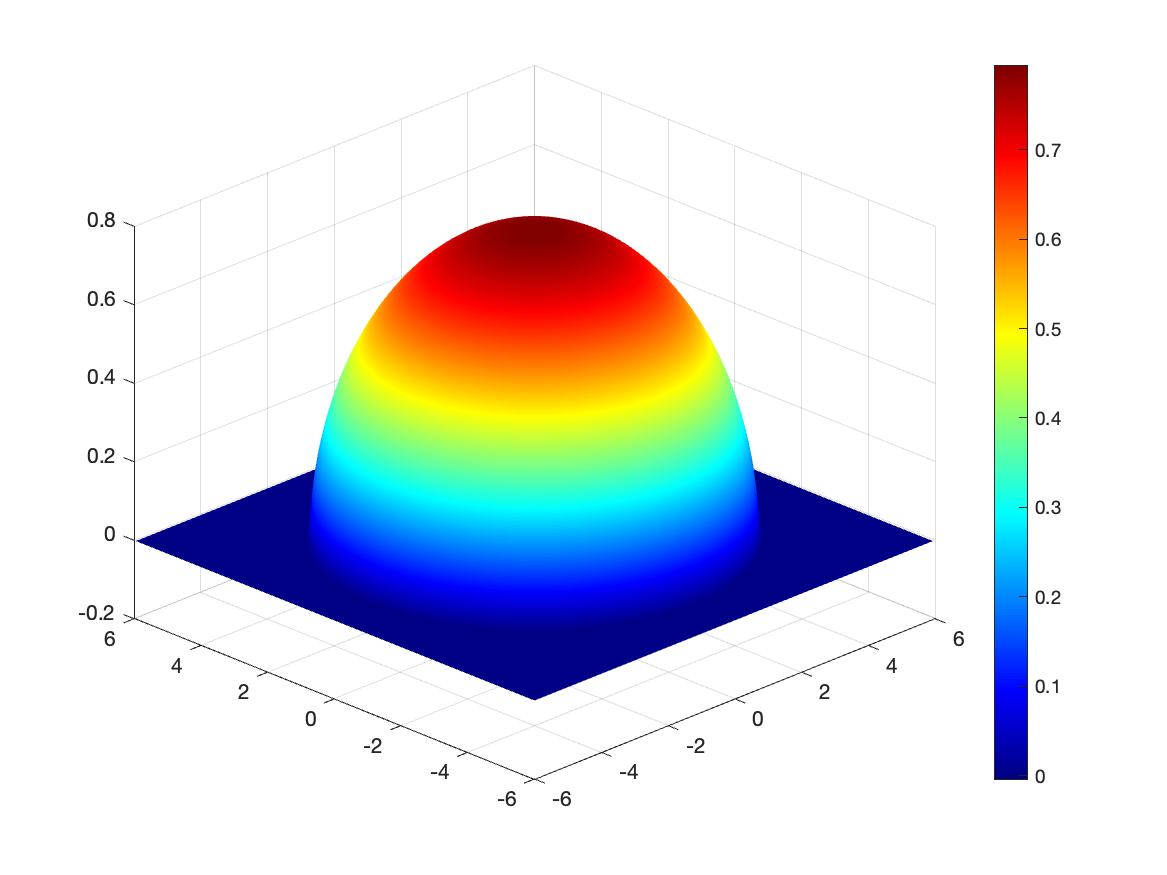}
  \caption{Surface, $m=3$}
 \end{subfigure}
 \begin{subfigure}[b]{0.3\textwidth}
  \includegraphics[width=\textwidth]{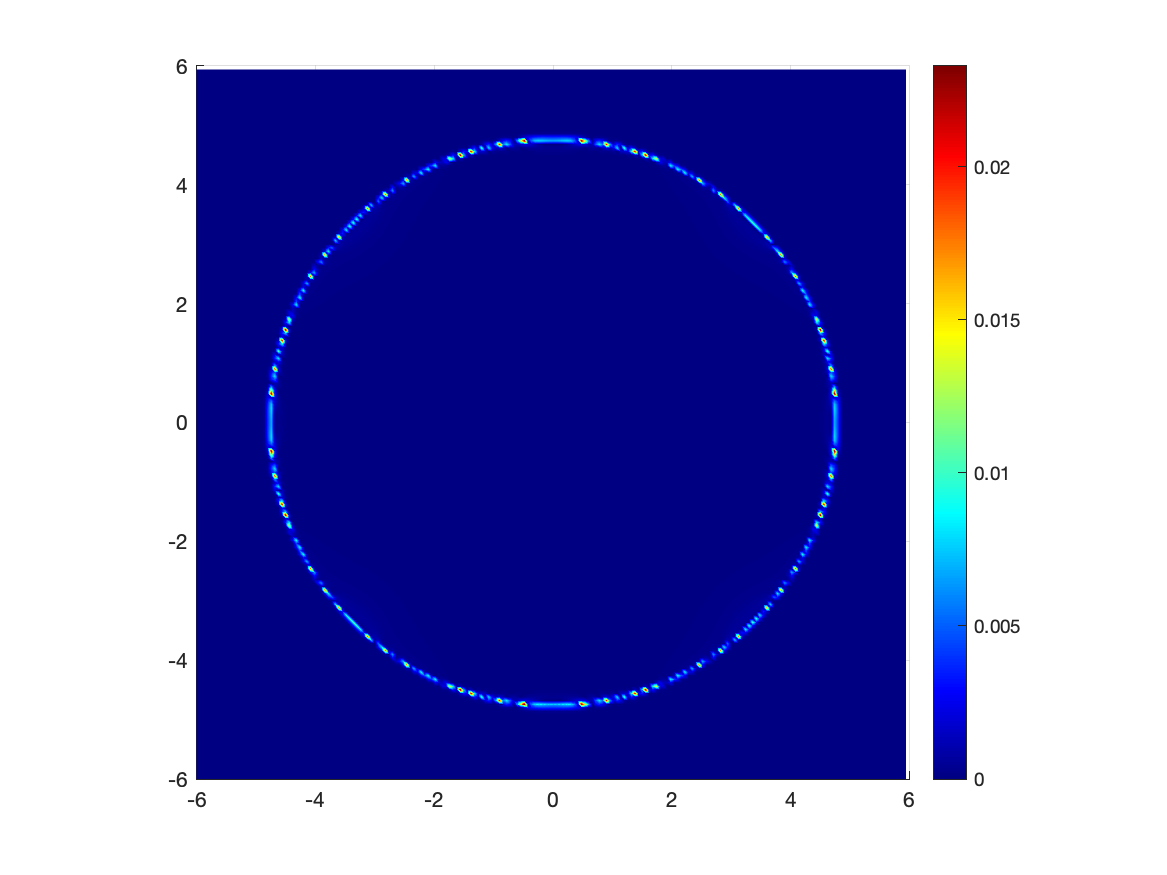}
  \caption{Error, $m=3$}
 \end{subfigure}
 \begin{subfigure}[b]{0.3\textwidth}
  \includegraphics[width=\textwidth]{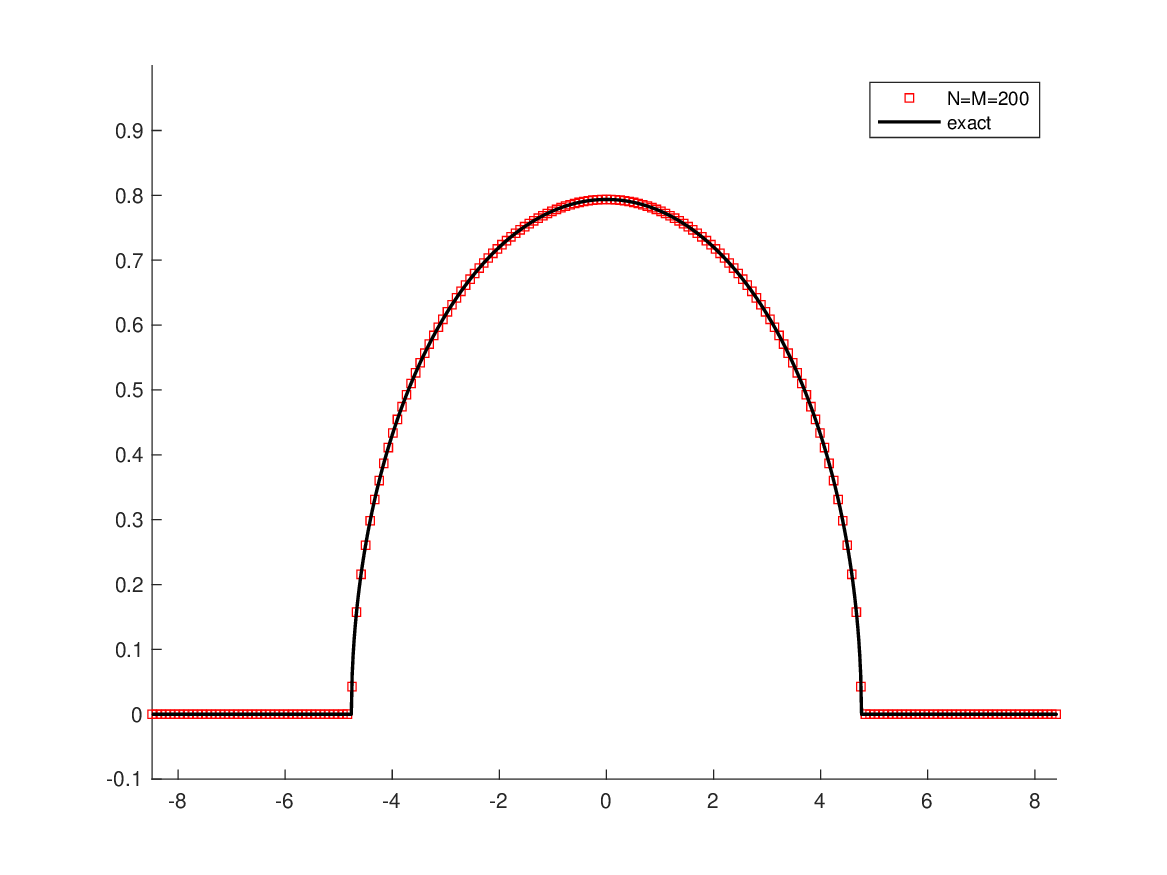}
  \caption{Slice, $m=3$}
 \end{subfigure}

 \begin{subfigure}[b]{0.3\textwidth}
  \includegraphics[width=\textwidth]{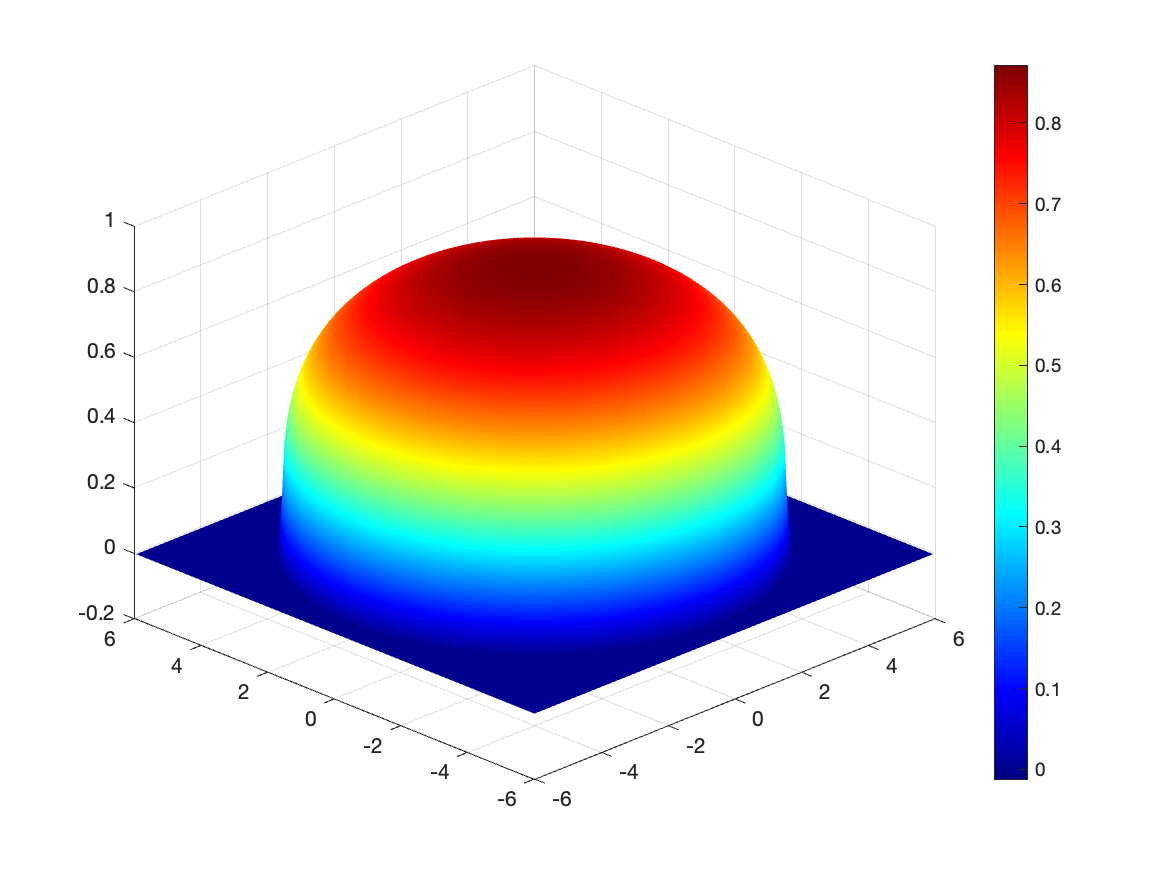}
  \caption{Surface, $m=5$}
 \end{subfigure}
 \begin{subfigure}[b]{0.3\textwidth}
  \includegraphics[width=\textwidth]{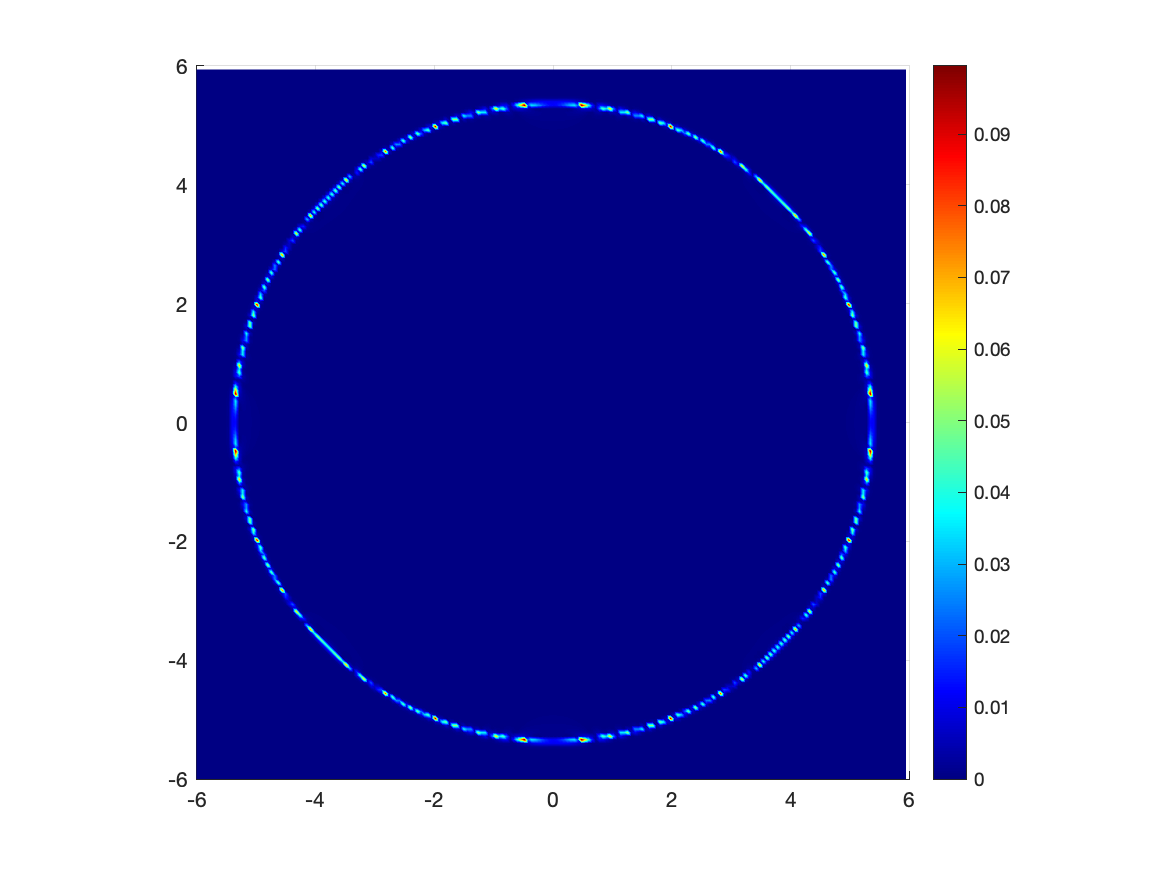}
  \caption{Error, $m=5$}
 \end{subfigure}
 \begin{subfigure}[b]{0.3\textwidth}
  \includegraphics[width=\textwidth]{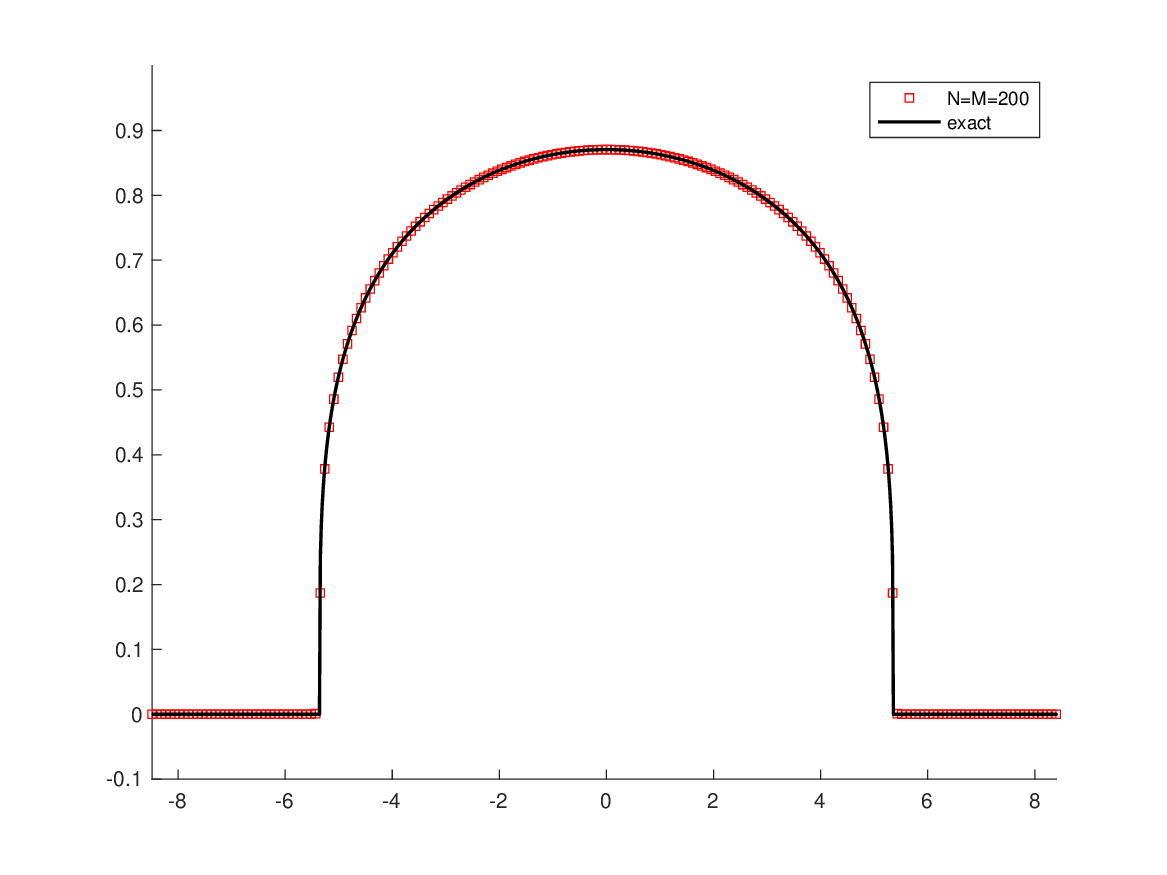}
  \caption{Slice, $m=5$}
 \end{subfigure}

 \begin{subfigure}[b]{0.3\textwidth}
  \includegraphics[width=\textwidth]{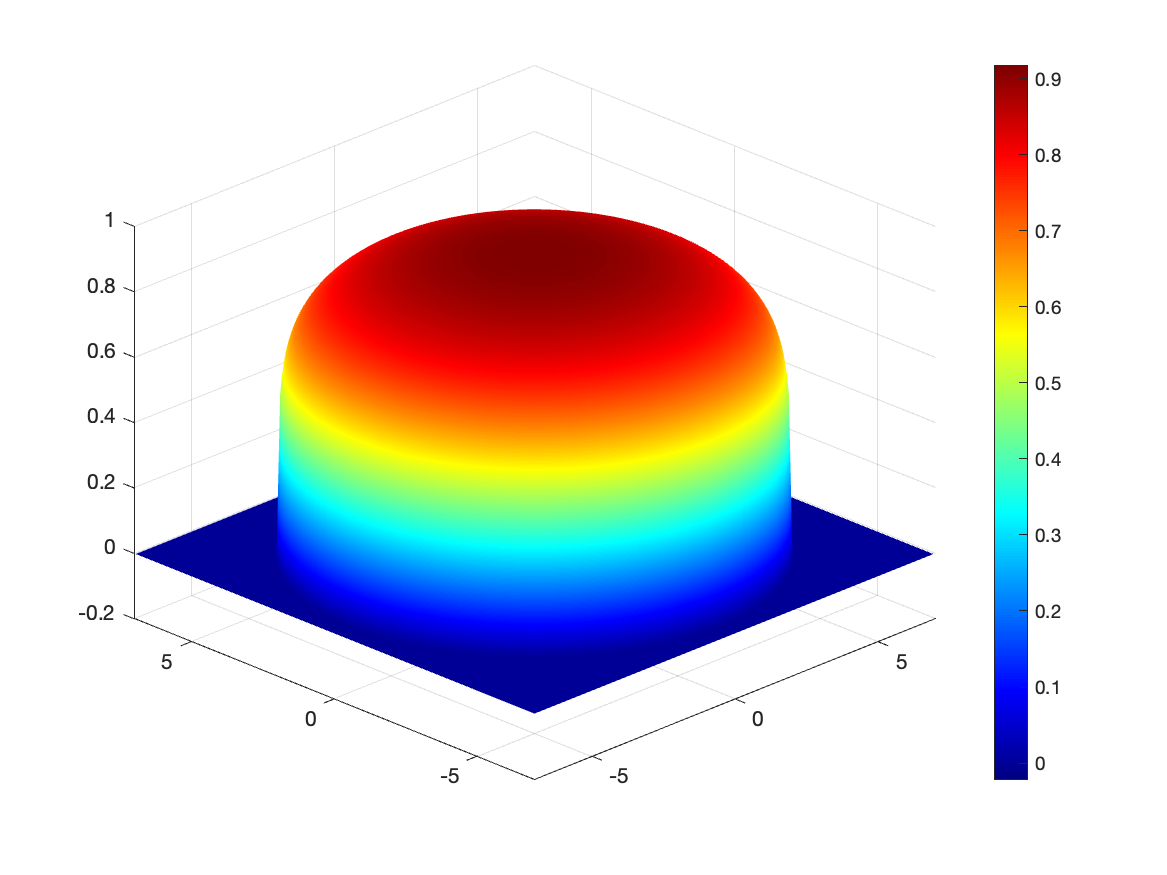}
  \caption{Surface, $m=8$}
 \end{subfigure}
 \begin{subfigure}[b]{0.3\textwidth}
  \includegraphics[width=\textwidth]{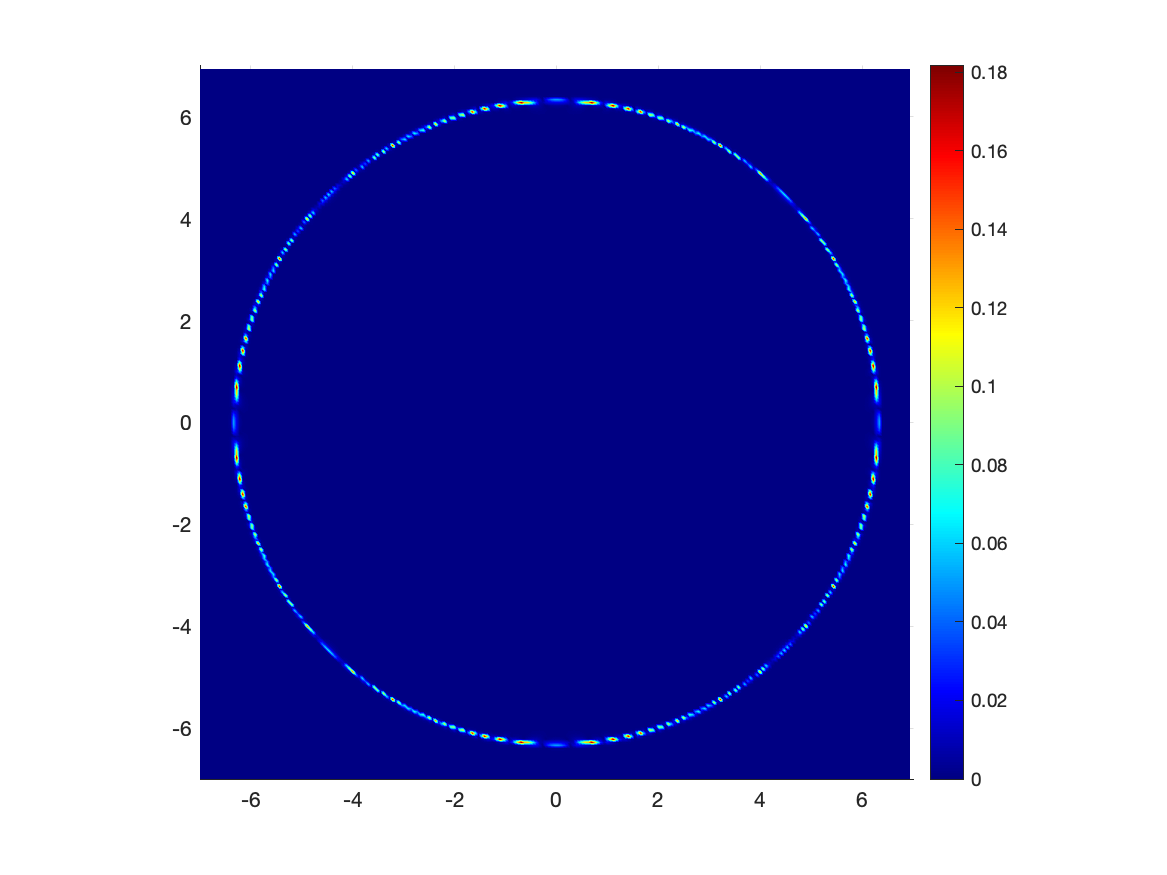}
  \caption{Error, $m=8$}
 \end{subfigure}
 \begin{subfigure}[b]{0.3\textwidth}
  \includegraphics[width=\textwidth]{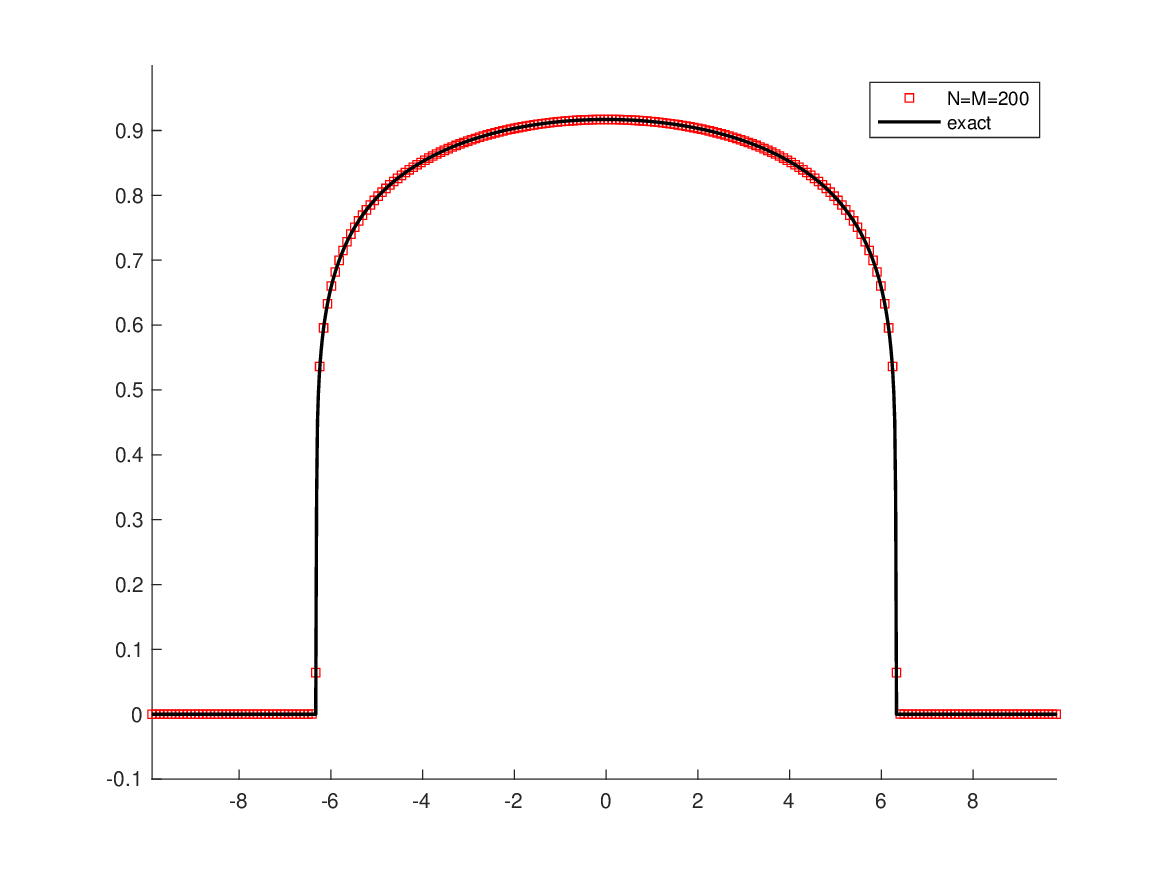}
  \caption{Slice, $m=8$}
 \end{subfigure}
 \caption{\textbf{Example \ref{ex:Barenblatt2D}.} 
 Numerical approximations at $T=2$ to the Barenblatt solutions of the two-dimensional PME with different values of $m$. (a),(d),(g),(j): surface plots of the numerical solutions; (b), (e), (h), (k): contour plots of the numerical errors; (c), (f), (i), (l): 1D slice-plots of the numerical solutions along $x=y$.
 The ETD-RK4 multi-resolution A-WENO6 method is used
 on the grid with $N\times M=200\times200$ for the $m=2,3,5$ cases and on the grid with $N\times M=233\times233$ for the $m=8$ case.
 The time-step size is $\Delta t=0.5\Delta x$.
 }
 \label{fig:2DBarenblatt}
\end{figure}

\begin{table}[!htbp]
\centering
\begin{tabular}{cccccccccccc}
\toprule[1.5pt]
\multicolumn{7}{c}{$m=2$}\\
\midrule
Time & MRWENO4 & & MRWENO6 & & MRWENO8\\
\cline{2-3} \cline{4-5} \cline{6-7}
marching & $\Delta t/\Delta x$  & CPU (s) & $\Delta t/\Delta x$  & CPU (s) & $\Delta t/\Delta x$  & CPU (s)\\
\midrule
ETD-RK3 & $1.0$ & $2.6\times10^{1}$ & $1.0$ & $3.0\times10^{1}$ & $0.9$ & $3.5\times10^{1}$ \\
SSP-ERK3 & $0.0129$ & $8.8\times10^{2}$ & $0.0129$ & $8.9\times10^{2}$ & $0.0129$ & $1.1\times10^{3}$ \\
ETD-RK3-linear & $1.0$ & $2.5\times10^{1}$ & $0.9$ & $3.0\times10^{1}$ & $0.9$ & $3.4\times10^{1}$ \\
SSP-IRK3-linear & $0.2$ & $4.4\times10^{2}$ & $0.2$ & $6.7\times10^{2}$ & $0.2$ & $1.1\times10^{3}$\\
\midrule
\multicolumn{7}{c}{$m=3$}\\
\midrule
Time & MRWENO4 & & MRWENO6 & & MRWENO8\\
\cline{2-3} \cline{4-5} \cline{6-7}
marching & $\Delta t/\Delta x$  & CPU (s) & $\Delta t/\Delta x$  & CPU (s) & $\Delta t/\Delta x$  & CPU (s)\\
\midrule
ETD-RK3 & $1.0$ & $3.2\times10^{1}$ & $1.0$ & $3.5\times10^{1}$ & $1.0$ & $4.0\times10^{1}$ \\
SSP-ERK3 & $0.01$ & $1.1\times10^{3}$ & $0.01$ & $1.2\times10^{3}$ & $0.01$ & $1.4\times10^{3}$ \\
ETD-RK3-linear & $1.0$ & $3.0\times10^{1}$ & $1.0$ & $3.4\times10^{1}$ & $1.0$ & $3.8\times10^{1}$ \\
SSP-IRK3-linear & $0.2$ & $4.5\times10^{2}$ & $0.1$ & $1.3\times10^{3}$ & $0.1$ & $1.8\times10^{3}$ \\
\midrule
\multicolumn{7}{c}{$m=5$} \\
\midrule
Time & MRWENO4 & & MRWENO6 & & MRWENO8\\
\cline{2-3} \cline{4-5} \cline{6-7}
marching & $\Delta t/\Delta x$  & CPU (s) & $\Delta t/\Delta x$  & CPU (s) & $\Delta t/\Delta x$  & CPU (s)\\
\midrule
ETD-RK3 & $1.0$ & $3.9\times10^{1}$ & $1.0$ & $4.4\times10^{1}$ & $1.0$ & $5.0\times10^{1}$ \\
SSP-ERK3 & $0.0067$ & $1.7\times10^{3}$ & $0.0067$ & $1.8\times10^{3}$ & $0.0067$ & $2.1\times10^{3}$ \\
ETD-RK3-linear & $1.0$ & $3.8\times10^{1}$ & $1.0$ & $4.3\times10^{1}$ & $1.0$ & $4.8\times10^{1}$ \\
SSP-IRK3-linear & $0.1$ & $9.7\times10^{2}$ & $0.09$ & $1.7\times10^{3}$ & $0.08$ & $2.3\times10^{3}$ \\
\midrule
\multicolumn{7}{c}{$m=8$} \\
\midrule
Time & MRWENO4 & & MRWENO6 & & MRWENO8\\
\cline{2-3} \cline{4-5} \cline{6-7}
marching & $\Delta t/\Delta x$  & CPU (s) & $\Delta t/\Delta x$  & CPU (s) & $\Delta t/\Delta x$  & CPU (s)\\
\midrule
ETD-RK3 & $1.0$ & $4.9\times10^{1}$ & $1.0$ & $5.6\times10^{1}$ & $1.0$ & $6.4\times10^{1}$ \\
SSP-ERK3 & $0.0045$ & $2.6\times10^{3}$ & $0.0045$ & $2.7\times10^{3}$ & $0.0045$ & $3.1\times10^{3}$ \\
ETD-RK3-linear & $1.0$ & $4.8\times10^{1}$ & $1.0$ & $5.5\times10^{1}$ & $0.9$ & $6.3\times10^{1}$ \\
SSP-IRK3-linear & $0.06$ & $1.9\times10^{3}$ & $0.05$ & $3.2\times10^{3}$ & $0.04$ & $6.1\times10^{3}$ \\
\bottomrule[1.5pt]
\end{tabular}
\caption{\textbf{Example \ref{ex:Barenblatt2D}.} 
Maximum ratios of the time-step sizes to the spatial grid size and the corresponding CPU times of different methods in the computation of the Barenblatt solutions of the two-dimensional PME. MRWENO$2r$ stands for the $2r$-th order multi-resolution A-WENO discretization in space.
SSP-IRK3-linear and ETD-RK3-linear indicate that the corresponding linear spatial discretizations of the MRWENO schemes are used.}
\label{tab:2Dbarenblatt}
\end{table}

\end{exmp}

\begin{exmp}\label{ex:ex9}
\textbf{Merging cones.}

We solve the 2D PME \eqref{eq:PME2D} with $m=2$ and the initial condition \begin{equation}
u(x,y,0)=
\begin{cases}
\exp(\frac{-1}{6-(x-2)^2-(y+2)^2}),& (x-2)^2+(y+2)^2<6,\\
\exp(\frac{-1}{6-(x+2)^2-(y-2)^2}),& (x+2)^2+(y-2)^2<6,\\
0, & \text{otherwise},
\end{cases}
\end{equation}
which contains two cones on the domain $\Omega=[-10, 10]^2$. The homogeneous Dirichlet boundary conditions are applied.

The solution is computed up to $T=4$ using the proposed ETD-RK4 multi-resolution A-WENO6 scheme on the grid with $N\times M=100\times100$. The time-step size is taken as $\Delta t=0.3\Delta x$. The numerical solution profiles at the time $t=0, 0.5, 1,$ and $4$ are shown in Figure \ref{fig:ex9_ETDRK4}, which illustrates the merging process of two cones. Similar to the 1D problem (Example 4), the sharp
wave fronts of the solution of this 2D PME are captured stably with high resolution in the simulation.

\begin{figure}[!htbp]
 \centering
 \begin{subfigure}[b]{0.3\textwidth}
  \includegraphics[width=\textwidth]{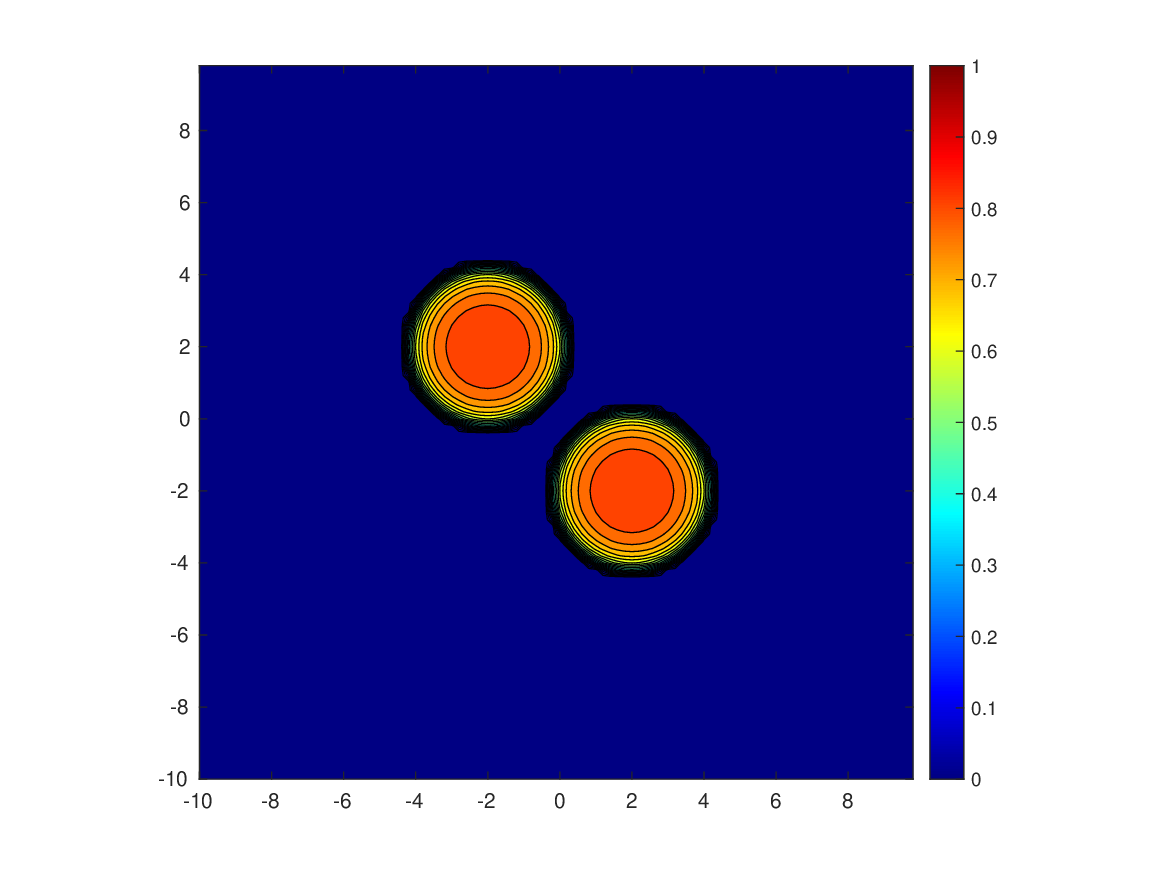}
  \caption{Contour at $t=0$}
 \end{subfigure}
 \begin{subfigure}[b]{0.3\textwidth}
  \includegraphics[width=\textwidth]{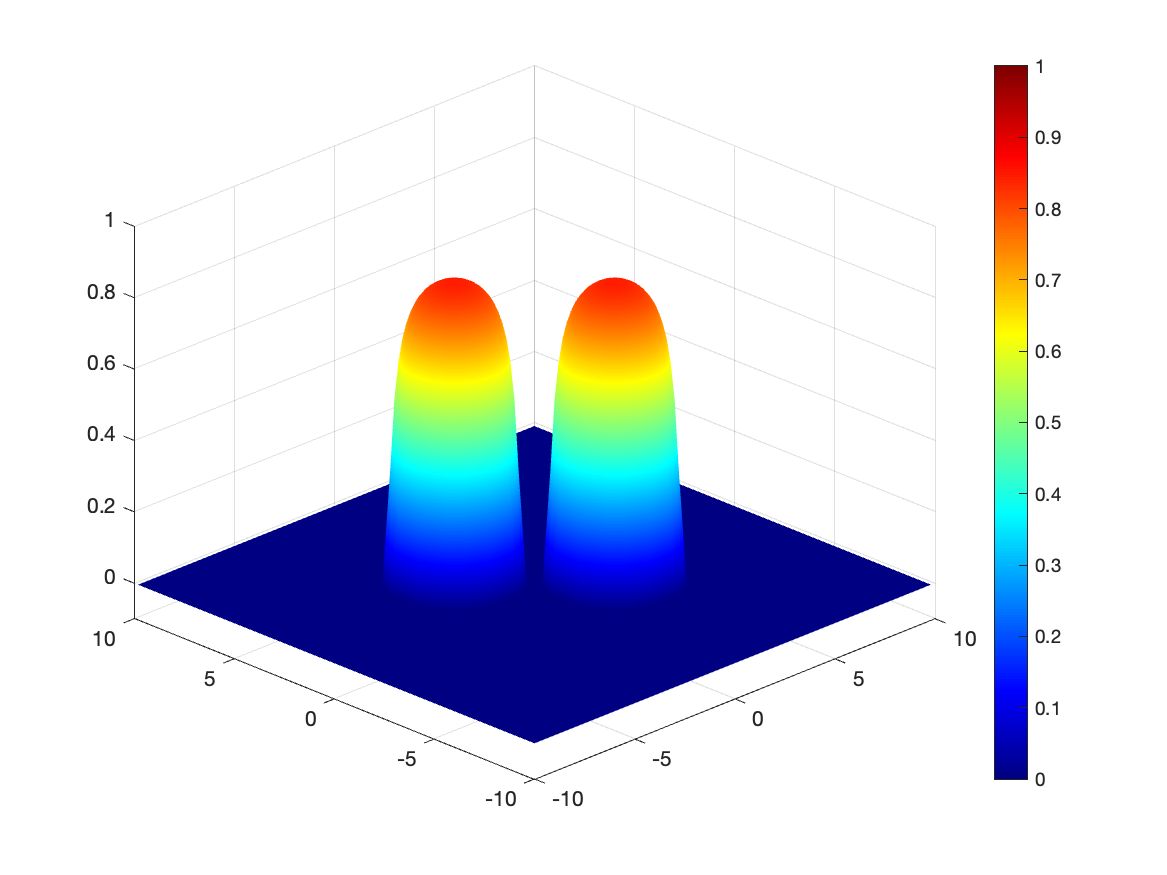}
  \caption{Surface at $t = 0$}
 \end{subfigure}
 
 \begin{subfigure}[b]{0.3\textwidth}
  \includegraphics[width=\textwidth]{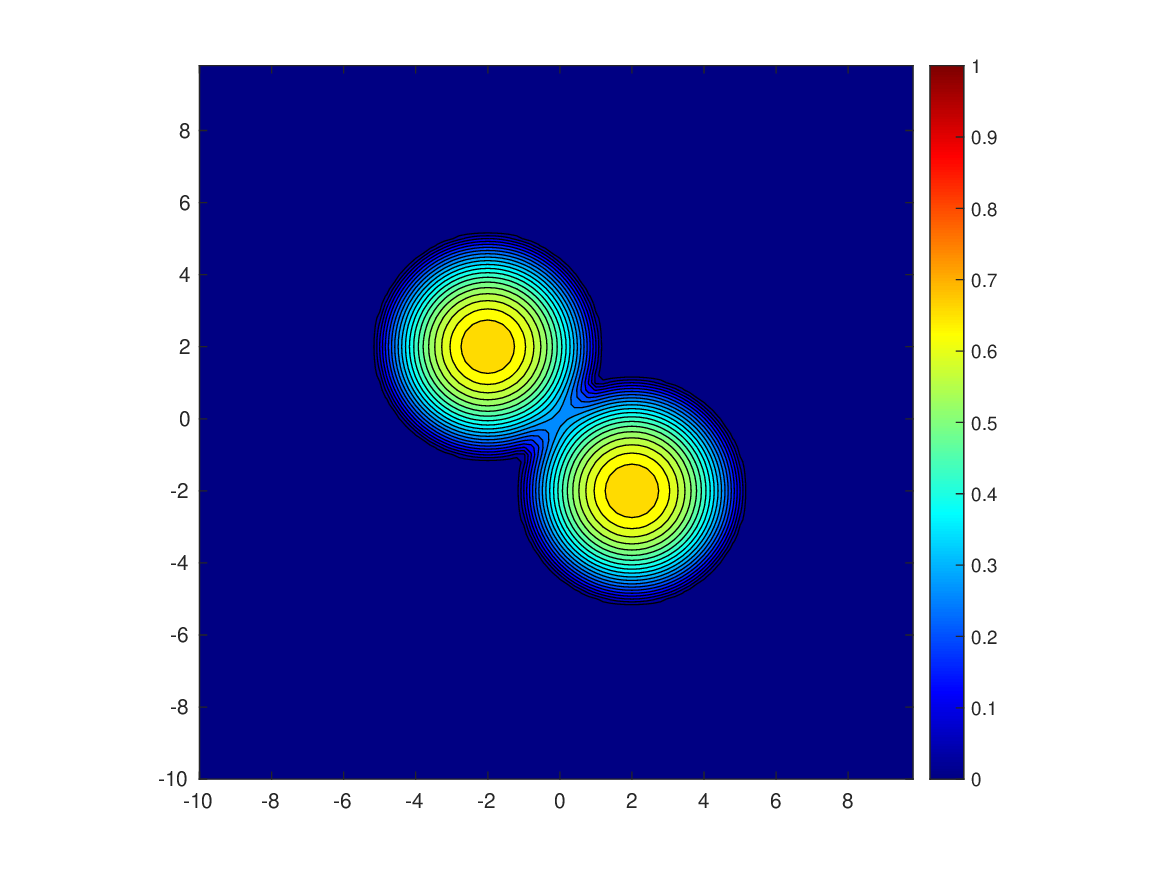}
  \caption{Contour at $t = 0.5$}
 \end{subfigure}
 \begin{subfigure}[b]{0.3\textwidth}
  \includegraphics[width=\textwidth]{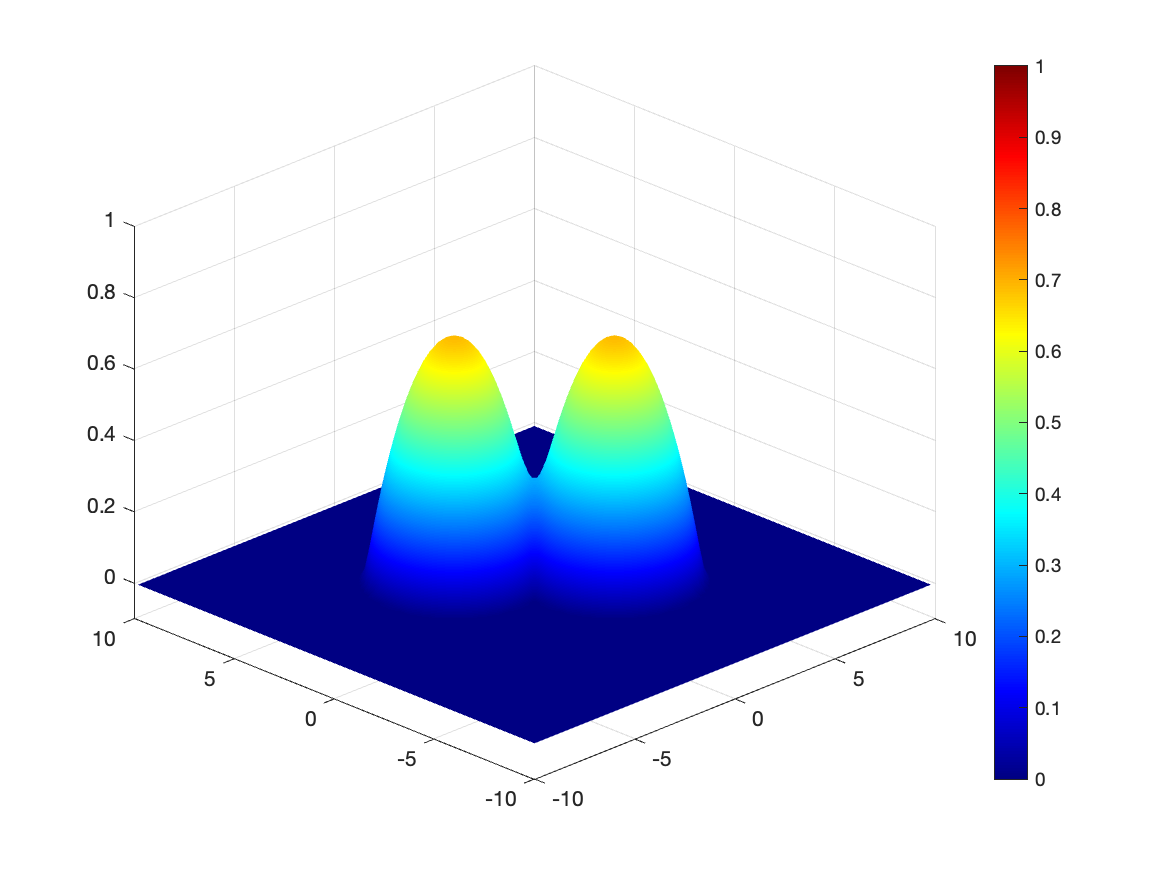}
  \caption{Surface at $t = 0.5$}
 \end{subfigure}

 \begin{subfigure}[b]{0.3\textwidth}
  \includegraphics[width=\textwidth]{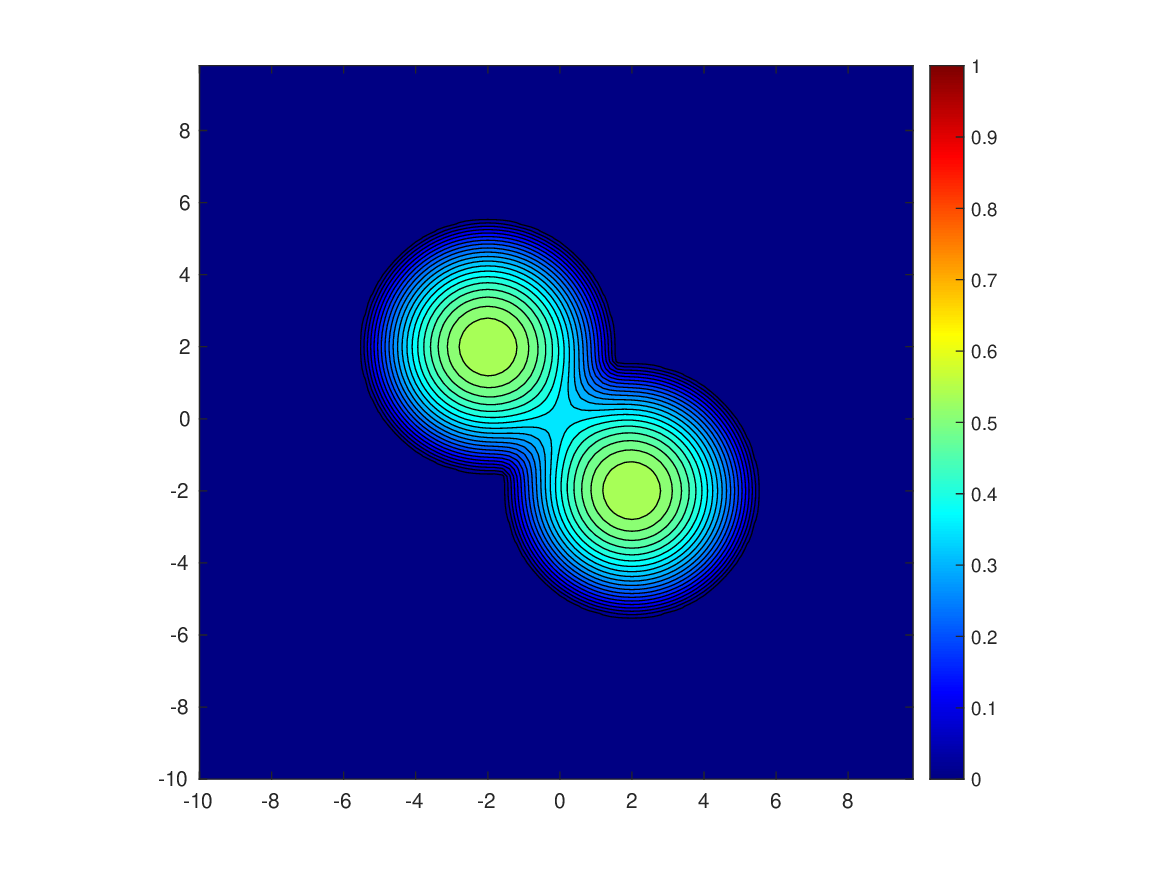}
  \caption{Contour at $t = 1$}
 \end{subfigure}
 \begin{subfigure}[b]{0.3\textwidth}
  \includegraphics[width=\textwidth]{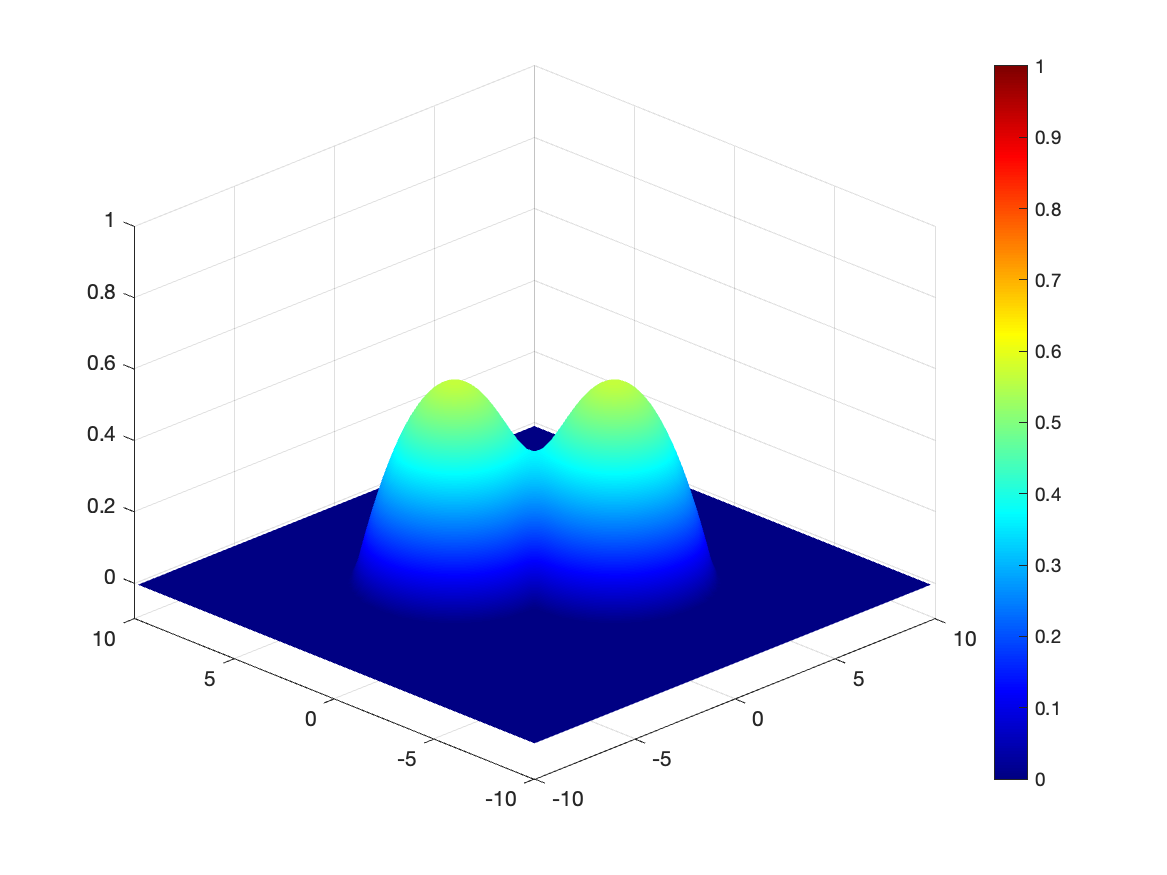}
  \caption{Surface at $t = 1$}
 \end{subfigure}

 \begin{subfigure}[b]{0.3\textwidth}
  \includegraphics[width=\textwidth]{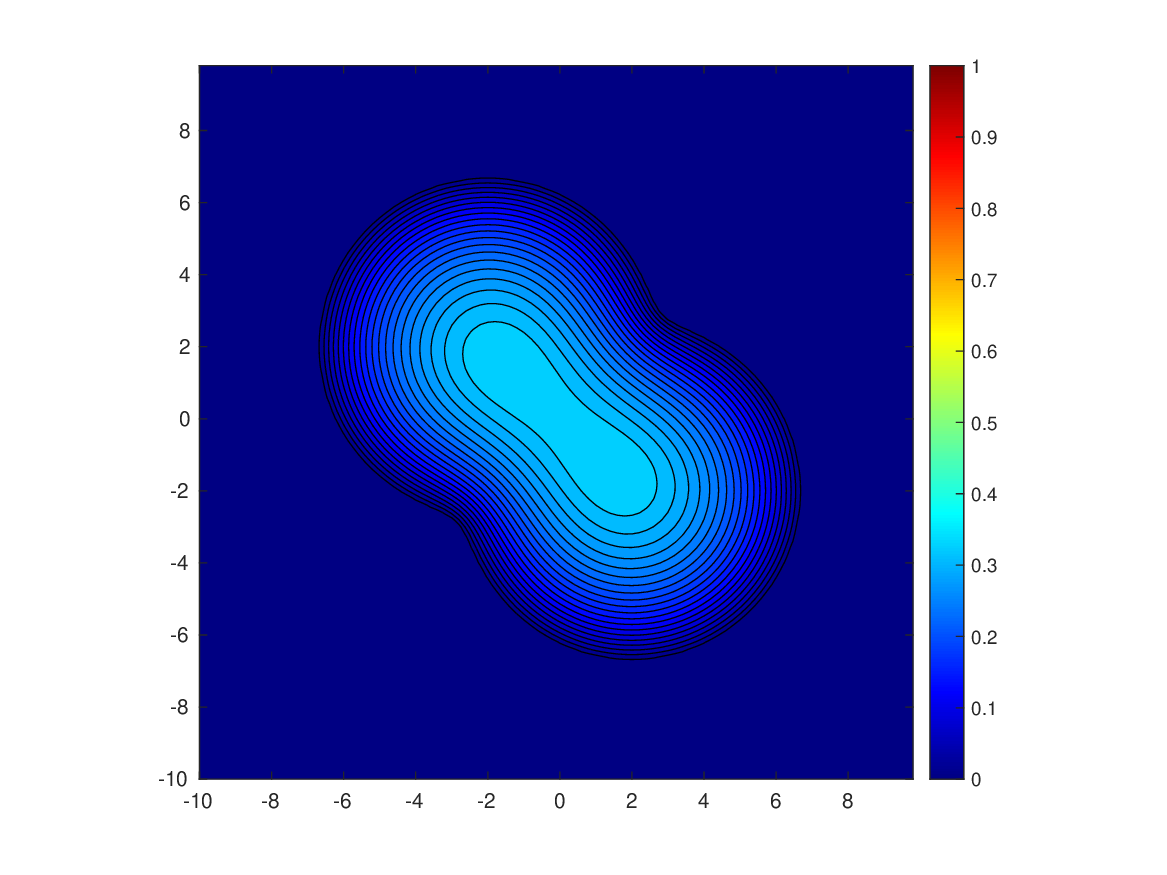}
  \caption{Contour at $t=4$}
 \end{subfigure}
 \begin{subfigure}[b]{0.3\textwidth}
  \includegraphics[width=\textwidth]{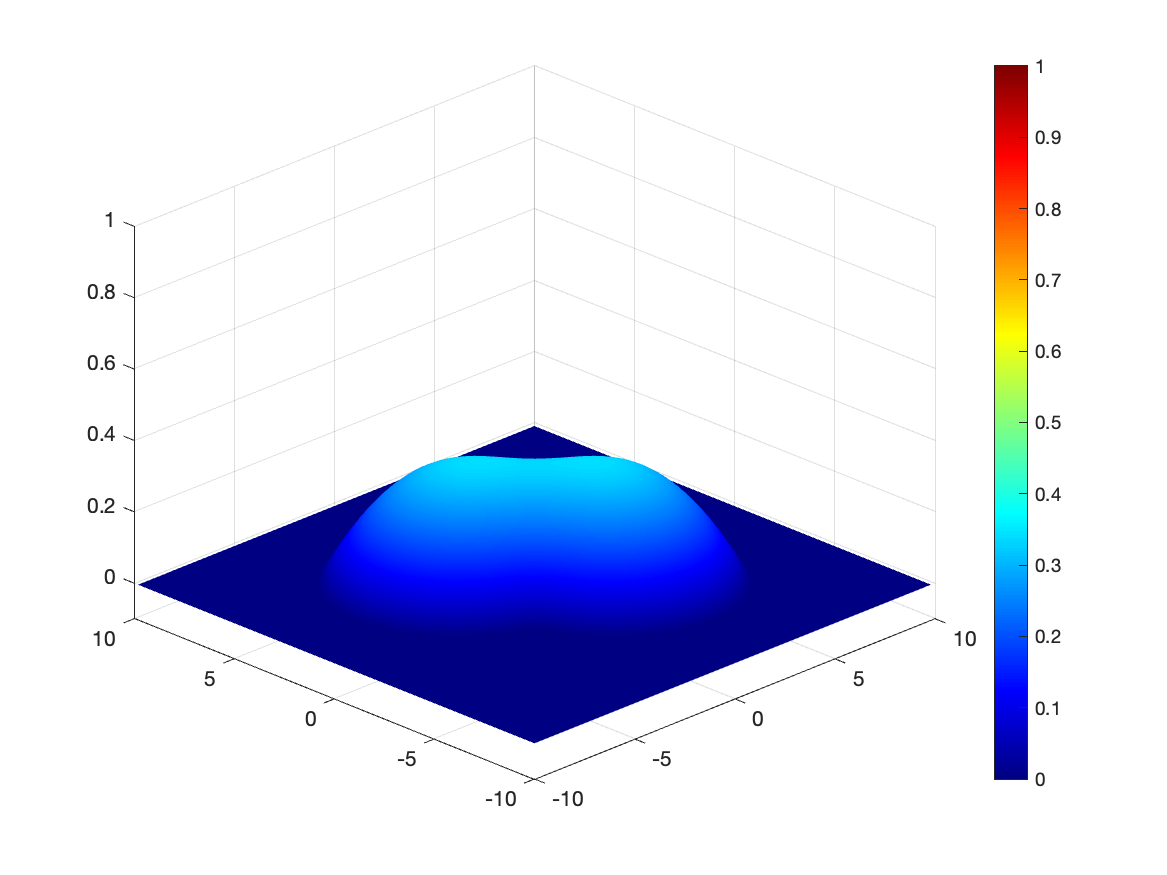}
  \caption{Surface at $t=4$}
 \end{subfigure}
 \caption{\textbf{Example \ref{ex:ex9}.}  
 Numerical solutions for the problem of the merging cones, modeled by the two-dimensional PME. The ETD-RK4 multi-resolution A-WENO6 scheme is use on a grid with $N\times M=100\times100$. The time-step size is $\Delta t=0.3\Delta x$.
 }
 \label{fig:ex9_ETDRK4}
\end{figure}

\end{exmp}

\begin{exmp}\label{ex:ex10}
\textbf{2D Buckley-Leverett equation.}

We consider the two-dimensional viscous Buckley-Leverett equation
\begin{equation}
u_t+f_1(u)_x+f_2(u)_y=\epsilon(u_{xx}+u_{yy}),
\end{equation}
where $f_1(u)$ and $f_2(u)$ are the fluxes without and with gravitational effects as in the one-dimensional case \eqref{eq:BL1DNG} and \eqref{eq:BL1DGR} respectively, and $\epsilon=0.01$.

We solve the problem on the domain $\Omega=[-\frac32, \frac32]^2$ with the initial condition 
\begin{equation}
u(x,y,0)=
\begin{cases}
1, & x^2+y^2<\frac12,\\
0, &\text{otherwise},
\end{cases}
\end{equation}
and the homogeneous Dirichlet boundary condition. The proposed  ETD-RK4 multi-resolution A-WENO6 scheme is applied in the simulation and the computation is performed on the grid with $N\times M=120\times120$. The time-step size is taken as $\Delta t=0.3\Delta x$. The obtained numerical solution at $T=0.5$ is shown in Figure \ref{fig:ex10}. Similar to the 1D example (Example 5), the large gradients of the solution are captured stably with high resolution and without numerical oscillation in the simulation of this 2D problem, which verifies the nonlinear stability of the proposed method.

\begin{figure}[!htbp]
	\centering
	\begin{subfigure}[b]{0.45\textwidth}
		\includegraphics[width=\textwidth]{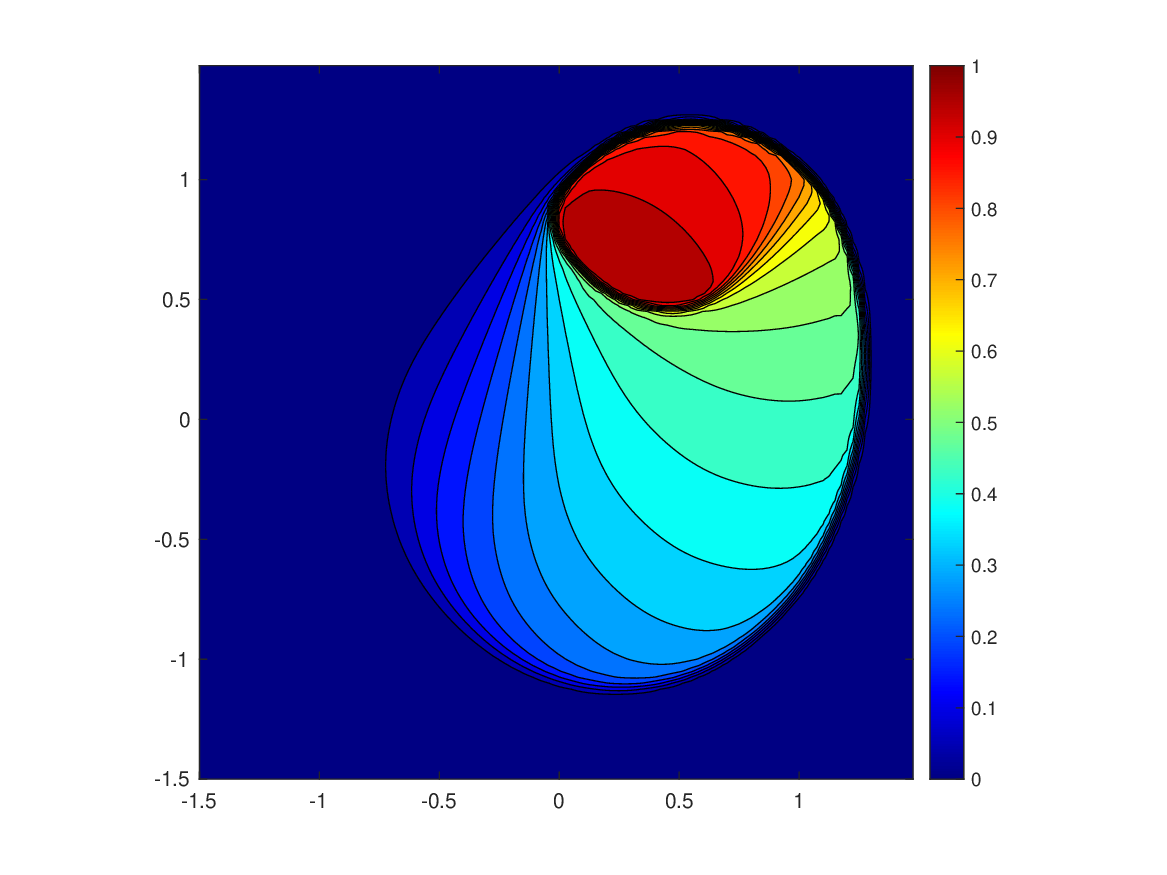}
		\caption{Contour}
	\end{subfigure}
	\begin{subfigure}[b]{0.45\textwidth}
		\includegraphics[width=\textwidth]{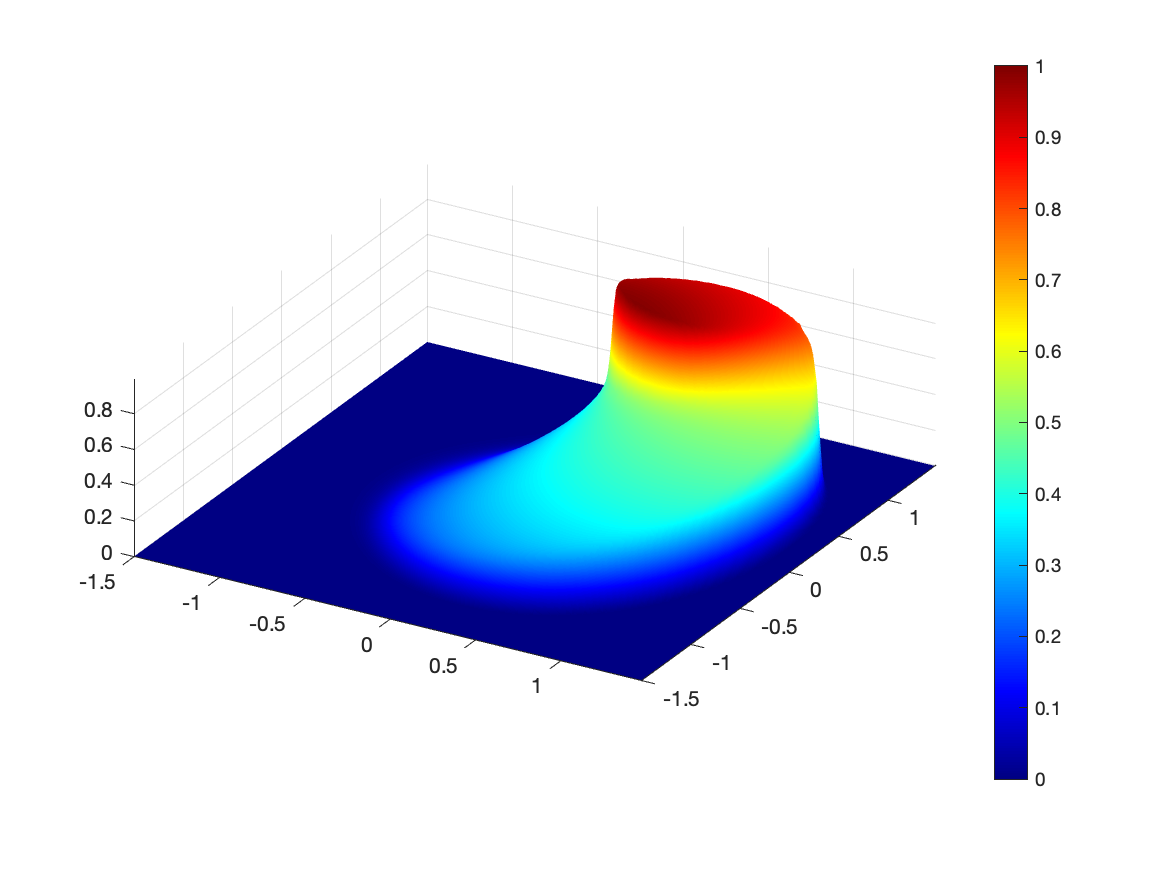}
		\caption{Surface}
	\end{subfigure}
	\caption{\textbf{Example \ref{ex:ex10}.} 
 Numerical solution of the two-dimensional Buckley–Leverett equation at $T=0.5$. The ETD-RK4 multi-resolution A-WENO6 scheme is used to solve the problem on a grid with $N\times M=120\times120$. The time-step size is $\Delta t=0.3\Delta x$.
 }
	\label{fig:ex10}
\end{figure}

\end{exmp}

\begin{exmp}\label{ex:ex11}
\textbf{A 2D strongly degenerate convection-diffusion equation.}

We consider the two-dimensional strongly degenerate parabolic convection-diffusion equation
\begin{equation}
u_t+(u^2)_x+(u^2)_y=\epsilon(\nu(u)u_x)_x+\epsilon(\nu(u)u_y)_y,
\end{equation}
where $\nu(u)$ is the same as in the one-dimensional case \eqref{eq:degevis}.

The problem is solved on the domain $\Omega=[-\frac32, \frac32]^2$ with $\epsilon=0.1$, the initial condition
\begin{equation}
u(x,y,0)=
\begin{cases}
1, & (x+\frac12)^2+(y+\frac12)^2<\frac{4}{25},\\
-1, &(x-\frac12)^2+(y-\frac12)^2<\frac{4}{25},\\
0,&\text{otherwise},
\end{cases}
\end{equation}
and the homogeneous Dirichlet boundary condition.
the ETD-RK4 multi-resolution A-WENO6 scheme is used to perform the simulation on the computational grid with $N\times M=120\times120$. The time-step size is taken as $\Delta t=0.1\Delta x$.
The obtained numerical solution at $T=0.5$ is shown in Figure \ref{fig:ex11}. Again, similar to the 1D example (Example 6), The sharp interfaces of the solution are captured stably with high resolution and without numerical oscillation in the simulation of this 2D degenerate convection-diffusion problem, which verifies the nonlinear stability of the proposed 
numerical method.

\begin{figure}[!htbp]
	\centering
	\begin{subfigure}[b]{0.45\textwidth}
		\includegraphics[width=\textwidth]{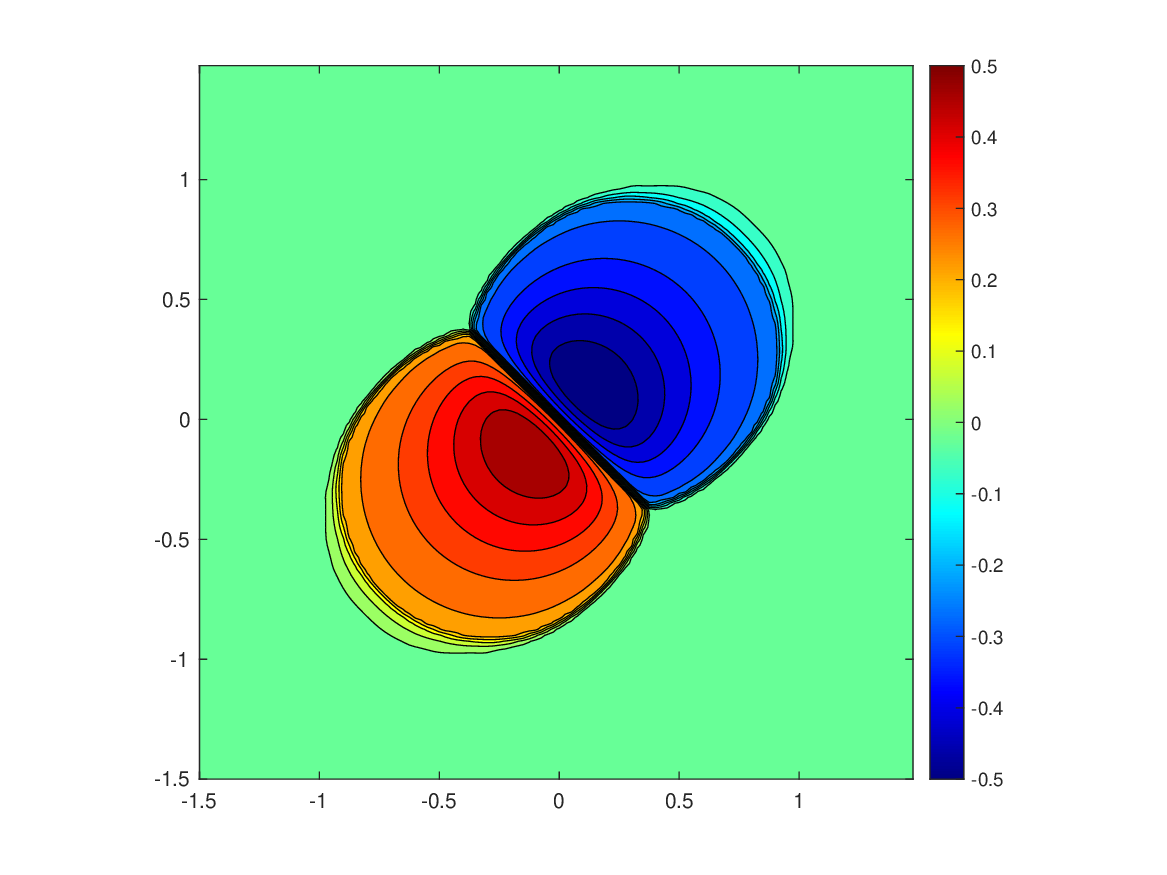}
		\caption{Contour}
	\end{subfigure}
	\begin{subfigure}[b]{0.45\textwidth}
		\includegraphics[width=\textwidth]{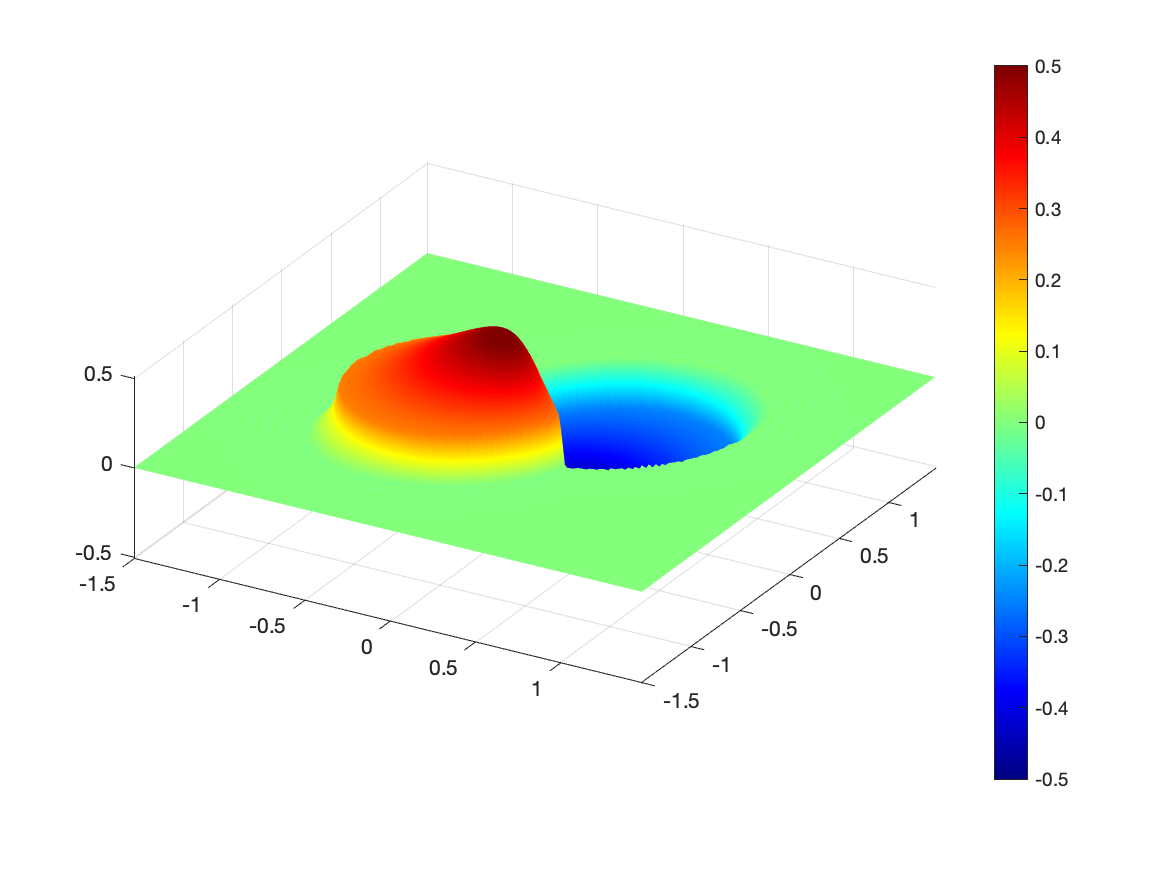}
		\caption{Surface}
	\end{subfigure}
	\caption{\textbf{Example \ref{ex:ex11}.} 
 Numerical solution of the 2D strongly degenerate convection-diffusion equation at $T=0.5$. The ETD-RK4 multi-resolution A-WENO6 scheme is used on a grid with $N\times M=120\times120$. The time-step size is $\Delta t=0.1\Delta x$.
 }
	\label{fig:ex11}
\end{figure}

\end{exmp}

\section{Conclusions}\label{sect:4}

High-order WENO methods have been developed in the literature to effectively solve nonlinear degenerate parabolic equations with high resolution. Since the sophisticated nonlinear properties and high-order accuracy of WENO methods require more operations than many other schemes, their computational costs increase significantly when they are applied to stiff degenerate parabolic PDEs and the time-step sizes are small. How to achieve fast computations of high-order WENO methods is a very important question. To deal with this issue, in this paper we apply 
the ETD-RK methods, a class of accurate exponential integrators, to the high-order multi-resolution alternative finite difference  WENO methods for solving degenerate parabolic equations. A novel and effective semilinearization approach, namely replacing the exact Jacobian of a high-order WENO scheme with that of the corresponding high-order linear scheme, is proposed to efficiently form the linear stiff part in applying the exponential integrators.  
Extensive numerical experiments are performed to demonstrate 
the effectiveness of this new approach, and verify high-order accuracy, nonlinear stability and high efficiency of the developed ETD-RK multi-resolution alternative WENO methods.
The ETD-RK methods resolve the stiffness of the nonlinear degenerate parabolic PDEs very well, and the desired large time-step size simulations of $\Delta t \sim O (\Delta x)$ are achieved. The comparisons with some commonly used explicit and implicit SSP-RK methods show that the proposed methods are more
efficient in solving the nonlinear degenerate parabolic PDEs,   especially the two-dimensional problems. We expect that the proposed novel semilinearization approach in applying the exponential integrators to a high-order nonlinear spatial discretization can be extended to other high-order nonlinear schemes besides the WENO schemes. This interesting topic will be 
investigated in our next research. 

%%%%%%%%%%%%%%%%%%%%%%%%%%%%%%%%%%%%%%%%%%%%%%%%%%%%%%%%%%%%%%%%%%%%

%\newpage
\bigskip
\noindent {\bf Conflict of Interest:} 
The authors declare that they have no known competing financial interests or personal
relationships that could have appeared to influence the work reported in this paper.

\appendix

\section{The smoothness indicators of the multi-resolution WENO in Section \ref{Sect:DiffDisc}} \label{appd:diffusion}

The expressions of the smoothness indicators $\beta_1, \beta_2, \beta_3$ and $\beta_4$ in the definition \eqref{eq:smoothindDiff} are given as follows:
\begin{equation*}
\begin{split}
\beta_1=&(-g_i + g_{i+1})^2,\\
\beta_2=&\frac{781}{720}(g_{i-1}-3g_i+3g_{i+1}-g_{i+2})^2 + \frac{13}{48}(g_{i-1}-g_i-g_{i+1}+g_{i+2})^2 +(g_{i-1}-g_i)^2,\\
\beta_3=&\frac{21520059541}{19838649600}(g_{i-2} - 5g_{i-1} + 10g_i - 10g_{i+1} + 5g_{i+2} - g_{i+3})^2\\
&+\frac{1}{440858880}(1851g_{i-2} - 31123g_{i-1} + 84114g_i - 84114g_{i+1} + 31123g_{i+2} - 1851g_{i+3})^2\\
&+\frac{1}{2246400}(131g_{i-2} - 1173g_{i-1} + 1042g_i + 1042g_{i+1} - 1173g_{i+2} + 131g_{i+3})^2\\
&+\frac{1421461}{5241600}(g_{i-2} - 3g_{i-1} + 2g_i + 2g_{i+1} - 3g_{i+2} + g_{i+3})^2 + (-g_i+g_{i+1})^2,\\
\beta_4=&\frac{1}{326918592000}(459034864256g_{i-3}^2+21743036840504g_{i-2}^2+193082473956456g_{i-1}^2\\
&-633842107028865g_{i-1}g_i+533907688202000g_i^2+610844549719320g_{i-1}g_{i+1}\\
&-1054387388310025g_ig_{i+1}+533907688202000g_{i+1}^2-346465395978597g_{i-1}g_{i+2}\\
&+610844549719320g_ig_{i+2}-633842107028865g_{i+1}g_{i+2}+193082473956456g_{i+2}^2\\
&+107421495993504g_{i-1}g_{i+3}-192940877965535g_ig_{i+3}+204591754773560g_{i+1}g_{i+3}\\
&-127805625375939g_{i+2}g_{i+3}+21743036840504g_{i+3}^2-7g_{i-2}(18257946482277g_{i-1}\\
&-29227393539080g_i+27562982566505g_{i+1}-15345927999072g_{i+2}\\
&+4682081208019g_{i+3}-605170517992g_{i+4})+g_{i-3}(-6214446276409g_{i-2}\\
&+17764726801752g_{i-1}-27790531210295g_{i}+25709223617840g_{i+1}\\
&-14082592044087g_{i+2}+4236193625944g_{i+3}-540644243257g_{i+4})\\
&-14082592044087g_{i-1}g_{i+4}+25709223617840g_ig_{i+4}-27790531210295g_{i+1}g_{i+4}\\
&+17764726801752g_{i+2}g_{i+4}-6214446276409g_{i+3}g_{i+4} +459034864256g_{i+4}^2).
\end{split}
\end{equation*}

\section{The smoothness indicators of the multi-resolution WENO in Section \ref{Sect:ConvDisc}}
The expressions of the smoothness indicators $\beta_1, \beta_2, \beta_3$ and $\beta_4$ in the definition \eqref{eq:WENOweightsconv} are given as follows: 
\begin{equation*}
\begin{split}
\beta_1=&\frac{121}{300}(4u_{i-1}^2 - 13u_{i-1}u_{i} + 13u_{i}^2 + 5u_{i-1}u_{i+1} - 13u_{i}u_{i+1} + 4u_{i+1}^2),\\
\beta_2=&\frac{1}{67200000}(112756316u_{i-2}^2 + 1657473113u_{i-1}^2 + 3613771547u_{i}^2 - 4707412996u_{i}u_{i+1} \\
&+ 1657473113u_{i+1}^2 + u_{i-1}(-4707412996u_{i} + 2846027902u_{i+1} - 631012985u_{i+2}) \\
&- 37u_{i-2}(22231031u_{i-1} - 29557877u_i + 17054405u_{i+1} - 3632623u_{i+2})\\
&+ 1093641449u_iu_{i+2} - 822548147u_{i+1}u_{i+2} + 112756316u_{i+2}^2),\\
\beta_3=&\frac{1}{696729600000000}(1191368301143900u_{i-3}^2 + 40597776375544695u_{i-2}^2 + 247305166240620450u_{i-1}^2 \\
&- 645469279961828850u_{i-1}u_{i} + 435991844146445900u_{i}^2 + 462392655194742375u_{i-1}u_{i+1}\\
&- 645469279961828850u_{i}u_{i+1} + 247305166240620450u_{i+1}^2 - 173309501486101245u_{i-1}u_{i+2} \\
&+ 248679888594363540u_{i}u_{i+2} - 196807231970740065u_{i+1}u_{i+2} + 40597776375544695u_{i+2}^2 \\
&- 1111u_{i-3}(12262548688269u_{i-2} - 28726145190345u_{i-1} + 35285736074690u_{i}\\
&- 24003850977690u_{i+1} + 8594065910145u_{i+2} - 1267677095269u_{i+3})\\
&+ 26668278436213590u_{i-1}u_{i+3} - 39202452778980590u_{i}u_{i+3} + 31914747306473295u_{i+1}u_{i+3} \\
&- 13623691592666859u_{i+2}u_{i+3} + 1191368301143900u_{i+3}^2 - 3u_{i-2}(65602410656913355u_{i-1} \\
&- 82893296198121180u_{i} + 57769833828700415u_{i+1} - 21137663643408778u_{i+2} + 3182669075390365u_{i+3})),\\
\beta_4=&\frac{1}{66952927641600000000000}(118023030647523865761268u_{i-4}^2 + 7310547749226640675750391u_{i-3}^2 \\
&+ 87932799446502525538448131u_{i-2}^2 - 346770742441690854589241870u_{i-2}u_{i-1}\\
&+ 348437796897631753060606543u_{i-1}^2 + 420998749597335849842664020u_{i-2}u_{i} \\
&- 861675614228526685698528040u_{i-1}u_{i} + 542828069875679609578330025u_{i}^2\\
&- 322534893789595815857580562u_{i-2}u_{i+1} + 671504849402026017193551778u_{i-1}u_{i+1} \\
&- 861675614228526685698528040u_{i}u_{i+1} + 348437796897631753060606543*u_{i+1}^2 \\
& + 152497445486154126541234954u_{i-2}u_{i+2} - 322534893789595815857580562u_{i-1}u_{i+2} \\
&+ 420998749597335849842664020u_{i}u_{i+2} - 346770742441690854589241870u_{i+1}u_{i+2} \\
&+ 87932799446502525538448131u_{i+2}^2 - 40748949157768827653691322u_{i-2}u_{i+3} \\
&+ 87450773089630621117847210u_{i-1}u_{i+3} - 115964286531853100210189720u_iu_{i+3} \\
&+ 97177290555967138846616342u_{i+1}u_{i+3} - 50211686059457479298460454u_{i+2}u_{i+3} \\
&+ 7310547749226640675750391u_{i+3}^2 + u_{i-3}(-50211686059457479298460454u_{i-2} \\
&+ 97177290555967138846616342u_{i-1} - 115964286531853100210189720u_{i} \\
&+ 87450773089630621117847210u_{i+1} - 40748949157768827653691322u_{i+2} \\
&+ 10740612518336469754371154u_{i+3} - 1227581961831224200169527u_{i+4}) \\
&- 11111u_{i-4}(165355769190611079815u_{i-3} - 556722292425786917777u_{i-2} \\
&+ 1058534363076305872403u_{i-1} - 1243189747760266986565u_{i} \\
&+ 923938536129339626101u_{i+1} - 424690494183694671475u_{i+2} \\
&+ 110483481399624174257u_{i+3} - 12465257465998279583u_{i+4}) \\
&+ 4718736080875031494758725u_{i-2}u_{i+4} - 10265881074933092585608211u_{i-1}u_{i+4} \\
&+ 13813081287364326487723715u_{i}u_{i+4} - 11761375308140834548269733u_{i+1}u_{i+4} \\
&+ 6185741391142918443420247u_{i+2}u_{i+4} - 1837267951476879707824465u_{i+3}u_{i+4} \\
&+ 118023030647523865761268u_{i+4}^2).
\end{split}
\end{equation*}

\end{document}